\documentclass[12pt]{article}
\usepackage{fullpage}
\usepackage{psfig}
\usepackage{latexsym}
\usepackage{amssymb}
\newtheorem{prop}{Proposition}[section]
\newtheorem{cor}{Corollary}[section]
\newtheorem{lemma}{Lemma}[section]

\newtheorem{definition}{Definition}[section]
\newcommand{\pfend}{{\large $\Box$}}

\newcommand{\no}{\nonumber}
\newcommand{\D}{\partial}
\newcommand{\be}{\begin{equation}}
\newcommand{\ee}{\end{equation}}
\newcommand{\ba}{\begin{eqnarray}}
\newcommand{\ea}{\end{eqnarray}}
\newcommand{\R}{\mathbb{R}}
\newcommand{\C}{\mathbb{C}}
\newcommand{\bI}{\mathbb{I}}
\newcommand{\comment}[1]{}
\newcommand{\op}[1]{{\cal #1}}
\newcommand{\vc}[1]{{\vec{#1}}}
\newcommand{\mat}[1]{{\bf #1}}
\begin{document}
\title{Metastability of Breather Modes of Time Dependent Potentials}
\author{P.~D.~Miller\thanks{Address: Department of Mathematics and
Statistics, P. O. Box 28M, Monash University, VIC 3800 AUSTRALIA.  Email:
{\tt millerpd@mail.maths.monash.edu.au}}\\ {\small Institute for
Advanced Study, Princeton}\\{\small Monash University} \and
A.~Soffer\thanks{Address: Department of Mathematics, Hill Center,
Rutgers University, 110 Frelinghuysen Road, Piscataway, NJ 08854-8019.
Email: {\tt soffer@math.rutgers.edu} }\\{\small Rutgers University}
\and M.~I.~Weinstein\thanks{Address: Mathematical Sciences Research,
Bell Labs 2C-358, Lucent Technologies, 600 Mountain Avenue, Murray
Hill, NJ 07974.  Email: {\tt miw@research.bell-labs.com} } \\ {\small
Bell Labs}\\{\small University of Michigan} }

\date{20 August 1999\\ Submitted to {\em Nonlinearity}\\
Revised 10 February 2000}

\maketitle

\begin{abstract}
We study the solutions of linear Schr\"odinger equations in which the
potential energy is a periodic function of time and is sufficiently
localized in space.  We consider the potential to be close to one that
is time periodic and yet explicitly solvable. A large family of such
potentials has been constructed and the corresponding Schr\"odinger
equation solved by Miller and Akhmediev \cite{physd}.  Exact bound
states, or breather modes, exist in the unperturbed problem and are
found to be generically metastable in the presence of small periodic
perturbations.  Thus, these states are long-lived but eventually
decay. On a time scale of order $\epsilon^{-2}$, where $\epsilon$ is a
measure of the perturbation size, the decay is exponential, with a
rate of decay given by an analogue of Fermi's golden rule.  For times
of order $\epsilon^{-1}$ the breather modes are frequency shifted.
This behavior is derived first by classical multiple scale expansions,
and then in certain circumstances we are able to apply the rigorous
theory developed by Soffer and Weinstein \cite{ionization} and
extended by Kirr and Weinstein \cite{KW} to justify the expansions and
also provide longer time asymptotics that indicate eventual dispersive
decay of the bound states with behavior that is algebraic in time.  As
an application, we use our techniques to study the frequency
dependence of the guidance properties of certain optical waveguides.
We supplement our results with numerical experiments.
\end{abstract}

\noindent{\bf AMS Classification Scheme:  37K55, 35C20, 35C15, 35Q55, 81Q05, 81Q15}

\renewcommand{\theequation}{\arabic{section}.\arabic{equation}}

\section{Introduction and Overview}
\label{sec:intro}
\setcounter{equation}{0}
We are interested in the initial value problem for the linear
Schr\"odinger equation in one space dimension
\be i\D_t f = \left(\ -\frac{1}{2}\D_x^2 +
V(x,t)\ \right)f\ \equiv\ \op{H}(t)\ f\
\label{eq:generalschrod}
\ee 
Here, $f$ is a complex-valued function of $x\in\R$ and
$t\in\R$.  We assume that $V(x,t)$ is a smooth real-valued
potential energy function which is sufficiently localized in space
(for example, of Schwartz class). In our specific applications,
$V(x,t)$ will be taken to be a periodic function of
$t$. However, the techniques we use can be adapted for to more general
time-dependence \cite{KW}.  Note that (\ref{eq:generalschrod}) is
a nonautonomous Hamiltonian system: 
\be i\D_t f\ =\ {\delta 
h\over\delta f^*}[f,f^*,t],\no 
\ee 
where $h$ denotes the
Hamiltonian energy: 
\be h[f,f^*,t]\ =\ \int \left({1\over2}|\partial_x
f|^2\ +\ V(x,t)|f|^2\right)\ \ dx.\no
\ee 
If $V(x,t)$ is not
independent of $t$, $h$ is not a conserved integral of the
flow. On the other hand, since the potential $V$ is real-valued, the
flow defined by (\ref{eq:generalschrod}) is always unitary in $L^2(\R)$, {\em
i.e.}  
\be \int |f(x,t)|^2\ dx\ =\ \int |f(x,0)|^2\
dx, \ \ t\in\R .  \label{eq:unitaryflow} 
\ee

In applications, it is often natural to decompose $V(x,t)$ as: \ba
V(x,t)\ =\ V_0(x,t)\ +\ W(x,t),\no \ea where $V_0(x,t)$ denotes an
{\em unperturbed potential}, and $W(x,t)\doteq
V(x,t)-V_0(x,t)$, denotes a small perturbation. Thus, \ba
\op{H}(t)\ \equiv\ \op{H}_0(t)\ +\ \op{W}(t),\label{eq:Hdecomp} \ea
and (\ref{eq:generalschrod}) can be rewritten as: \be i\D_t f\ =\
\left(\ \op{H}_0(t) + \op{W}(t)\ \right) f \label{eq:schrod1}\ee Here,
we have denoted the multiplication operator $f\mapsto W(x,t)f$ by
$\op{W}(t)$.  The choice of $V_0(x,t)$ is often dictated by some
{\em a priori} knowledge of the solutions of the unperturbed system
\ba i\D_t f\ =\ \op{H}_0(t) f\label{eq:unperturbedschrod}.  \ea

A problem of importance is then to contrast the detailed dynamics of
solutions to (\ref{eq:schrod1}) with those of the unperturbed system
(\ref{eq:unperturbedschrod}).  In particular, {\em if
(\ref{eq:unperturbedschrod}) has bound state solutions (breather
modes, or solutions having finite energy and not decaying as
$|t|\to\infty$) do they persist in the perturbed dynamical system
(\ref{eq:schrod1})?}

The simplest variant of this problem is the case where the unperturbed
part is stationary, {\em i.e.}  $V_0(x,t)\ =\ V_0(x)$. Suppose the
operator $\op{H}_0$ has an $L^2$ eigenfunction.  The unitary evolution
of the spatially localized eigenfunction is time-periodic and
represents a bound state solution of the unperturbed Schr\"odinger
equation (\ref{eq:unperturbedschrod}).  The perturbed model (in this
and in the more general case when $\op{H}_0$ has multiple discrete
eigenvalues) is related to the problem of ionization of an atom by a
time-dependent electromagnetic field
\cite{LandauLifshitz,GalindoPascal} and the problem of describing the
effects of weak inhomogeneities on the propagation of continuous waves
in optical fibers \cite{M74}.  Using a time-dependent method developed in the
context of (i) quantum resonances and the perturbation theory of
embedded eigenvalues in the continuous spectrum \cite{resonances} and (ii)
resonances and radiation damping of bound nonlinear wave equations
\cite{rdamping}, Soffer and Weinstein studied the metastability of
such states \cite{ionization}.  Generalizations of this theory for
handling multifrequency perturbations \cite{KW} and the interference
of multiple bound states in the unperturbed problem \cite{KWmultimode}
have been explored by Kirr and Weinstein.  Based on the observation
that the mechanism for instability of the bound state is coupling of
the bound state to the continuous spectral modes, the analysis was
carried out at the level of the coupled equations for the bound state
and dispersive components of the solution.  Under general hypotheses
on unperturbed Hamiltonian (local energy decay estimates on the
unitary propagator $e^{-it\op{H}_0}$) this equivalent dynamical system
was studied and it was shown that a bound state is generically
unstable but long-lived. The lifetime is given by a formula analogous
to the Fermi golden rule \cite{ionization}.

In this paper we consider the case where the unperturbed Hamiltonian
is genuinely time-dependent. A physical application of the theory we
develop, in the context of frequency detuning in periodically
modulated optical waveguides \cite{BesleyOL,SAM}, will be presented in 
\S \ref{sec:application}.  Let $\op{H}_0(t)\ =\ -{1\over2}\D_x^2\ +\
V_0(x,t)$, where $V_0(x,t)$ is smooth, periodic in $t$ with the same
period for each $x$ and of sufficiently rapid decay for large $x$
for each $t\in\R$. The particular choices of $V_0(x,t)$ we
consider in this paper belong to a large family of very special,
so-called {\em separable}, time-dependent potentials, $V_0(x,t),\
x\in\R$, studied by Miller and Akhmediev \cite{physd}.  The separable
potential $V_0(x,t)$ can be chosen to be time-periodic, in which case
the unperturbed problem supports exact bound states (breather modes)
and the initial value problem for (\ref{eq:unperturbedschrod}) can be
solved exactly. That is, a complete set of eigenmodes and generalized
eigenmodes can be explicitly displayed with respect to which the
dynamics of (\ref{eq:unperturbedschrod}) is diagonal.  This class of
potentials is intimately connected with the soliton theory of
completely integrable multicomponent cubic nonlinear Schr\"odinger
equations \cite{Manakov,FT}.

The existence of such exact breather modes in the unperturbed
time-periodic problem is quite remarkable and we believe that this is
a highly non-generic phenomenon\footnote{The scarcity of breather
solutions of {\em nonlinear} wave equations defined on a spatial
continuum of infinite extent has been extensively explored in the
setting of perturbations of the completely integrable sine-Gordon
equation; see, for example, \cite{BMW,Bnr,D,Kich}.  The connection
with linear nonautonomous problems can be made by viewing a breather
solution of a nonlinear dynamical problem as a bound state of a linear
problem with a given (self-consistent) potential.}.  Indeed, from a
general dynamical systems perspective, (\ref{eq:unperturbedschrod})
with such a choice of $V_0(x,t)$, may be viewed as a parametrically
forced wave equation (here we are actually considering the
time-periodic function $V_0(x,t)$ itself to be the sum of a
time-independent part and a time-periodic modulation).  One therefore
expects that the presence of resonances will perturb the Floquet
multipliers (corresponding to bound states) off the unit circle as in
the elementary example of Mathieu's equation \cite{Arnold}.  The
persistence of breather solutions under the time-periodic perturbation
would imply the non-departure of a Floquet multiplier from the unit
circle {\em to all orders} in the size of the perturbation. The fact
that infinitely many such conditions hold for these special separable
potentials is no doubt linked to the infinite sequence of symmetries
and time invariants enjoyed by the completely integrable nonlinear
flow that underpins the construction of the separable potentials (see
\ref{app:separable} for more details).  Of course, this is only a
heuristic picture. In fact, the perturbation theory of the Floquet
multipliers is complicated by the fact that they are embedded in the
continuous spectrum which covers the unit circle. However, spectral
deformation methods have been developed for some classes of models
that could well be adapted here.  Relevant technical details can be
found in \cite{CyconEtAl,HislopSigal}\footnote{ For nonlinear wave
equations defined on an infinite lattice ({\em e.g.} discrete
sine-Gordon, discrete $\phi^4$), breather solutions can be constructed
for sufficiently large lattice spacing; see for example
\cite{MacAu94}.  The radiative decay of such discrete breathers, for
sufficiently small lattice spacing, is expected to be governed by a
mechanism of the kind studied in this paper; see also \cite{rdamping}.
Related to this are results for the dynamics of kinks of discrete
nonlinear wave equations, in which the techniques of this paper have
been used to study the ``pinning'' of discrete kinks on a lattice
site.  This pinning is marked by the slow radiative decay of spatially
localized and time-periodic or quasiperiodic oscillations about a
static kink \cite{KevrWe00}.  }.

We want to make our motivation for pursuing deformations of these
admittedly rather special periodic potentials very clear.  First of
all, the problem is relevant to the analysis of optical waveguides.  In
the paraxial approximation, the slowly varying envelope of a highly
oscillatory electric field in a dielectric medium with inhomogeneous
dielectric properties (index of refraction) satisfies an
equation of form (\ref{eq:generalschrod}). Here, $t$ denotes the
longitudinal variable, the direction of propagation, and $x$ is the
transverse spatial variable.  For inhomogeneous index profiles
corresponding to exactly solvable potentials, light {\em of a
particular frequency} propagates as a non-attenuating bound state mode
in these wave guides.  However, if the light frequency deviates from
the ``integrable frequency'' the propagating wave will be governed by
the perturbed equation (\ref{eq:schrod1}). Thus the question of
whether such modes persist and if not what their lifetime is for the
perturbed dynamics naturally arises.  We will give more details about
this problem in \S \ref{sec:application}.

But it is also true that the study of perturbed separable periodic
potentials is important in general terms.  Given an arbitrary
time-periodic potential in the Schr\"odinger equation, one wants to
study the corresponding dynamics using perturbation theory.  In doing
so, the first question one must address is that of finding a
``nearby'' problem that can be solved exactly.  We simply take the
point of view that many periodic potentials will be closer to a
separable periodic potential (in a sense that can be made precise)
than to any time-independent potential.  

In any case, with the explicit spectral theory associated with
$V_0(x,t)$ in hand, our goal is to carry out a detailed analytical
study of the coupled mode dynamics induced by a time-dependent
perturbation $W(x,t)$.  We establish the generic metastabilty of the
exact bound states associated with separable periodic potentials
$V_0(x,t)$ and obtain a detailed picture of the dynamics.

The paper is structured as follows. In \S \ref{sec:separable} we
first review the construction of time-dependent exactly solvable
potentials \cite{physd} from a set of discrete data, and then show how
the initial value problem for such {\em separable} potentials can be
solved explicitly.  We then describe how properties of the separable
potentials depend on the choice of the discrete data generating them.
Next, we derive by projection onto an orthonormal basis the general
coupled mode equations which arise when a separable potential is
perturbed by some arbitrary correction $W(x,t)$.  This section
will then conclude with a detailed derivation of the {\em
two-soliton} time-periodic potential and its associated explicit
complete set of bound states and generalized eigenfunctions. More
details about the separable potentials described in 
\S \ref{sec:separable} are given in \ref{app:separable}.  

In \S \ref{sec:periodic} Floquet theory is then used to map the
coupled mode equations to a system associated with a time-dependent
perturbation of an {\em autonomous} system, a situation analyzed in
detail in \cite{ionization} and \cite{KW}.  In \S \ref{sec:analysis},
we then describe the dynamics of solutions of the coupled mode
equations for even time-periodic perturbations $W(x,t)$ of a separable
two-soliton even time-periodic potential $V_0(x,t)$ ($V_0(x,t)$ and
$W(x,t)$ both share even parity in $x$ and have the same temporal
period).  In particular, we study the initial value problem when the
initial condition is a pure bound state of the unperturbed problem.
First, we study the small time behavior of the coupled mode equations
(without requiring $W(x,t)$ to be small) and deduce that the bound
state amplitude behaves as $1-Ct^2$ for some constant $C$ and
interpret this result in the context of the theory of ideal
measurements in quantum mechanics (the ``watched pot'' effect).  We
then assume the perturbation $W(x,t)$ to be small, of size $\epsilon$,
and seek the behavior of the bound state amplitude over intermediate
times of order $\epsilon^{-1}$ and $\epsilon^{-2}$ using the classical
method of multiple scales.

We show the existence of a perturbation-induced frequency shift of the
breather mode evident on timescales of order $\epsilon^{-1}$ and
exponential decay of the bound state mode amplitude on timescales of
order $\epsilon^{-2}$.  The condition for the decay constant to be
nonzero is a direct analog of the ``Fermi golden rule''.  

Then, using the transformation to an autonomous system found in \S
\ref{sec:periodic}, we show how the rigorous theory developed for
multifrequency perturbations of autonomous systems by Kirr and
Weinstein \cite{KW} can be applied in some cases to justify the
multiple scales calculation, and to provide more detailed information
about the infinite time behavior of the solution.  This analysis
completes the portrait of the dynamics, showing that the exponential
decay is ultimately washed out in a sea of dispersive waves, at which
point the decay becomes algebraic in time.

Having described the theory, in \S \ref{sec:application} we consider
an application of the analysis to a problem of frequency detuning in
planar optical waveguides.  Finally, in \S \ref{sec:numerics} the
prediction of an exponential decay constant $\Gamma$ for the bound
state mode amplitude found in section \ref{sec:analysis} is compared
to numerical simulations of the perturbed time-dependent Schr\"odinger
equation.

A detailed description of the theory of separable potentials, at once
summarizing for completeness and also further developing the results
of \cite{physd}, can be found in \ref{app:separable}.  In
\ref{app:localdecay} the reader will find the proofs of the decay
estimates that we will use in \S \ref{sec:analysis} in order to
apply the results of \cite{KW}.

\noindent{\bf Regarding notation:}
We will use the inner product
\begin{equation}
\langle f(\cdot),g(\cdot)\rangle = \int_{-\infty}^\infty f(x)^*g(x)\,dx
\end{equation}
on $L^2(\R)$.  Occasionally, the angled brackets will denote the inner
product in more general Hilbert spaces.  Linear operators will be denoted
with calligraphic letters, vectors with arrows, and matrices
with boldface letters.  We will often use the function defined by:
\begin{equation}
\langle x\rangle\doteq (1+x^2)^{1/2}\,.
\end{equation}
Complex conjugation will be denoted with stars, and time averages will
be denoted with bars.

\section{Exactly solvable time-dependent potentials}
\label{sec:separable}
\setcounter{equation}{0}
In this section we recall for our purposes a class of time-dependent
potentials $V_0(x,t)$ related to $M$-soliton solutions of certain
completely integrable nonlinear flows. Because of the intimate
connection of these potentials to integrable systems, it is possible
to explicitly derive the spectral representation associated with such
potentials \cite{physd}. This section is divided into five
parts. First, the direct construction of separable potentials from a
set of discrete data $\cal D$ is outlined.  Then, we show how the same
discrete data $\cal D$ gives rise to formulas for a complete set of
modes for the time-dependent Schr\"odinger equation corresponding to
the separable potential $V_0(x,t)$ and how this basis is easily used
to express the general solution of the initial value problem.  We then
give a qualitative description of the kinds of functions $V_0(x,t)$
one can obtain from this procedure.  As we ultimately want to consider
perturbations of $V_0(x,t)$, we next show how to use the basis of
solutions to the unperturbed problem to derive the coupled mode
equations which trivialize the unperturbed dynamics and lay bare the
perturbative effects.  Finally, we specialize to the case of an even
periodic potential corresponding to a two-soliton solution of the cubic
nonlinear Schr\"odinger equation.  As one might anticipate, the
evenness (in $x$) of the potential leads to some simplifications in
the spectral representation.

\subsection{Separable time-dependent potentials.}
\label{sec:sub-separable}
Let us present the construction of the family of time-dependent
potentials that we will consider in this paper, and describe their
properties with respect to the linear Schr\"odinger equation.  More details
can be found in \ref{app:separable}.

Each potential we shall consider will be specified by a certain set of
discrete data.  Let $N$ and $M$ be independent natural numbers.  A set
of discrete data $\cal D$ consists of $M$ distinct complex numbers
$\lambda_1,\dots, \lambda_M$ in the upper half-plane, and $M$ vectors
$\vc{g}^{(1)},\dots,\vc{g}^{(M)}$ in $\C^N$.

The discrete data $\cal D$ is used to build a potential function $V_0(x,t)$
in the following way.  Introduce the scalar expression
\begin{equation}
a(x,t,\lambda)=\left(\lambda^M + \sum_{p=0}^{M-1}\lambda^pa^{(p)}(x,t)\right)
e^{-2i(\lambda x + \lambda^2 t)}\,,
\end{equation}
and the $N$-component vector expression
\begin{equation}
\vc{b}(x,t,\lambda)=\sum_{p=0}^{M-1}\lambda^p\vc{b}^{(p)}(x,t)\,.
\end{equation}
In these expressions, the coefficients $a^{(p)}(x,t)$ and
$\vc{b}^{(p)}(x,t)$ are undetermined functions of $x$ and $t$.  They
will now be determined with the use of the discrete data $\cal D$.
For $k=1,\dots,M$, we insist that $a(x,t,\lambda)$ and
$\vc{b}(x,t,\lambda)$ satisfy the relations:
\begin{equation}
\begin{array}{rcl}
a(x,t,\lambda_k) &=& \vc{g}^{(k)\dagger} \vc{b}(x,t,\lambda_k)\,,\\
\vc{b}(x,t,\lambda_k^*) &=& -a(x,t,\lambda_k^*)\vc{g}^{(k)}\,.
\end{array}
\label{eq:linsys}
\end{equation}
These equations amount to a square linear inhomogeneous system of
algebraic equations for the coefficient functions $a^{(p)}(x,t)$ and
the components of $\vc{b}^{(p)}(x,t)$.  We will soon illustrate this
procedure with a concrete example.  From the solution of this linear
system, the potential function connected with the discrete data $\cal
D$ is given in terms of the components of $\vc{b}^{(M-1)}(x,t)$ by
\begin{equation}
V_0(x,t)\doteq -4\sum_{n=1}^N \Bigg| b_n^{(M-1)}(x,t)\Bigg|^2\,.
\label{eq:potentialdef}
\end{equation}
This function $V_0(x,t)$ is a genuinely time-dependent potential well.
Furthermore, it can be shown that $V_0(x,t)$ is in the Schwartz space
as a function of $x$, and its $L^1$ norm is constant in $t$.  The
latter follows from the fact that $V_0(x,t)$ can be viewed as the
self-consistent nonlinear potential for an $N$-component cubic
nonlinear Schr\"odinger equation, which conserves the sum of the $L^2$
norms of the $N$ field components, which are proportional to the
$b_n^{(M-1)}(x,t)$ for $n=1,\dots,N$.

\subsection{Solution of the linear Schr\"odinger equation with a separable
potential.}  Along with the potential function $V_0(x,t)$,
this construction starting from the discrete data $\cal D$ also
provides all of the solutions of the corresponding linear
Schr\"odinger equation \cite{physd}.  These are built from the function
$a(x,t,\lambda)$ as follows.  For all real $\lambda$, set
\begin{equation}
\Psi_{\rm d}(x,t,\lambda)\doteq\left(\pi\prod_{k=1}^M|\lambda-\lambda_k|^2\right)^{-1/2}a(x,t,\lambda)\,,
\end{equation}
and let the functions $\Psi_{{\rm b},1}(x,t),\dots,\Psi_{{\rm
b},M}(x,t)$ be the result of applying the Gram-Schmidt procedure (in
$L^2(\R)$ with respect to $x$) to the functions
$a(x,t,\lambda_1^*),\dots,a(x,t,\lambda_M^*)$ at any fixed value of
$t$.  Then we have \cite{physd}: 
\begin{enumerate}
\item
Each function $\Psi_{\rm d}(x,t,\lambda)$ for $\lambda\in\R$ and each
function $\Psi_{{\rm b},k}(x,t)$ is a solution of the linear
Schr\"odinger equation with potential $V_0(x,t)$.  The fact that the
$L^2$ inner product is an invariant of the evolution shows that the
functions $\Psi_{{\rm b},k}(x,t)$ do not depend on the choice of the
time $t$ at which they are obtained from the Gram-Schmidt process.
\item
For any fixed $t$, these functions form an orthonormal basis of $L^2(\R)$.
\end{enumerate}
These facts show us how to solve the initial value problem for
the linear Schr\"odinger equation for the potential $V_0(x,t)$.
Namely, to find the solution of
\begin{equation}
i\partial_t f + \frac{1}{2}\partial_x^2 f -V_0(x,t)f=0\,,\hspace{0.3 in}
f(x,0)=f_0(x)\in L^2(\R)\,,
\end{equation}
one simply projects the initial data onto the basis at $t=0$ by defining:
\begin{equation}
\hat{f}(\lambda)\doteq \langle \Psi_{\rm d}(\cdot,0,\lambda),f_0(\cdot)\rangle
\,,\hspace{0.3 in}
\hat{f}_k\doteq\langle\Psi_{{\rm b},k}(\cdot,0),f_0(\cdot)\rangle\,,
\end{equation}
and then recovers the solution as a superposition of modes:
\begin{equation}
f(x,t)=\sum_{k=1}^M\hat{f}_k\Psi_{{\rm b},k}(x,t) +\int_{-\infty}^\infty
\hat{f}(\lambda)\Psi_{\rm d}(x,t,\lambda)\,d\lambda\,.
\end{equation}
If the potential $V_0(x,t)$ is slightly perturbed, it may still be convenient 
to expand in this basis, but then the coefficients $\hat{f}(\lambda)$ and
$\hat{f}_1,\dots,\hat{f}_M$ will become time-dependent.

\subsection{Qualitative description of separable potentials.}
Let us describe the types of potential functions $V_0(x,t)$ that can
be obtained by this procedure.  In the generic case when the real
parts of the parameters $\lambda_1,\dots,\lambda_M$ are all distinct,
these potentials have the form of a collision among $M$ moving
potential wells.  That is, as $t\rightarrow \pm\infty$,
\begin{equation}
V_0(x,t)\sim\sum_{k=1}^M V_0^{(k)\pm}(x,t)\,,
\end{equation}
where the individual wells have the form
\begin{equation}
V_0^{(k)\pm}(x,t)=-4\rho_k^2{\rm sech}^2 (2\rho_k(x+2\sigma_k t)-\delta_k^\pm)
\,,
\end{equation}
where $\delta_k^\pm$ are constants that depend on the vectors
$\vc{g}^{(1)}, \dots,\vc{g}^{(M)}$, and where
$\lambda_k=\sigma_k+i\rho_k$.  Considered in isolation, each well
carries a single bound state.  When the wells collide for finite $t$,
the bound states can become mixed, and a state $f(x,t)$ that is bound
in a single well as $t\rightarrow -\infty$ will have a component in
each well as $t\rightarrow +\infty$.  The associated scattering matrix
can be computed exactly \cite{multiport}.

If some of the parameters $\lambda_k$ share the same real part
$\sigma$, then the asymptotics of the potential $V_0(x,t)$ in the
frame moving with velocity $-2\sigma$ will no longer be stationary,
but will be generally quasiperiodic.  In particular, if all of the
parameters $\lambda_k$ are purely imaginary, then the potential
$V_0(x,t)$ will generally be a quasiperiodic function of the time $t$.  This is
clear because taking $\lambda_k=i\rho_k$ with $\rho_k$ real and positive
ensures that the only time dependence that enters into the computation
of $V_0(x,t)$ is via the exponentials $\exp(\pm 2i\rho_k^2 t)$.  Such a
potential is automatically quasiperiodic.  We can further ensure that
the potential function is strictly periodic by making the frequencies
commensurate.  This will be true\footnote{For $M>2$.  The potential is
always periodic if $M=2$ and is stationary if $M=1$.} if the
parameters $\rho_k$ have the form
\begin{equation}
\rho_k=\sqrt{n_k\frac{\Omega_0}{2}+\Delta}\,,
\label{eq:comensurate}
\end{equation}
where $\Omega_0$ is some fundamental frequency and $n_k$ are distinct
integers.  This choice ensures that the frequencies 
\begin{equation}
\omega_{jk}\doteq 2\rho_j^2-2\rho_k^2=(n_j-n_k)\Omega_0
\end{equation}
are all integer multiples of $\Omega_0$.  Only the frequency
differences $\omega_{jk}$ are important because the potential is given
as a sum of absolute values (\ref{eq:potentialdef}).  In fact, it can
be seen from the form of the linear system (\ref{eq:linsys})
that
\begin{equation}
b_n^{(p)}(x,t)=e^{2i(\rho_1^2 + \dots +\rho_M^2)t}
G_{n,p}(\{e^{i\omega_{jk}t}\},x)\,,
\end{equation}
where $G_{n,p}$ is, for each fixed $x$, a rational function of the
exponentials $\exp(i\omega_{jk}t)$.  The sufficiency of the relations
(\ref{eq:comensurate}) to guarantee time periodicity of $V(x,t)$ with
fundamental frequency $\Omega_0$ is then clear from
(\ref{eq:potentialdef}).

%The integers $n_k$ can be chosen arbitrarily as long as they are distinct.
%One example is to take
%\begin{equation}
%n_k=\frac{k\cdot (k-1)}{2}\,,\hspace{0.3 in}\Delta=\frac{\Omega_0}{16}\,,
%\end{equation}
%which arranges that the $\rho_k$ are positive half-integer multiples of
%$\sqrt{\Omega_0}/2$:
%\begin{equation}
%\rho_k=\left(k-\frac{1}{2}\right)\frac{\sqrt{\Omega_0}}{2}\,,\hspace{0.1 in}
%k=1,2,3,\dots\hspace{0.2 in}\Rightarrow\hspace{0.2 in}
%\omega_{jk}\doteq 2\rho_j^2-2\rho_k^2=
%\frac{(j-k)(j+k-1)}{2}\Omega_0\,.
%\end{equation}

\subsection{Perturbed separable potentials and coupled mode equations.}
As we have already suggested, the explicit basis of exact solutions
derived in the previous subsection forms a natural coordinate system
in which to study perturbed problems.  Let $W(x,t)$ be a correction to
the potential energy, so that the equation becomes
\begin{equation}
if_t\ =\  \left(\ - \frac{1}{2}\D^2_x\ +\ V_0(x,t)\ \right) f\  +\  W(x,t) f\ 
 =\ \op{H}_0(t) f\ +\ \op{W}(t) f\,.
\label{eq:perturb}
\end{equation}
Here, $V_0(x,t)$ is a separable time-dependent potential built from
the discrete data ${\cal
D}=\{\lambda_1,\dots,\lambda_M,\vc{g}^{(1)},\dots,\vc{g}^{(M)}\}$.
So we use completeness to express $f(x,t)$ for each fixed $t$ in terms
of the basis of solutions of the unperturbed problem:
\begin{equation}
f(x,t)=\sum_{k=1} B_{{\rm b},k}(t)\Psi_{{\rm b},k}(x,t) +
\int_{-\infty}^\infty B_{\rm
d}(t,\lambda)\Psi_{\rm d}(x,t,\lambda)\,d\lambda\,.
\end{equation}
In the absence the the perturbation $W(x,t)$, the mode amplitudes
$B_{{\rm b},k}$ and $B_{\rm d}(t,\lambda),\ \lambda\in\R$ are governed
by the equations $\D_t B_{{\rm b},k}=0,\ \D_t B_{\rm d}(t,\lambda)=0$.
In the presence of a perturbation $W(x,t)$ {\em coupled mode
equations} can be derived by projecting (\ref{eq:perturb}) onto the
basis elements $\Psi_{{\rm b},k}(x,t)$ and $\Psi_{\rm
d}(x,t,\lambda)$.  This yields the system of {\em coupled mode
equations}:
\begin{equation}
\begin{array}{rcl}
i\D_t \vc{B}_{\rm b}(t) &=& \displaystyle
\mat{M}(t)\vc{B}_{\rm b}(t) + \int_{-\infty}^\infty B_{\rm d}(t,\lambda)
\vc{N}(t,\lambda)\,d\lambda\,,\\\\
i\D_t B_{\rm d}(t,\eta) &=&\displaystyle
\vc{N}(t,\eta)^\dagger \vc{B}_{\rm b}(t) +
\int_{-\infty}^\infty K(t,\eta,\lambda)B_{\rm d}(t,\lambda)\,d\lambda\,,
\end{array}
\label{eq:coupledmode}
\end{equation}
for the coefficients of $f(x,t)$, where $\vc{B}_{\rm b}(t)$ is the vector
of bound state amplitudes $B_{{\rm b},k}(t)$, and where the {\em matrix
elements} of the perturbation $W(x,t)$ are explicitly given by
\begin{equation}
\begin{array}{rcl}
M_{kj}(t) &=& \langle \Psi_{{\rm b},k}(\cdot,t),
W(\cdot,t)\Psi_{{\rm b},j}(\cdot,t) \rangle\,,\\\\ N_k(t,\lambda) &=&
\langle \Psi_{{\rm b},k}(\cdot,t), W(\cdot,t) \Psi_{\rm
d}(\cdot,t,\lambda)\rangle\,,\\\\ K(t,\eta,\lambda) &=& \langle
\Psi_{\rm d}(\cdot,t,\eta), W(\cdot,t) \Psi_{\rm
d}(\cdot,t,\lambda)\rangle\,,
\end{array}
\label{eq:matrixelts}
\end{equation}
where $N_k(t,\lambda)$ are the components of the vector
$\vc{N}(t,\lambda)$ and $M_{kj}(t)$ are the elements of the matrix
$\mat{M}(t)$.  In particular, it follows that the matrix $\mat{M}(t)$
is Hermitian and the scalar kernel $K(t,\eta,\lambda)$ is Hermitian
symmetric because $W(x,t)$ is real.  With the unperturbed problem
exactly diagonalized in this way, this system is a useful starting
point for perturbation theory.

\subsection{Even two-soliton periodic potentials.}
In this subsection, we illustrate the procedures described above in
some detail with an example that is important in applications and that
will guide the subsequent discussion.  We consider the case $N=1$ and
$M=2$, and accordingly introduce the expressions
\begin{equation}
a(x,t,\lambda)=(\lambda^2 + a^{(1)}(x,t)\lambda + a^{(0)}(x,t))e^{-2i(\lambda
x + \lambda^2 t)}\,,\hspace{0.3 in} 
b(x,t,\lambda)=b^{(1)}(x,t)\lambda + b^{(0)}(x,t)\,.
\end{equation}
Because $N=1$ these are both scalar expressions, and we have at the
moment four complex valued unknown functions, $a^{(0)}(x,t)$,
$a^{(1)}(x,t)$, $b^{(0)}(x,t)$ and $b^{(1)}(x,t)$.  To find these, we
introduce the discrete data $\lambda_1$, $\lambda_2$, $g^{(1)}$, and
$g^{(2)}$ (again, here the $g^{(k)}$ are complex scalars because
$N=1$).  The linear equations (\ref{eq:linsys}) then become
\begin{equation}
\begin{array}{rcl}
\displaystyle
(\lambda_1^2 + a^{(1)}(x,t)\lambda_1 + a^{(0)}(x,t))e^{-2i(\lambda_1x +
\lambda_1^2 t)} &=& g^{(1)*}(b^{(1)}(x,t)\lambda_1 + b^{(0)}(x,t))\\\\
\displaystyle
(\lambda_2^2 + a^{(1)}(x,t)\lambda_2 + a^{(0)}(x,t))e^{-2i(\lambda_2x +
\lambda_2^2 t)} &=& g^{(2)*}(b^{(1)}(x,t)\lambda_2 + b^{(0)}(x,t))
\end{array}
\end{equation}
and
\begin{equation}
\begin{array}{rcl}
b^{(1)}(x,t)\lambda_1^* + b^{(0)}(x,t) &=&
\displaystyle
-g^{(1)}(\lambda_1^{*2} +a^{(1)}(x,t)
\lambda_1^* + a^{(0)}(x,t))e^{-2i(\lambda_1^*x + \lambda_1^{*2}t)}\\\\
b^{(1)}(x,t)\lambda_2^* + b^{(0)}(x,t) &=&
\displaystyle
-g^{(2)}(\lambda_2^{*2} +a^{(1)}(x,t)
\lambda_2^* + a^{(0)}(x,t))e^{-2i(\lambda_2^*x + \lambda_2^{*2}t)}\,.
\end{array}
\end{equation}
Given the discrete data $\cal D$, one can solve these equations for
$a^{(0)}(x,t)$, $a^{(1)}(x,t)$, $b^{(0)}(x,t)$ and $b^{(1)}(x,t)$, say
by Cramer's rule, and thus obtain explicit expressions in terms of
exponential functions.

\begin{figure}[h]
\begin{center}
\mbox{\psfig{file=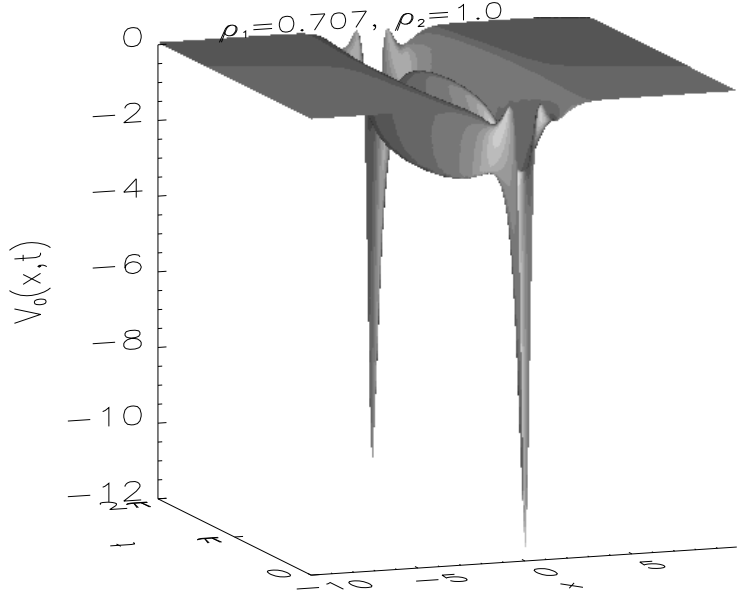,width=4.5 in}}
\end{center}
\caption{\em The potential well $V_0(x,t)$ for $\rho_1=1/\sqrt{2}$ and
$\rho_2=1$.  The phase parameters are
$\theta_1=\theta_2=0$.  This well is time-periodic with period $L=2\pi$
and even in $x$.}
\label{fig:R2potential}
\end{figure}

\begin{figure}
\begin{center}
\mbox{\psfig{file=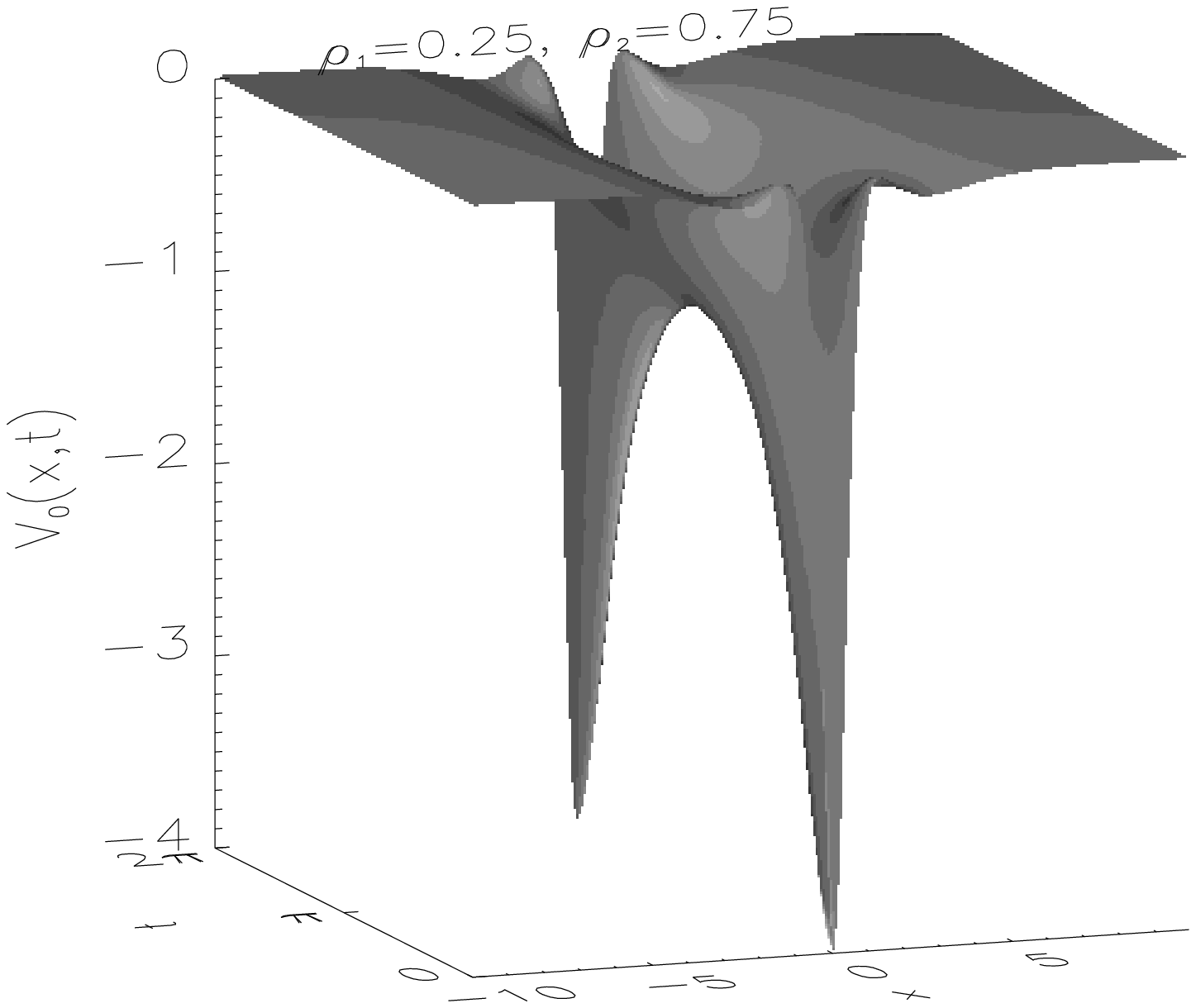,width=4.5 in}}
\end{center}
\caption{\em The potential well $V_0(x,t)$ for $\rho_1=1/4$ and
$\rho_2=3/4$.  The phase parameters are
$\theta_1=\theta_2=0$.  This well is time-periodic with period $L=2\pi$
and even in $x$.}
\label{fig:R9potential}
\end{figure}

Specializing to the case of $\lambda_1=i\rho_1$, $\lambda_2=i\rho_2$
(we assume without loss of generality that $\rho_2>\rho_1$), we obtain
a time-periodic potential function, since the parameters $\lambda_k$
are pure imaginary and then the commensurability condition is
automatically satisfied for $M=2$.  Furthermore choosing
$g^{(1)}=e^{i\theta_1}$ and $g^{(2)}=e^{i\theta_2}$ ensures that the
potential function is even in $x$.  Indeed, we then
find that with $s=\rho_2+\rho_1$ and $d=\rho_2-\rho_1$,
\begin{equation}
b^{(1)}(x,t)=2sd\frac{\displaystyle\rho_1\cosh(2\rho_2 x)e^{2i\rho_1^2 t + i\theta_1}-
\rho_2\cosh(2\rho_1 x)e^{2i\rho_2^2 t+i\theta_2}}{d^2\cosh(2sx) + s^2\cosh(2dx)-
4\rho_1\rho_2\cos(2sdt+\theta_2-\theta_1)}\,.
\end{equation}
The potential function is then given by
\begin{equation}
V_0(x,t)=-4|b^{(1)}(x,t)|^2\,,
\end{equation}
which is easily seen to be periodic in $t$ with period $L=\pi/(sd)$, and
an even function of $x$.
The shapes of these time-periodic potential wells are shown in
Figures~\ref{fig:R2potential} and \ref{fig:R9potential} for $\theta_1=\theta_2=0$ and two different choices of the parameters $\rho_1$ and $\rho_2$.  
From the solution of the same
linear system, we also find
\begin{equation}
a^{(0)}(x,t)=\rho_1\rho_2\frac{(\rho_1+\rho_2)^2S_1S_2 +\rho_1^2e^{-i\omega t-i(\theta_2-\theta_1)} + \rho_2^2e^{i\omega t+i(\theta_2-\theta_1)} - 2\rho_1\rho_2 C}
{2\rho_1\rho_2\cos(\omega t+(\theta_2-\theta_1)) - (\rho_1^2 + \rho_2^2)C +
(\rho_1+\rho_2)^2S_1S_2}\,,
\end{equation}
and
\begin{equation}
a^{(1)}(x,t)=i\frac{(\rho_1^2-\rho_2^2)\rho_1C_2S_1 + (\rho_2^2-\rho_1^2)\rho_2
C_1S_2}{2\rho_1\rho_2\cos(\omega t+(\theta_2-\theta_1)) - (\rho_1^2 + \rho_2^2)C +
(\rho_1+\rho_2)^2S_1S_2}\,,
\end{equation}
where $S_k\doteq \sinh(2\rho_k x)$, $C_k\doteq \cosh(2\rho_k x)$, and
$C\doteq \cosh(2(\rho_1+\rho_2)x)$, and
where the frequency is $\omega=2\pi/L=2sd$.  Note that $a^{(0)}(x,t)$ is an
even function of $x$, while $a^{(1)}(x,t)$ is odd.
We may then write the mode function $a(x,t,\lambda)$ in the form
\begin{equation}
\begin{array}{rcl}
a(x,t,\lambda)&=&\displaystyle
\left((\lambda^2 + a^{(0)}(x,t))\cosh(-2i\lambda x) +
\lambda a^{(1)}(x,t)\sinh(-2i\lambda x)\right)e^{-2i\lambda^2 t} +\\\\
&&\left((\lambda^2 + a^{(0)}(x,t))\sinh(-2i\lambda x) + \lambda a^{(1)}(x,t)
\cosh(-2i\lambda x)\right)e^{-2i\lambda^2 t}\,,
\end{array}
\end{equation}
in which the first term is even in $x$ and the second term is odd in $x$.
Also, it is clear that $a(-x,t,\lambda)=a(x,t,-\lambda)$.

\begin{figure}[h]
\begin{center}
\mbox{\psfig{file=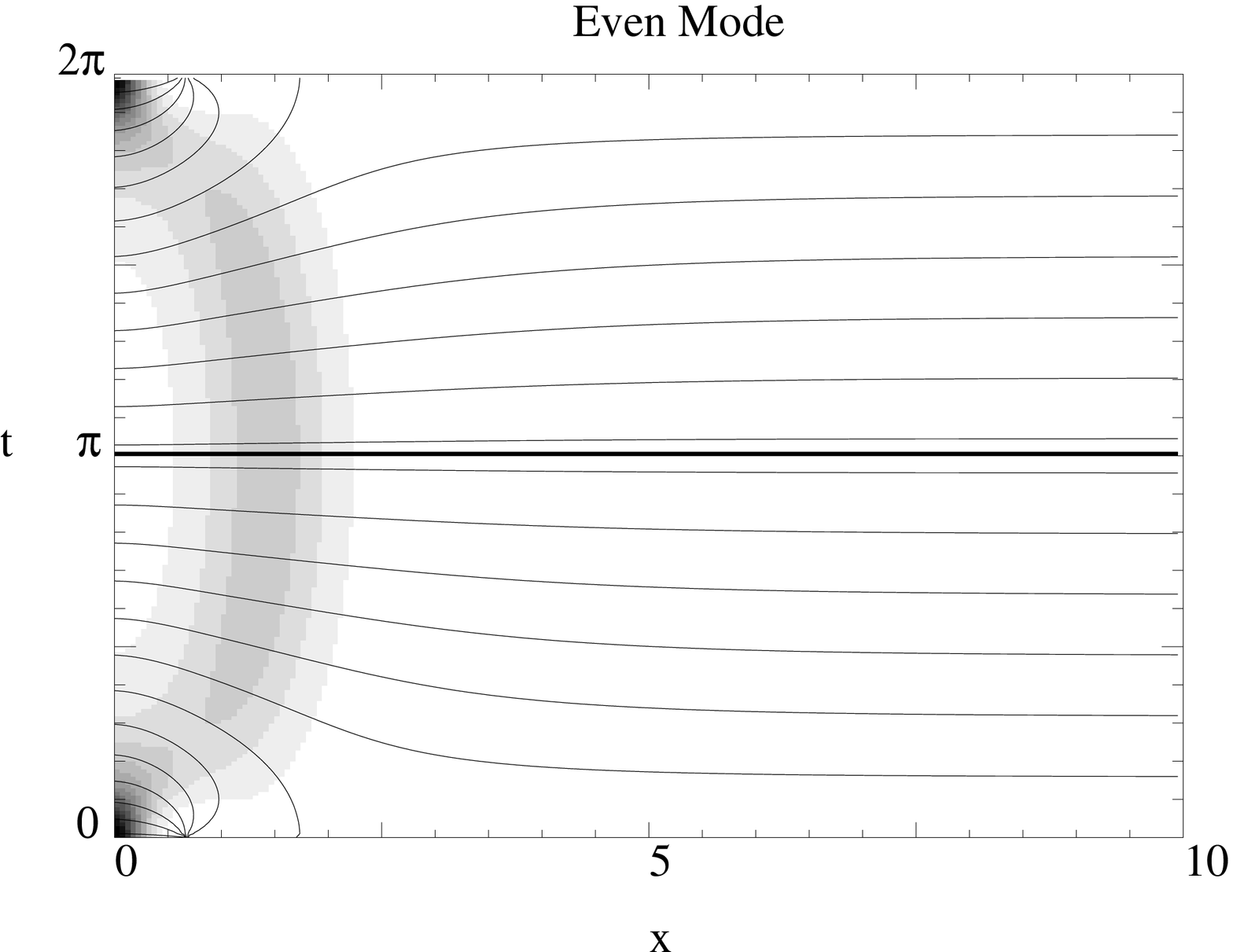,width=3 in}
\psfig{file=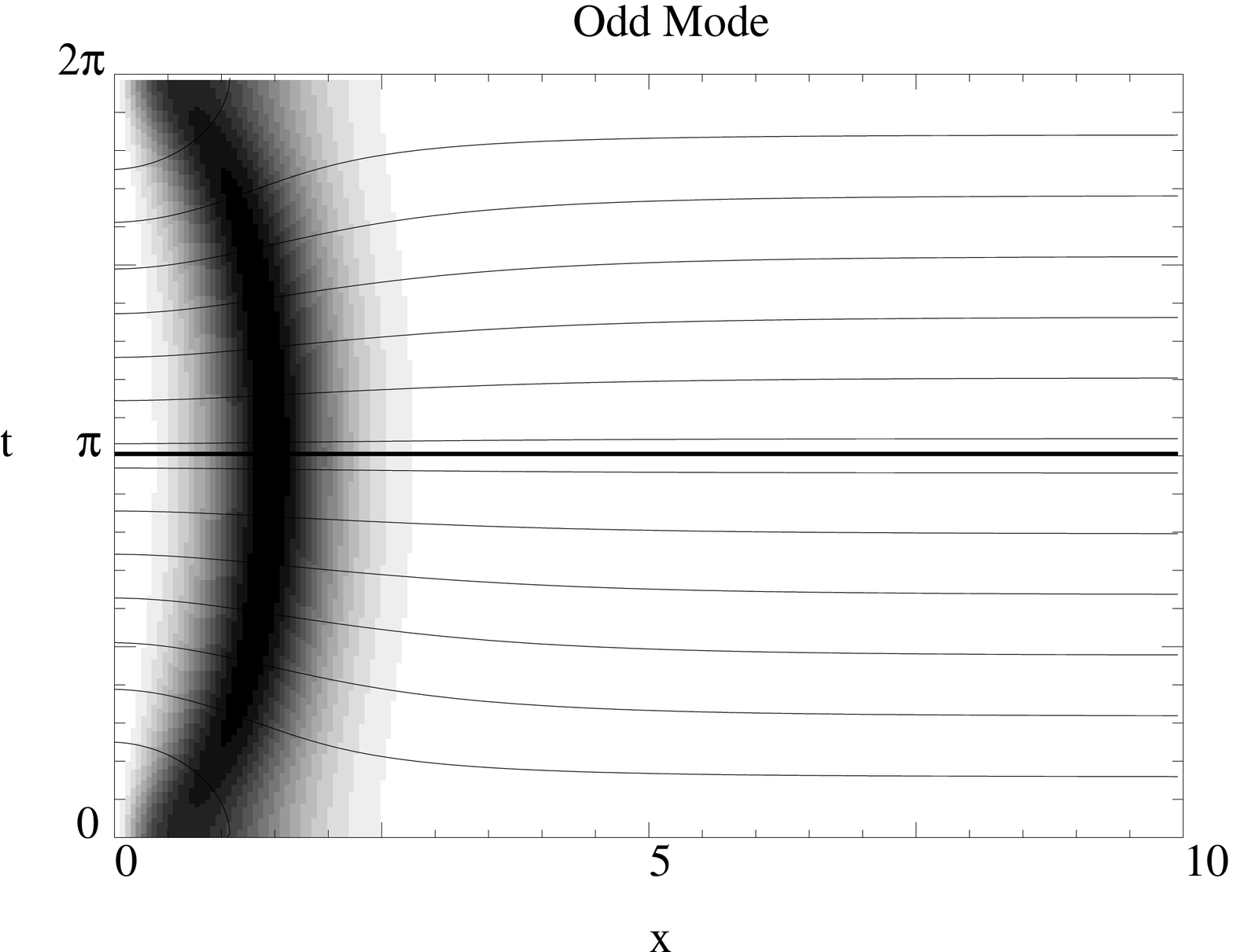,width=3 in}}
\end{center}
\caption{\em Equal phase contours for the even (left) and odd (right)
modes superimposed on a density plot of the corresponding square
modulus.  The parameter values are $\rho_1=1/\sqrt{2}$ and $\rho_2=1$.
For these values of $\rho_1$ and $\rho_2$, the Floquet multiplier 
is equal to $1$, and therefore these are periodic functions of $t$.}
\label{fig:R2modes}
\end{figure}

\begin{figure}[h]
\begin{center}
\mbox{\psfig{file=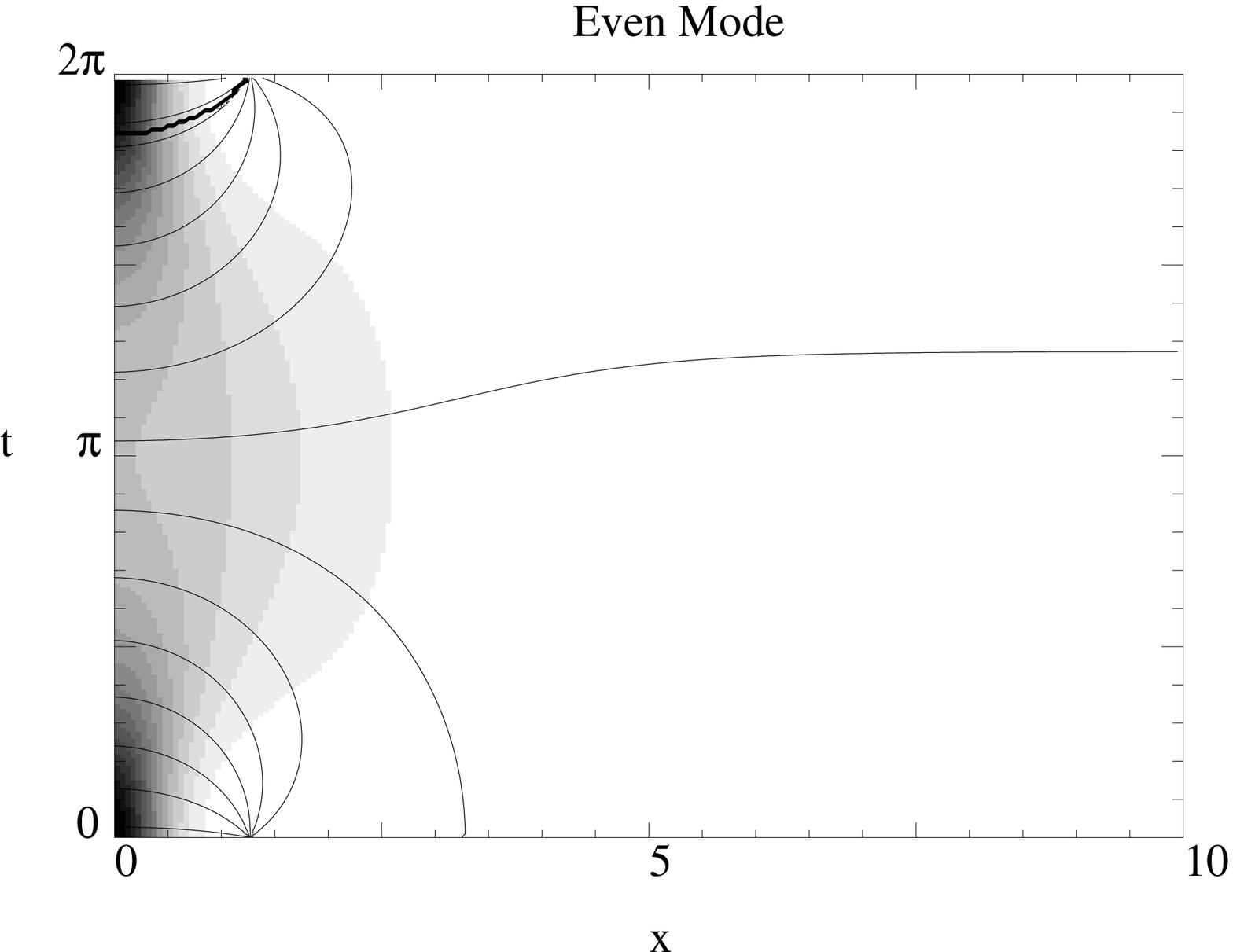,width=3 in}
\psfig{file=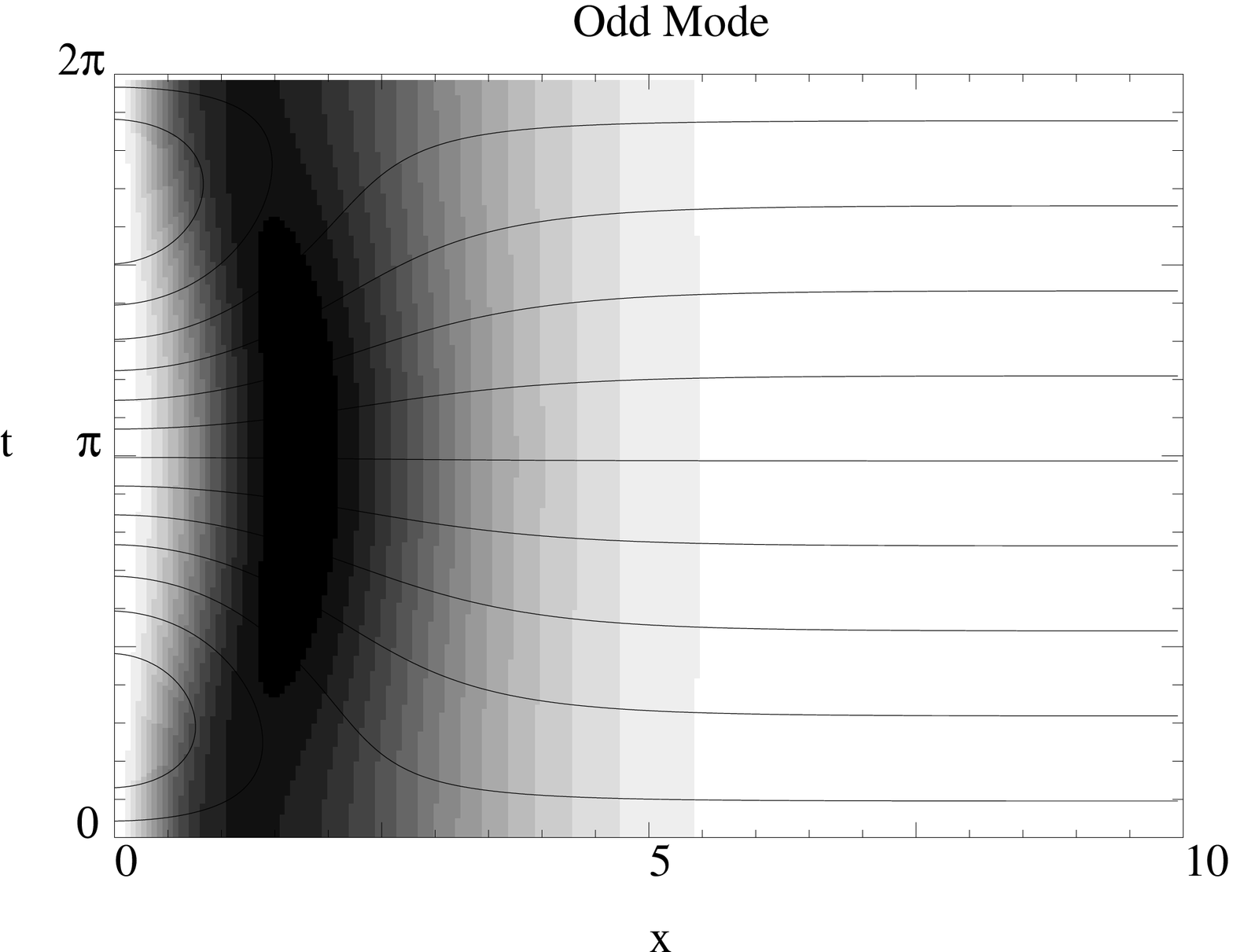,width=3 in}}
\end{center}
\caption{\em Equal phase contours for the even (left) and odd (right)
modes superimposed on a density plot of the corresponding square
modulus.  The parameter values are $\rho_1=1/4$ and $\rho_2=3/4$.
Here, the Floquet multiplier is not equal to $1$, and the modes are
not periodic in $t$, although they have Bloch form.}
\label{fig:R9modes}
\end{figure}

A particularly convenient orthonormal basis of the two-dimensional space
of bound states is given by the formulas
\begin{equation}
\begin{array}{rcl}
\Psi^{\rm (e)}_{\rm b}(x,t)&=&\displaystyle
\frac{1}{\sqrt{4(\rho_1+\rho_2)}}\left[
\frac{2}{\rho_1-\rho_2}a(x,t,-i\rho_1)+\frac{2}{\rho_2-\rho_1}a(x,t,-i\rho_2)\right]\,,\\
\\ \Psi^{\rm (o)}_{\rm b}(x,t)&=&\displaystyle
\frac{1}{\sqrt{4(\rho_1+\rho_2)}}\left[
\sqrt{\frac{\rho_2}{\rho_1}}\frac{2}{\rho_1-\rho_2}a(x,t,-i\rho_1) +
\sqrt{\frac{\rho_1}{\rho_2}}\frac{2}{\rho_2-\rho_1}a(x,t,-i\rho_2)\right]\,.
\end{array}
\end{equation}
In this case, the even symmetry of the potential $V_0(x,t)$ guarantees
that we may choose one basis element to be even and the other to be
odd; we are using superscripts ``(e)'' and ``(o)'' to refer to even
and odd functions of $x$ respectively.  These bound state solutions of
the linear Schr\"odinger equation are shown
in Figures~\ref{fig:R2modes} and \ref{fig:R9modes}.  
The two bound state modes are Bloch functions in $t$, with the same
Floquet multiplier $\exp(2i\beta_{\rm
b}L)=\exp(2i\rho_1^2L)=\exp(2i\rho_2^2L)$.  Note that, in reference to
the remark made at the end of \S \ref{sec:sub-separable}, the function
$\Psi_{\rm b}^{\rm (e)}(x,t)$ is proportional to
$\psi=2ib^{(1)}(x,t)$, which is a two-soliton solution of the
nonlinear Schr\"odinger equation
\begin{equation}
i\partial_t\psi +\frac{1}{2}\partial_x^2\psi + |\psi|^2\psi = 0\,.
\end{equation}
Correspondingly, $V_0(x,t)=-|\psi|^2$ is the self-consistent potential.
%The Floquet multiplier of this subspace is
%$\exp(2i\beta_{\rm b}L)=\exp(2i\rho_1^2L)=\exp(2i\rho_2^2L)$.  

It will also be useful to decompose the continuum into odd and even
parts.  Using the fact that $a(x,t,-\lambda)=a(-x,t,\lambda)$, define
\begin{equation}
\begin{array}{rcl}
\displaystyle 
\Psi^{\rm (e)}_{\rm d}(x,t,\lambda)&\doteq &\displaystyle
\frac{1}{\sqrt{2\pi (\lambda^2 + \rho_1^2)
(\lambda^2 + \rho_2^2)}}(a(x,t,\lambda)+a(x,t,-\lambda))\,,\\
\\
\displaystyle \Psi^{\rm (o)}_{\rm d}(x,t,\lambda)&\doteq &
\displaystyle \frac{1}{\sqrt{2\pi (\lambda^2 + \rho_1^2)
(\lambda^2 + \rho_2^2)}}(a(x,t,\lambda)-a(x,t,-\lambda))\,,
\end{array}
\end{equation}
where $\lambda \ge 0$.  The solutions $\Psi^{\rm e}_{\rm
d}(x,t,\lambda)$ and $\Psi^{\rm (o)}_{\rm d}(x,t,\lambda)$ are also 
Bloch functions in $t$ of period $L$ with Floquet multiplier
$\exp(-2i\lambda^2L)$.

These solutions of the unperturbed problem have the following inner
products \cite{physd}:
\begin{equation}
\begin{array}{rcl}
\displaystyle \langle \Psi^{\rm (e)}_{\rm b}(\cdot,t),\Psi^{\rm
(e)}_{\rm d}(\cdot,t,\lambda) \rangle &=& 0\,,\\\\ \displaystyle
\langle \Psi^{\rm (e)}_{\rm b}(\cdot,t),\Psi^{\rm (e)}_{\rm
b}(\cdot,t)\rangle & =& 1\,,\\\\ \displaystyle \langle \Psi^{\rm (e)}_{\rm
d}(\cdot,t,\lambda),\Psi^{\rm (e)}_{\rm d}(\cdot,t,\eta) \rangle
&=&\displaystyle \delta(\lambda-\eta)\,.
\end{array}
\end{equation}
In this latter relation it is assumed that both $\lambda$ and $\eta$
are positive real.  Similar relations hold among the odd solutions,
and of course everything even is orthogonal to everything odd.  If the
perturbation $W(x,t)$ is also even in $x$, then this observation will
allow us to treat the even and odd parts of the field $f(x,t)$ in
isolation to all orders in perturbation theory.

In our subsequent analysis of the coupled mode equations for this
family of periodic potentials, we shall assume that the perturbation
$W(x,t)$ is also an even function of $x$, and thus restrict attention
to the subspace of initial conditions, $f(x,0)$ which are either even
or odd in $x$. By the spatial symmetry of $V=V_0+W$, $f(x,t)$ has the
same parity as $f(x,0)$.  In analogy with the above derivation of
coupled mode equations, we can then expand $f(x,t)$ in terms of (even or odd)
modes of the unperturbed problem:
\begin{equation}
f(x,t)=B^{(\alpha)}_{\rm b}(t)\Psi^{(\alpha)}_{\rm b}(x,t)+
 \int_0^\infty B^{(\alpha)}_{\rm d}(t,\lambda)
\Psi^{(\alpha)}_{\rm d}(x,t,\lambda)\,d\lambda ,
\end{equation}
where $\alpha = {\rm e}$ if $f(x,0)$ is even and $\alpha = {\rm o}$ if
$f(x,0)$ is odd. Coupled mode equations for the amplitudes
$B^{(\alpha)}_{\rm b}(t),\ B^{(\alpha)}_{\rm d}(t,\lambda),
\lambda\in\R$ analogous to those derived in the absence of any
particular symmetry can then be derived by projecting the dynamical
system (\ref{eq:perturb}) onto these even and odd basis modes.  In
\S \ref{sec:periodic} we show, by using the Floquet factorization
of the unitary evolution associated with the unperturbed dynamics,
that these coupled mode equations can be reexpressed as the following
system ({\em c.f.} the system (\ref{eq:periodiccoupledmode})) which is
more amenable to our techniques:
\begin{equation}
\begin{array}{rcccl}
i\D_t A_{\rm b}(t) &+&2\beta_{\rm b} A_{\rm b}(t)&=& \displaystyle
M(t)A_{\rm b}(t) + \int_0^\infty N^{\rm (p)}(t,\lambda)A_{\rm
d}(t,\lambda)\,
d\lambda\,,\\\\
i\D_t A_{\rm d}(t,\eta) &-&2\eta^2A_{\rm d}(t,\eta)&=&\displaystyle
N^{\rm (p)}(t,\eta)^*A_{\rm b}(t) +
\int_0^\infty K^{\rm (p)}(t,\eta,\lambda)A_{\rm d}(t,\lambda)\,d\lambda\,,
\end{array}
\label{eq:system}
\end{equation}
where
\begin{equation}
A_{\rm b}(t)\doteq B_{\rm b}(t)e^{2i\beta_{\rm b}t}\,,\hspace{0.3 in}
A_{\rm d}(t,\lambda)\doteq B_{\rm d}(t,\lambda)e^{-2i\lambda^2 t}\,,
\end{equation}
and where the scalar coefficients, all periodic with period $L$, are
\begin{equation}
\begin{array}{rcl}
M(t) &=& \langle \Psi_{\rm b}(\cdot,t), W(\cdot,t)\Psi_{\rm
b}(\cdot,t)
\rangle \\\\
N^{\rm (p)}(t,\lambda) &=& \langle \Psi_{\rm b}(\cdot,t), W(\cdot,t)
\Psi_{\rm d}(\cdot,t,\lambda)\rangle e^{2i(\lambda^2+\beta_{\rm
b})t}\\\\
K^{\rm (p)}(t,\eta,\lambda) &=& \langle \Psi_{\rm d}(\cdot,t,\eta),
 W(\cdot,t)\Psi_{\rm d}(\cdot,t,\lambda)\rangle
e^{2i(\lambda^2-\eta^2)t}\,.
\end{array}
\label{eq:MNKs}
\end{equation} 

\noindent{\bf Remark:} To avoid cumbersome formulae, we have omitted the
superscripts (o) and (e), with the understanding that the
amplitudes correspond to either one type or the other, depending on
the parity of $f(x,0)$.  \pfend

The system (\ref{eq:system}) may be viewed as that governing a family
of oscillators: a single discrete oscillator whose amplitude is $A_{\rm
b}(t)$ coupled to a continuum of oscillators with amplitudes $A_{\rm
d}(t,\eta), \ \eta\in\R_+$.

In \S \ref{sec:analysis} we shall analyze the coupled mode system
(\ref{eq:system}), and determine the detailed asymptotic behavior of its
solutions for small $W(x,t)$ over different timescales. 

\section{Coupled mode equations for periodic potentials.}
\label{sec:periodic}
\setcounter{equation}{0} Consider a dynamical system of the form
(\ref{eq:schrod1}), where both the unperturbed and the perturbed
potential are time-periodic with the same period $L$. We will
encounter a concrete example of such a problem in \S
\ref{sec:application}.  Floquet theory \cite{Arnold} suggests the
introduction of a new time-periodic basis, with respect to which the
problem (\ref{eq:perturb}) becomes a periodic perturbation of {\em
autonomous} Hamiltonian system. This change of basis transforms the
problem at hand into one similar to that treated in \cite{ionization}
and \cite{KW}.  Similar methods are used along with resonance theory
in a weakly nonlinear setting in \cite{Sigal}.

\subsection{Floquet factorization.}
Let $\op{U}(t)$ denote the unitary evolution operator (or propagator)
of the unperturbed problem, so that for any $L^2(\R)$ function $f(x)$,
$f(x,t)= \op{U}(t)f(x)$ is the solution of the unperturbed problem
with $f(x,0)=f(x)$.  As a consequence of the periodicity, the
evolution operator can be factored into two operators on $L^2(\R)$:
\begin{equation}
\op{U}(t)= \op{P}(t)e^{-it\op{B}}\,
\end{equation}
where $\op{P}(t+L)=\op{P}(t)$, and $\op{B}$ is independent of $t$.
This factorization can be motivated by the observation that by
periodicity, there is an operator $\op{M}$ satisfying
$\op{U}(t+L)=\op{U}(t)\op{M}$, and that by setting $t=0$ in fact one
has $\op{M}=\op{U}(L)$.  Since $\op{U}(L)$ is unitary, one can find a
self-adjoint operator $\op{B}$ such that
$\op{M}=\op{U}(L)=e^{-iL\op{B}}$.  This operator $\op{B}$ in turn
defines the abelian unitary group $e^{-it\op{B}}$.  Now it is easy to
see that $\op{P}(t)=\op{U}(t)e^{it\op{B}}$ is a unitary operator
satisfying $\op{P}(t+L)=\op{P}(t)$.

Let $\op{W}(t)$ be the operator of multiplication by the correction
to the potential $W(x,t)$, and set $y(x,t)=\op{P}(t)^\dagger f(x,t)$.
Then, the perturbed equation (\ref{eq:schrod1}) becomes
\begin{equation}
i\D_ty-\op{B}y = \tilde{\op{W}}(t)y\,,
\label{eq:transf}
\end{equation}
where $\tilde{\op{W}}(t)\doteq \op{P}(t)^\dagger \op{W}(t)\op{P}(t)$,
a ``dressing'' of $\op{W}(t)$.  The form of (\ref{eq:transf}) is
similar to the type of problem treated in \cite{ionization} and
\cite{KW}.  The ``unperturbed Hamiltonian'' $\op{B}$ is
time-independent, and self-adjoint.  The perturbation
$\tilde{\op{W}}(t)$ is localized, self-adjoint, and time-periodic
because the periods of $\op{P}(t)$ and $\op{W}(t)$ are equal.
Typically the perturbation contains frequency components at all
overtones of the fundamental frequency, and thus the version of the
theory described in \cite{KW} is most appropriate.  The key
qualitative difference between the present situation and that treated
in \cite{KW} is that here the unperturbed operator can have multiple
bound states.  We will soon introduce a symmetry that removes this
difficulty from the scope of this paper.  However, the methods of
\cite{ionization} and \cite{KW} can be extended to give results on
radiation damping due to the coupling of multiple discrete modes to
the continuum for a general class of spatially localized and
time-dependent perturbations \cite{KWmultimode}.

In fact, one can simplify the problem even further by invoking the
spectral theorem for the self-adjoint operator $\op{B}$.  This
guarantees the existence of an isomorphism $\op{V}:L^2(\R)\rightarrow
L^2(\Sigma,d\mu)$ to the space of square integrable functions on some
set $\Sigma$ with measure $d\mu$, such that
$\op{V}\op{B}=\op{T}\op{V}$ where $\op{T}$ is a real diagonal operator
on $L^2(\Sigma,d\mu)$ ({\em i.e. } an operator of multiplication by a
function from $\Sigma$ to $\R$).  Setting $z(t)=\op{V}y(t)$,
we find the equation
\begin{equation}
i\D_tz-\op{T}z = \op{V}\tilde{\op{W}}(t)\op{V}^\dagger z\,.
\label{eq:transfdiag}
\end{equation}
In quantum mechanics, making the transformation from
(\ref{eq:perturb}) to (\ref{eq:transfdiag}) to facilitate the study of
perturbations is known as going from the {\em Schr\"odinger picture}
into the {\em interaction picture}.

In the particular example we will analyze in detail, arising from even
perturbations of the two-soliton even potential described at the end
of \S \ref{sec:separable}, the operator $\op{B}$ has a single
degenerate eigenvalue of $-2\beta_{\rm b}<0$ of geometric multiplicity
two.  By restricting separately to even and odd spaces of initial
conditions (possible because the potential $V_0(x,t)+W(x,t)$ is
symmetric in $x$) the problem is reduced to one which, formally, is
precisely of the type studied in \cite{KW}.  We may then apply the
methods developed in \cite{KW} (subject to some appropriate hypotheses)
without modification.

We now use our explicit knowledge developed in \S
\ref{sec:separable} of the unitary propagator $\op{U}(t)$
corresponding to a time-periodic separable potential $V_0(x,t)$ to
find the operators $\op{P}(t)$ and $\op{B}$, and then to diagonalize
$\op{B}$.  This effectively implements the program described above
and casts the perturbed problem (\ref{eq:perturb}) in a form
(\ref{eq:transfdiag}) more suitable for analysis.  We begin with the
observation that each element of the basis of solutions of the
unperturbed problem is a Bloch function or Floquet mode.  We have:
\begin{equation} 
\Psi_{\rm
d}(x,t+L,\lambda)=e^{-2i\lambda^2L}\Psi_{\rm d}(x,t,\lambda)
\label{eq:PsitplusL} 
\end{equation} 
where $\lambda$ is arbitrary real.  Also, we have 
\begin{equation} 
a(x,t+L,-i\rho_k)=e^{2i\rho_k^2L}a(x,t,-i\rho_k)\,.
\label{eq:PhitplusL}
\end{equation} 
Note that the commensurability relations (\ref{eq:comensurate}) imply
quite generally that the Floquet multipliers $\exp(2i\rho_k^2L)$ {\em
are all equal}.  This generalizes the observation made above in the
context of the two-soliton potentials.  This means that the entire
$M$-dimensional subspace of bound states consists of degenerate
Floquet modes.  In particular, the elements of any orthonormal basis
$\{\Psi_{{\rm b},k}(x,t)\,,k=1,\dots,M\}$ have the same Floquet
multiplier\footnote{Recall that the Floquet exponents are not unique
but that the Floquet multipliers are. Identification of the Floquet
exponents with a single number $\beta_{\rm b}>0$ amounts to a
particular choice of branch of the logarithm.} $\exp(2i\beta_{\rm
b}L)$.  It now follows from (\ref{eq:PsitplusL}) and
(\ref{eq:PhitplusL}) that the functions defined by
\begin{equation}
\begin{array}{rcl}
\Psi^{\rm (p)}_{\rm d}(x,t,\lambda)&\doteq &
e^{2i\lambda^2 t}\Psi_{\rm d}(x,t,\lambda)\,,\\\\
\Psi^{\rm (p)}_{{\rm b},k}(x,t)&\doteq &
e^{-2i\beta_{\rm b}t}\Psi_{{\rm b},k}(x,t)\,,
\end{array}
\label{eq:periodicdef}
\end{equation}
are time-periodic with period $L$, as denoted by the superscript
``(p)''.  

As described in \S \ref{sec:separable}, the solution of the
unperturbed problem with initial data $f(x)$ is expanded as
\begin{equation}
\begin{array}{rcl}
f(x,t)&=&\op{U}(t)f(x)\\\\
&=&\displaystyle \sum_{k=1}^M\langle\Psi_{{\rm
b},k}(\cdot,0),f(\cdot)\rangle \Psi_{{\rm
b},k}(x,t)+\int_{-\infty}^\infty \langle \Psi_{\rm
d}(\cdot,0,\lambda),f(\cdot)\rangle \Psi_{\rm
d}(x,t,\lambda)\,d\lambda\\\\ &=&\displaystyle
\sum_{k=1}^M\langle\Psi^{\rm (p)}_{{\rm b},k}(\cdot,0),e^{2i\beta_{\rm
b}t}f(\cdot)\rangle \Psi^{\rm (p)}_{{\rm
b},k}(x,t)+\int_{-\infty}^\infty \langle \Psi^{\rm (p)}_{\rm
d}(\cdot,0,\lambda),e^{-2i\lambda^2 t}f(\cdot)\rangle \Psi^{\rm
(p)}_{\rm d}(x,t,\lambda)\,d\lambda\,.
\end{array}
\end{equation}
We now use the completeness relation at $t=0$ to factor $\op{U}(t)$
as $\op{P}(t)e^{-it\op{B}}$ where
\begin{equation}
e^{-it\op{B}}f(x)=\sum_{k=1}^M\langle \Psi^{\rm (p)}_{{\rm
b},k}(\cdot,0), e^{2i\beta_{\rm b}t}f(\cdot)\rangle \Psi^{\rm (p)}_{{\rm
b},k}(x,0) + \int_{-\infty}^\infty \langle \Psi^{\rm (p)}_{\rm
d}(\cdot,0,\lambda), e^{-2i\lambda^2 t}f(\cdot)\rangle \Psi^{\rm
(p)}_{\rm d}(x,0,\lambda)\,d\lambda\,,
\end{equation}
and 
\begin{equation}
\op{P}(t)g(x)=\sum_{k=1}^M\langle \Psi^{\rm (p)}_{{\rm
b},k}(\cdot,0),g(\cdot)\rangle \Psi^{\rm (p)}_{{\rm b},k}(x,t)
+\int_{-\infty}^\infty \langle \Psi^{\rm (p)}_{\rm
d}(\cdot,0,\lambda),g(\cdot)\rangle \Psi^{\rm (p)}_{\rm
d}(x,t,\lambda)\,d\lambda\,.
\label{eq:pdef}
\end{equation}
We have used, several times, the fact that at $t=0$ there is no
distinction between the basis elements and their periodic counterparts
defined by (\ref{eq:periodicdef}).  It is easy to see that
$f(x,t)=\op{U}(t)f(x)=\op{P}(t)e^{-it\op{B}}f(x)$ is the solution of
the unperturbed initial value problem with data $f(x)\in L^2(\R)$ and
that $\op{P}(t)$ is periodic with period $L$ and
$\op{U}(L)=e^{-iL\op{B}}$.  

The generator of the abelian unitary group
$e^{-it\op{B}}$ is
\begin{equation}
\begin{array}{rcl}
\displaystyle \op{B}f(x)&=&\displaystyle 
i\frac{d}{dt}e^{-it\op{B}}f(x)\Bigg|_{t=0}\\\\
&=&\displaystyle
\sum_{k=1}^M\langle \Psi^{\rm (p)}_{{\rm b},k}(\cdot,0),-2\beta_{\rm b}f(\cdot)
\rangle \Psi^{\rm (p)}_{{\rm b},k}(x,0) + \int_{-\infty}^\infty
\langle \Psi^{\rm (p)}_{\rm d}(\cdot,0,\lambda),2\lambda^2 f(\cdot)\rangle
\Psi^{\rm (p)}_{\rm d}(x,0,\lambda)\,d\lambda\\\\
&=&\displaystyle
\sum_{k=1}^M\langle \Psi_{{\rm b},k}(\cdot,0),-2\beta_{\rm b}f(\cdot)
\rangle \Psi_{{\rm b},k}(x,0) + \int_{-\infty}^\infty
\langle \Psi_{\rm d}(\cdot,0,\lambda),2\lambda^2 f(\cdot)\rangle
\Psi_{\rm d}(x,0,\lambda)\,d\lambda\,.
\end{array}
\end{equation}
In the last step we have dropped the superscripts ``(p)'' since
everything is evaluated at $t=0$.  
This formula for the self-adjoint
operator $\op{B}$ makes clear its spectral decomposition.  The
isomorphism $\op{V}$ takes a function $g(x)\in L^2(\R)$ to a function
$A_{\rm d}(\lambda)$ for $\lambda\in\R$ and a set of $M$ numbers
$A_{{\rm b},k}$ for $k=1,\dots,M$ defined by
\begin{equation}
A_{\rm d}(\lambda)\doteq
\langle \Psi^{\rm (p)}_{\rm d}(\cdot,0,\lambda),g(\cdot)
\rangle\,,
\end{equation}
and for $k=1,\dots,M$,
\begin{equation}
A_{{\rm b},k}\doteq
\langle \Psi^{\rm (p)}_{{\rm b},k}(\cdot,0),g(\cdot)\rangle\,.
\end{equation}
The diagonal operator $\op{T}$ is then simply defined by
\begin{equation}
\op{T}\left[\begin{array}{c}A_{\rm d}(\lambda)\\A_{{\rm b},1}\\
\vdots\\A_{{\rm b},M}\end{array}\right]=
\left[\begin{array}{cccc}
2\lambda^2 &&&\\
&-2\beta_{\rm b}&&\\
&&\ddots&\\
&&&-2\beta_{\rm b}\end{array}\right]
\left[
\begin{array}{c}A_{\rm d}(\lambda)\\
A_{{\rm b},1}\\\vdots\\A_{{\rm b},M}
\end{array}\right]\,.
\end{equation}

It is now easy to use the definition of the unitary periodic operator
$\op{P}(t)$ and the unitary isomorphism $\op{V}$, along with the
completeness relation to compute the dressed operator
$\op{V}\tilde{\op{W}}(t)\op{V}^\dagger$ and thus write the perturbed problem
(\ref{eq:perturb}) in the simple form (\ref{eq:transfdiag}).  The
dynamical unknowns are in the range of $\op{V}$, the space
$L^2(\Sigma,d\mu)$, and are given in terms of $f(x,t)$, the solution of
(\ref{eq:perturb}), by
\begin{equation}
\begin{array}{rcccl}
A_{{\rm b},k}&=&(\op{V}\op{P}(t)^\dagger f(\cdot,t))_{{\rm b},k}&=&
\langle \Psi^{\rm (p)}_{{\rm b},k}(\cdot,t),f(\cdot,t)\rangle\,,\\
A_{\rm d}(\lambda)&=&(\op{V}\op{P}(t)^\dagger f(\cdot,t))_{\rm d}(\lambda)&=&
\langle \Psi^{\rm (p)}_{\rm d}(\cdot,t,\lambda),f(\cdot,t)\rangle\,.
\end{array}
\end{equation}
When $f(x,t)$ satisfies (\ref{eq:perturb}), these quantities satisfy
the system
\begin{equation}
\label{eq:periodiccoupledmode}
\begin{array}{rcl}
i\partial_t \vc{A}_{\rm b}+2\beta_{\rm b}\vc{A}_{\rm b} &=&
\displaystyle \mat{M}(t)\vc{A}_{\rm b} +\int_{-\infty}^\infty
A_{\rm d}(\lambda)\vc{N}^{\rm (p)}(t,\lambda)\,d\lambda\,,\\\\
i\partial_t A_{\rm d}(\eta)-2\eta^2A_{\rm d}(\eta)&=&
\displaystyle \vc{N}^{\rm (p)}(t,\eta)^\dagger \vc{A}_{\rm b} +
\int_{-\infty}^\infty K^{\rm (p)}(t,\eta,\lambda)A_{\rm d}(\lambda)\,d\lambda
\,,
\end{array}
\end{equation}
where $\vc{A}_{\rm b}$ is the vector of components $A_{{\rm
b},1},\dots,A_{{\rm b},M}$, and the time-periodic matrix elements are
defined in terms of (\ref{eq:matrixelts}) by
\begin{equation}
\begin{array}{rcl}
\vc{N}^{\rm (p)}(t,\lambda)&\doteq &e^{2i(\lambda^2+\beta_{\rm
b})t}\vc{N}(t,\lambda)\,,\\\\ K^{\rm (p)}(t,\eta,\lambda)&\doteq
&e^{2i(\lambda^2-\eta^2)t}K(t,\eta,\lambda)\,.
\end{array}
\end{equation}
The periodicity of these matrix elements when $W(x,t)$ is periodic
with period $L$ is also clear from these explicit formulas and the
Bloch relations (\ref{eq:PsitplusL}) and (\ref{eq:PhitplusL}) for the
basis of solutions; these imply similar ones for the matrix elements
defined by (\ref{eq:matrixelts}).  We have
\begin{equation}
\begin{array}{rcl}
\mat{M}(t+L)&=&\mat{M}(t)\,,\\\\
\vc{N}(t+L,\lambda)&=&e^{-2i(\lambda^2+\beta_{\rm b})L}\vc{N}
(t,\lambda)\,,\\\\
K(t+L,\eta,\lambda)&=&e^{2i(\eta^2-\lambda^2)L}K(t,\eta,\lambda)\,.
\end{array}
\end{equation}
Of course, the right-hand side of (\ref{eq:periodiccoupledmode}) is
just the operator $\op{V}\tilde{\op{W}}(t)\op{V}^\dagger$ operating on
the dynamical unknowns.  Similarly, the perturbation operator
$\tilde{\op{W}}(t)$ operating in the space $L^2(\R)$ can be explicitly
written as
\begin{equation}
\begin{array}{rcl}
\tilde{\op{W}}(t)f(x)&=&\displaystyle
\sum_{k=1}^M\sum_{l=1}^M M_{k,l}(t)
\langle \Psi_{{\rm b},l}(\cdot,0),f(\cdot)\rangle
\Psi_{{\rm b},k}(x,0)\\\\
&+&\displaystyle
\sum_{k=1}^M\int_{-\infty}^\infty
N^{\rm (p)}_k(t,\eta)\langle\Psi_{\rm d}(\cdot,0,\eta),f(\cdot)
\rangle\Psi_{{\rm b},k}(x,0)\,d\eta\\\\
&+&\displaystyle
\sum_{l=1}^M\int_{-\infty}^\infty
N^{\rm (p)}_l(t,\lambda)^*\langle\Psi_{{\rm b},l}(\cdot,0),f(\cdot)
\rangle\Psi_{\rm d}(x,0,\lambda)\,
d\lambda\\\\
&+&\displaystyle
\int_{-\infty}^\infty\int_{-\infty}^\infty
K^{\rm (p)}(t,\lambda,\eta)\langle\Psi_{\rm d}(\cdot,0,\eta),
f(\cdot)\rangle\Psi_{\rm d}(x,0,\lambda)\,d\lambda\,d\eta\,.
\end{array}
\end{equation}

For the special choice of $V_0(x,t)$ discussed at the end of \S
\ref{sec:separable}, evenness implies that there is one bound state of
each parity.  If the perturbation $W(x,t)$ also has even symmetry in
$x$, the coupled mode system (\ref{eq:periodiccoupledmode}) reduces to
a system of the type (\ref{eq:system}) when the initial condition is
restricted to either even or odd parity.  It is easily checked that
the unknowns as defined above correspond exactly to those defined in
\S \ref{sec:separable} for the system (\ref{eq:system}).

\section{Analysis of the Coupled Mode Equations}
\label{sec:analysis}
\setcounter{equation}{0} In this section we study the structural
instability of the even and odd breather modes introduced at the end
of \S \ref{sec:separable} associated with the two-soliton
time-periodic even potentials.  We first give a simple argument valid
for short times that in the presence of a perturbation $W(x,t)$ to the
potential $V_0(x,t)$, the bound state begins to decay initially.  We
then seek to capture the dynamics for longer times, primarily to show
that this initial phase of decay does not reverse itself, but takes on
a different, exponentially decaying, character.  The decay will be
first calculated formally, using asymptotic expansions and the method
of multiple scales.  Then, using the results of Kirr and Weinstein
\cite{KW}, we show that at least in the odd case, it is possible to
make statements about the decay process that are valid globally in
time.  In particular, these arguments will rigorously justify the
formal results for the odd case, and will show that the exponential
decay model is only a valid approximation until it becomes smaller
than the dispersive part of the solution.  The bound state ultimately
dies algebraically in time, qualitatively indistinguishable from the
dispersive components of the solution to which it is orthogonal.

\subsection{Small time analysis.  The watched pot effect.}
A simple calculation carried out at the level of the coupled mode
equations (\ref{eq:system}) shows that the effect of the perturbation
is to cause the bound state to decay immediately both forward and
backward in time.  More complicated calculations will be required to
show that the decay does not stop or reverse for longer times,
although it takes on a different character.  The approach in the small
time analysis is simply to expand the solution in Taylor series:
\begin{equation}
\begin{array}{rcl}
A_{\rm b}(t)&=&A_{\rm b}(0) + c_1t + c_2 t^2 + O(t^3)\,,\\
A_{\rm d}(\lambda,t)&=&d_1(\lambda)t + O(t^2)\,,
\end{array}
\end{equation}
and use the (known) Taylor expansions of the matrix elements, in particular,
\begin{equation}
\begin{array}{rcl}
M(t)&=&M(0)+M'(0)t + O(t^2)\,,\\
N^{\rm (p)}(t,\lambda)&=&N^{\rm (p)}(0,\lambda) + O(t)\,.
\end{array}
\end{equation}
Substituting these series into (\ref{eq:system}), one finds:
\begin{equation}
\begin{array}{l}
\displaystyle\left[ic_1+2\beta_{\rm b}A(0)-M(0)A_{\rm b}(0)\right]\,\,+\\\\
\displaystyle
\hspace{0.3 in}\left[2ic_2 +2\beta_{\rm b}c_1-M(0)c_1-M'(0)A_{\rm b}(0)
-\int_0^\infty d_1(\lambda)N^{\rm (p)}(0,\lambda)\,d\lambda\right]t\,\,=\,\,
O(t^2)\,,
\\\\
id_1(\eta) -N^{\rm (p)}(0,\eta)^*A_{\rm b}(0)\,\, =\,\, O(t)\,.
\end{array}
\end{equation}
Solving for $c_1$ and $c_2$ yields an approximation for $A_{\rm b}(t)$, valid
for small $t$:
\begin{equation}
\begin{array}{l}
\displaystyle A_{\rm b}(t)=A_{\rm b}(0)\Bigg[1-i\left(M(0)
-2\beta_{\rm b}\right)t\\\\
\displaystyle\hspace{0.3 in} -\,\, \frac{1}{2}\left(iM'(0)+\left(M(0)-2\beta_{\rm b}\right)^2 +
\int_0^\infty |N^{\rm (p)}(0,\lambda)|^2\,d\lambda\right)t^2 + O(t^3)\Bigg]\,.
\end{array}
\end{equation}
It easily follows that
\begin{equation}
|A_{\rm b}(t)|^2 = |A_{\rm b}(0)|^2\left[1-t^2\int_0^\infty |N^{\rm (p)}(0,\lambda)|^2\,d\lambda + O(t^3)\right]\,.
\end{equation}
Note that smallness of the perturbation is not exploited in these
calculations.  This Taylor expansion shows that the initial phase of
the evolution is a process of radiative decay, since $|A_{\rm
b}(t)|^2<|A_{\rm b}(0)|^2$ for all nonzero $t$ in some neighborhood of
$t=0$.  The decay is symmetric in time.  

The fact that the decay is an order $O(t^2)$ effect is quite
general\footnote{In the general setting, the decay is a simple consequence of
the Cauchy-Schwarz inequality.  One supposes that $\op{U}(t)$ is the
unitary propagator of the possibly time-dependent unperturbed problem:
\begin{displaymath}
i\op{U}_t(t)\phi^0=\op{H}_0(t)\op{U}(t)\phi^0\,,
\end{displaymath}
for all states $\phi^0$.  One then considers the perturbed equation
\begin{displaymath}
i\psi_t=(\op{H}_0(t)+\op{W}(t))\psi\,,
\end{displaymath}
by setting $\psi(t)=\op{U}(t)\phi(t)$, giving the ``interaction picture'' equation
\begin{displaymath}
i\phi_t=\op{U}(t)^\dagger\op{W}(t)\op{U}(t)\phi\,,
\end{displaymath}
which one solves by Taylor series in $t$.  The result is:
\begin{displaymath}
\phi(t)=\left(\op{I}-i\op{W}(0)t + \frac{t^2}{2}\left(-i\op{W}'(0)-\op{W}(0)^2+[\op{H}_0(0),\op{W}(0)]\right)+O(t^3)\right)\phi(0)\,.
\end{displaymath}
The probability of remaining in the unperturbed state is then found to be (using self-adjointness of both $\op{W}(0)$ and $\op{H}_0(0)$)
\begin{displaymath}
\begin{array}{rcl}
|\langle\op{U}(t)\phi(0),\op{U}(t)\phi(t)\rangle|^2&=&
|\langle \phi(0),\phi(t)
\rangle|^2\\
&=&\|\phi(0)\|_2^4 - (\|\op{W}(0)\phi(0)\|_2^2\|\phi(0)\|_2^2-
|\langle\phi(0),\op{W}(0)\phi(0)\rangle|^2)t^2 + O(t^3)\,.
\end{array}
\end{displaymath}
This quantity is initially decreasing in time as a consequence of the
Cauchy-Schwarz inequality.}  and well-known in the perturbation
theory of stationary Schr\"odinger equations.  It has an interesting
interpretation in the quantum theory of ideal measurements, the
so-called ``watched pot effect''.  Suppose that an ideal measurement
is made at some point during the evolution of the wave function to
determine whether the state is bound, and the measurement yields a
positive result.  The probability of a positive result at time $t$ is
$|A_{\rm b}(t)|^2/|A_{\rm b}(0)|^2$.  The theory of ideal measurements
says that as a consequence of the measurement disturbing the system,
the wave function ``collapses'' upon a positive result to the bound
state, and evolution of the wave function according to the
Schr\"odinger equation continues from this ``reset'' bound state.  One
may then try to determine the asymptotic effect of making many such
measurements in a finite time interval.  In particular, we can ask
about the limiting probability of finding the system in the bound
state after {\em each} of $n$ ideal measurements performed at times
$t_n=T/n$, as $n\rightarrow \infty$.  After each positive result, the
wave function collapses and the experiment is restarted.  The
Schr\"odinger evolution takes place over short time intervals so it is
appropriate to replace the probability in each interval $p(t)$ by its
short-time approximation $p(t)=1-(\alpha t)^2+O(t^3)$.  The $n$
measurements are independent events, so the probability of always
finding the system bound after each measurement is simply
\begin{equation}
P_n=p(T/n)^n\,.
\end{equation}
Because the ``time slice'' decay probability $1-p(t)$ is quadratic in
$t$, $P_n$ tends to unity\footnote{The superlinear nature of the decay
probability is important.  If $p(t)=1-|\alpha t| + O(t^2)$, then $P_n$
tends to $e^{-|\alpha T|}<1$ instead.} as $n\rightarrow\infty$,
regardless of the value of $T$.  So if the measurements are performed
infinitely often, the decay of the bound state never occurs.  The
quantum ``watched pot'' never boils.

\subsection{Multiple scales analysis.}
\label{sec:sub-ms}
We begin the multiple scales analysis by assuming that the correction
$W(x,t)=W(x,t;\epsilon)$ to the potential energy has an expansion in a
small parameter, $\epsilon$ (see, for example, equation
(\ref{eq:expansionofW})\ ): \be W(x,t;\epsilon)\ =\ \epsilon W_1(x,t)\
+\ O(\epsilon^2).  \no\ee It then follows that the coupling
coefficient functions in (\ref{eq:system}) have formal expansions for
small $\epsilon$:
\begin{equation}
\begin{array}{rcl}
M(t)&=&\epsilon M_1(t) + \epsilon^2 M_2(t) + O(\epsilon^3)\,,
\\ 
N^{\rm (p)}(t,\lambda)&=&\epsilon N^{\rm (p)}_1(t,\lambda) + O(\epsilon^2)\,,\\ 
K^{\rm (p)}(t,\eta,\lambda)&=&\epsilon K^{\rm (p)}_1(t,\eta,\lambda) +
O(\epsilon^2)\,.
\end{array}
\label{eq:MNKexpansions}
\end{equation}
Here, $M_1(t), N^{\rm (p)}_1(t)$ and $K^{\rm (p)}_1(t,\eta,\lambda)$
correspond to the expressions for $M(t)$, $N^{(\rm p)}$ and $K^{(\rm
p)}$ in (\ref{eq:MNKs}) with $W$ replaced by $W_1$.

The amplitudes $A_{\rm b}(t;\epsilon)$ and $A_{\rm d}(t,\lambda;\epsilon)$
are assumed to have asymptotic expansions of the form
\begin{equation}
\begin{array}{rcl}
A_{\rm b}(t;\epsilon)&=&A_{\rm b}^{(0)}(T_0,T_1,T_2,\dots)+\epsilon
A_{\rm b}^{(1)}(T_0,T_1,T_2,\dots)+\epsilon^2A_{\rm b}^{(2)}(T_0,T_1,T_2,\dots)+ O(\epsilon^3)\,,\\
A_{\rm d}(t,\lambda;\epsilon)&=&A^{(0)}_{\rm d}(T_0,T_1,T_2,\dots,\lambda)+
\epsilon A_{\rm d}^{(1)}(T_0,T_1,T_2,\dots,\lambda) + O(\epsilon^2)\,,
\end{array}
\label{eq:expansion}
\end{equation}
where the $T_k\doteq \epsilon^kt$ are time scale variables.  Such
expansions of given functions $A_{\rm b}(t;\epsilon)$ and $A_{\rm
d}(t,\lambda;\epsilon)$ are highly nonunique.  However the guiding
principle of the method of multiple scales (see, for example,
\cite{KevorkianCole}) stipulates that the dependence of the various
terms on the ``slow'' times $T_1$, $T_2$, and so on is chosen so that
each term is uniformly bounded as a function of the ``fast'' time
$T_0$.  This procedure is quite systematic, and is supposed to keep
the error terms in any truncation uniformly small in time intervals
where $T_k$ is bounded for some $k$ as $\epsilon$ tends to zero.  We
will see by comparison with the rigorous results that this formal
procedure indeed works as advertised.  

One now substitutes these expansions into the system (\ref{eq:system})
and expands the time derivative operating on the expansion coefficients in
(\ref{eq:expansion}) according to the chain rule:
\begin{equation}
\partial_t=  \partial_{T_0}+\epsilon\partial_{T_1}+\epsilon^2\partial_{T_2}
+\dots\,.
\label{eq:chain}
\end{equation}
The coupling coefficients, being all periodic functions of $t$ with
period $L$ independent of $\epsilon$, are taken to be explicit
functions of $t=T_0$ {\em only}.  Substituting them into (\ref{eq:system})
along with the expansions (\ref{eq:expansion}) and the chain rule formula
(\ref{eq:chain}), and equating terms with the same powers of
$\epsilon$ leads to a hierarchy of
equations:
\begin{eqnarray}
\label{eq:order0}
O(1):&&\left\{\begin{array}{rcccl}
i\D_{T_0}A_{\rm b}^{(0)}&+&2\beta_{\rm b}A_{\rm b}^{(0)}&=&
0\,,\\\\
i\D_{T_0}A_{\rm d}^{(0)}(\eta)&-&2\eta^2A_{\rm d}^{(0)}(\eta)&=&0\,,
\end{array}\right.\\\nonumber\\
\label{eq:order1}
O(\epsilon):&&\left\{\begin{array}{rcccl}
i\D_{T_0}A_{\rm b}^{(1)}&+&2\beta_{\rm b}A_{\rm b}^{(1)}&=&
\displaystyle -i\D_{T_1}A_{\rm b}^{(0)}+M_1(T_0)A_{\rm b}^{(0)} \\\\
&&&&\displaystyle\hspace{0.2 in}+\,\, 
\int_0^\infty N_1^{\rm (p)}(T_0,\lambda)A_{\rm d}^{(0)}(\lambda)\,d\lambda\,, 
\\\\
i\D_{T_0}A_{\rm d}^{(1)}(\eta)&-&2\eta^2A_{\rm d}^{(1)}(\eta)&=&\displaystyle
-i\D_{T_1}A_{\rm d}^{(0)}(\eta)+N_1^{\rm (p)}(T_0,\eta)^*A_{\rm b}^{(0)}\\\\
&&&&\displaystyle\hspace{0.2 in} +\,\,
\int_0^\infty K_1^{\rm (p)}(T_0,\eta,\lambda)
A_{\rm d}^{(0)}(\lambda)\,d\lambda\,, 
\end{array}\right.\\\nonumber\\
\label{eq:order2}
O(\epsilon^2):&&\left\{\begin{array}{rcccl}
i\D_{T_0}A_{\rm b}^{(2)}&+&2\beta_{\rm b}A_{\rm b}^{(2)}&=&
\displaystyle 
-i\D_{T_1}A_{\rm b}^{(1)}-i\D_{T_2}A_{\rm b}^{(0)}\\\\
&&&&\displaystyle\hspace{0.2 in}+\,\,
M_2(T_0)A_{\rm b}^{(0)}+M_1(T_0)A_{\rm b}^{(1)}\\\\ 
&&&&\displaystyle\hspace{0.2 in}+\,\,
\int_0^\infty N_1^{\rm (p)}(T_0,\lambda)A_{\rm d}^{(1)}(\lambda)\,d\lambda\\\\
&&&&\displaystyle\hspace{0.2 in}+\,\,
\int_0^\infty N_2^{\rm (p)}(T_0,\lambda)A_{\rm d}^{(0)}(\lambda)\,d\lambda\,, 
\\\\
i\D_{T_0}A_{\rm d}^{(2)}(\eta)&-&2\eta^2A_{\rm d}^{(2)}(\eta)&=&\,\,\dots\,,
\end{array}\right.
\end{eqnarray}
and so on.  Our initial conditions are encoded in the expansions
(\ref{eq:expansion}) as $A_{\rm b}^{(0)}(0,0,0,\dots)=A_{{\rm b}0}$,
$A_{\rm d}^{(0)}(0,0,0,\dots,\lambda)=0$,
and for $j\ge 1$, $A_{\rm b}^{(j)}(0,0,0,\dots)= A_{\rm
d}^{(j)}(0,0,0,\dots,\lambda)=0$.  We now proceed to solve the
hierarchy sequentially.

Solving the equations (\ref{eq:order0}) at order $O(1)$ subject to the
initial conditions gives
\begin{equation}
A_{\rm b}^{(0)}=Ce^{2i\beta_{\rm b}T_0}\,,\hspace{0.3 in}
A_{\rm d}^{(0)}(\eta)=0\,,
\end{equation}
where $C=C(T_1,T_2,\dots)$ satisfies the initial condition
$C(0,0,\dots)=A_{{\rm b}0}$ but is otherwise undetermined at this
stage.

In the first of the two equations (\ref{eq:order1}) appearing at $O(\epsilon)$,
it is natural to make the substitution
\begin{equation}
A_{\rm b}^{(1)}=fe^{2i\beta_{\rm b}T_0}\,,
\end{equation}
which leads to the equation
\begin{equation}
\D_{T_0}f=-\D_{T_1}C-iM_1(T_0)C\,.
\end{equation}
Integrating with the use of the initial condition $f(T_0=0)=0$, and keeping
in mind that $T_1$ and $T_0$ are to be thought of as independent variables,
leads to the expression for $A_{\rm b}^{(1)}$:
\begin{equation}
A_{\rm b}^{(1)}=\left(-\D_{T_1}C\cdot T_0-iC\int_0^{T_0}M_1(s)\,ds\right)
e^{2i\beta_{\rm b}T_0}\,.
\end{equation}
We need for this correction to be bounded as a function of $T_0$ so
that the asymptotic expansion will be well-ordered for long times.
Since $M_1(s)$ is a periodic function of period $L$, this requirement
uniquely determines $\D_{T_1}C$:
\begin{equation}
\D_{T_1}C=-i\overline{M_1} C\,,\hspace{0.2 in}
\overline{M_1} \doteq \frac{1}{L}\int_0^L M_1(s)\,ds\,.
\end{equation}
Thus,
\begin{equation}
C=De^{-i\overline{M_1} T_1}\,,\hspace{0.2 in}D=D(T_2,\dots)\,,
\hspace{0.2 in}D(0,\dots)=A_{{\rm b}0}\,.
\end{equation}
Putting together what we have for the bound state amplitude at this time,
\begin{equation}
A_{\rm b}^{(0)}=De^{-i\overline{M_1} T_1}e^{2i\beta_{\rm b}T_0}\,,
\hspace{0.2 in}
A_{\rm b}^{(1)}=-iDe^{-i\overline{M_1} T_1}e^{2i\beta_{\rm b}T_0}
\int_0^{T_0}\left(M_1(s)-\overline{M_1}\right)\,ds\,.
\end{equation}
This has been the first application in our calculation of the guiding
principle of the method of multiple scales, that dependence of
expansion terms on ``slow'' times is chosen to ensure that the
expansion terms are uniformly bounded with respect to the ``fast''
time $T_0$.  Now we solve for the correction to the dispersive mode
amplitude at this order (in fact the leading term) using the second of
the equations (\ref{eq:order1}).  Substituting the expressions from
the previous order and using the initial conditions gives a unique
expression:
\begin{equation}
A_{\rm d}^{(1)}(\eta)=-iDe^{-i\overline{M_1} T_1}
\int_0^{T_0}N_1^{\rm (p)}(s,\eta)^*e^{2i\beta_{\rm b}s}e^{-2i\eta^2(T_0-s)}
\,ds\,.
\end{equation}

Continuing systematically with the equation (\ref{eq:order2}) for the
bound state amplitude correction at order $O(\epsilon^2)$, we substitute
all the expressions known thus far and observe the utility of the change of
variables
\begin{equation}
A_{\rm b}^{(2)}=he^{-i\overline{M_1} T_1}e^{2i\beta_{\rm b}T_0}\,.
\end{equation}
We find for $h$ the simple equation
\begin{equation}
\label{eq:heqn}
\D_{T_0}h=-D\cdot\left(M_1(T_0)-\overline{M_1}\right)\int_0^{T_0}
\left(M_1(s)-\overline{M_1}\right)\,ds-\D_{T_2}D-D\cdot\gamma(T_0)\,,
\end{equation}
where 
\be
\gamma(T_0)\doteq iM_2(T_0)+\int_0^\infty N_1^{\rm (p)}(T_0,\lambda)
\int_0^{T_0}e^{-2i(\lambda^2 +\beta_{\rm b})(T_0-s)}
N_1^{\rm (p)}(s,\lambda)^*\,ds\,d\lambda\,.\label{eq:gammaeqn}
\ee

Equation (\ref{eq:heqn}) can be analyzed as follows. By linearity, we
can express $h$ as a sum: $h=h_1+h_2$, where
\begin{equation}
\begin{array}{rcl}
\D_{T_0}h_1&=&\displaystyle -D\cdot\left(M_1(T_0)-\overline{M_1}\right)
\int_0^{T_0}\left(M_1(s)-\overline{M_1}\right)\,ds\,,\\\\
\D_{T_0}h_2&=&-\D_{T_2}D-D\cdot\gamma(T_0)\,,
\end{array}
\label{eq:heqns}
\end{equation}
and where we assume the initial conditions $h_1(T_0=0)=h_2(T_0=0)=0$.
Integrating the equation for $h_1$ exactly using the 
initial condition gives:
\begin{equation}
h_1= -{D\over2}\left(\int_0^{T_0}(M_1(s)-\overline{M_1})\, ds
\right)^2\,,
\label{eq:Ab12}
\end{equation}
which is periodic, and in particular bounded, by periodicity of
$M_1(T_0)$.

We now want to select the dependence of $D$ on the slow time $T_2$
such that $h_2$, as found from the second of equations
(\ref{eq:heqns}), is a bounded function of $T_0$.  Using the initial
condition to integrate the equation for $h_2$ with respect to $T_0$
while holding $T_2$ fixed gives
\begin{equation}
h_2=-T_0\partial_{T_2}D - D\int_0^{T_0}\gamma(s)\,ds\,.
\label{eq:formulah}
\end{equation}
Clearly, the possibility of choosing $D(T_2)$ so that
the expression (\ref{eq:formulah}) is bounded in $T_0$ depends on the
behavior of $\gamma(T_0)$ in the limits $T_0\rightarrow \pm\infty$.
We now study $\gamma(T_0)$ for large $|T_0|$.

We can compute the $s$-integral in (\ref{eq:gammaeqn}) exactly if we
introduce the Fourier series for the periodic function $N_1^{\rm
(p)}(T_0,\lambda)$:
\begin{equation}
N_1^{\rm (p)}(T_0,\lambda)=\sum_{k=-\infty}^\infty
N_{1,k}(\lambda)e^{2\pi i kT_0/L}\,.
\end{equation}
Note that in terms of the Fourier coefficients of $N^{\rm
(p)}(T_0,\lambda;\epsilon)$ itself (see (\ref{eq:MNKexpansions})), we
have
\begin{equation}
N_{1,k}(\lambda)=\lim_{\epsilon\rightarrow
0}\epsilon^{-1}N_k(\lambda;\epsilon)\,.
\end{equation}
Substituting the Fourier series into (\ref{eq:gammaeqn}), integrating
term by term with respect to $s$, and changing variables to
$\sigma=\lambda^2$, we arrive at
\begin{equation}
\gamma(T_0)=iM_2(T_0)+\sum_{n,k=-\infty}^\infty i\int_0^\infty
\frac{N_{1,n}(\sqrt{\sigma})N_{1,k}(\sqrt{\sigma})^*}{4\sqrt{\sigma}
(\sigma-\sigma_k)}\left[e^{-2i(\sigma-\sigma_n)T_0}-
e^{2\pi i(n-k)T_0/L}\right]\,d\sigma\,,
\end{equation}
where the resonances $\sigma_n$ are defined by 
\begin{equation}
\sigma_n\doteq \pi n/L - \beta_{\rm b}\,.
\label{eq:resonances}
\end{equation}
Note that for all terms having $\sigma_k>0$, the difference of the
exponentials in the integrand vanishes for $\sigma=\sigma_k$, so there
is no singularity.  Moreover, the Fourier coefficients
$N_{1,n}(\lambda)$ are by construction analytic functions of $\lambda$
for $\lambda$ in a sector including the real axis, and so the
quantities $N_{1,n}(\sqrt{\sigma})$ are analytic in a neighborhood of
the positive real $\sigma$ axis.  This property extends to the whole
integrand, and we may therefore deform the integration contour away
from the real axis in an effort to study the behavior for large
$|T_0|$ by a steepest descents type argument.

\begin{figure}[h]
\begin{center}
\mbox{\psfig{file=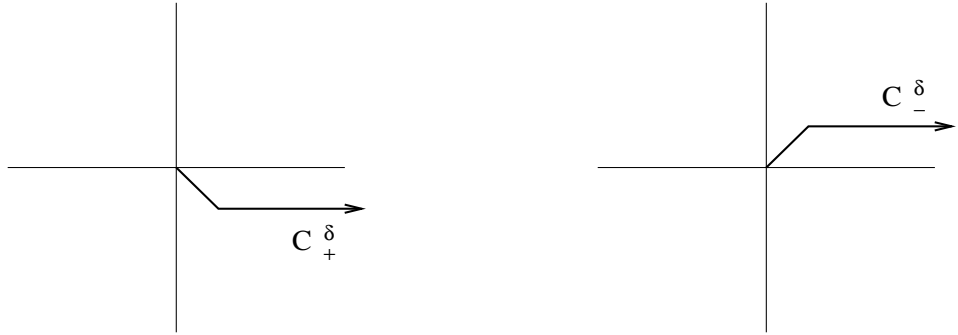,width=5 in}}
\end{center}
\caption{\em The deformed integration contours $C^\delta_+$ and $C^\delta_-$.}
\label{fig:cdelta}
\end{figure}

For positive $T_0$, we deform the contour into the lower half plane.
For $\delta>0$, let $C_+^\delta $ be the contour consisting of the
diagonal segment from $0$ to $(1-i)\delta$ followed by the horizontal
ray from $(1-i)\delta$ to $-i\delta +\infty$ (see
Figure~\ref{fig:cdelta}).  
We have
\begin{equation}
\begin{array}{rcl}
\gamma(T_0)&=&\displaystyle iM_2(T_0)+\sum_{n,k=-\infty}^\infty i\int_{C_+^\delta} 
\frac{N_{1,n}(\sqrt{\sigma})N_{1,k}(\sqrt{\sigma})^*}
{4\sqrt{\sigma}(\sigma-\sigma_k)}\left[e^{-2i(\sigma-\sigma_n)T_0}-
e^{2\pi i (n-k)T_0/L}\right]\,d\sigma\\\\
&=&\displaystyle
iM_2(T_0)+\sum_{n,k=-\infty}^\infty i\int_{C_+^\delta}\frac{N_{1,n}(\sqrt{\sigma})
N_{1,k}(\sqrt{\sigma})^*}{4\sqrt{\sigma}(\sigma-
\sigma_k)}e^{-2i(\sigma-\sigma_n)T_0}\,d\sigma\\\\
&&\displaystyle\hspace{0.2 in} -\,\,
\sum_{n,k=-\infty}^\infty i\int_{C_+^\delta}\frac{N_{1,n}(\sqrt{\sigma})
N_{1,k}(\sqrt{\sigma})^*}{4\sqrt{\sigma}(\sigma-
\sigma_k)}e^{2\pi i(n-k)T_0/L}\,d\sigma\\\\
&=&\gamma_0(T_0)+\gamma_1^+(T_0) + \gamma_2^+(T_0)\,,
\end{array}
\end{equation}
so that on the new contour $C_+^\delta$ the two integrals converge
independently.  The term $\gamma_0(T_0)$ is periodic in $T_0$ with
period $L$ and mean value $i\overline{M_2}$.  The term
$\gamma_2^+(T_0)$ is also a periodic function of $T_0$ with period
$L$.  Its mean value is given by the terms in the sum with $n=k$:
\begin{equation}
\overline{\gamma_2^+} = -\sum_{n=-\infty}^\infty i\int_{C_+^\delta}
\frac{|N_{1,n}(\sqrt{\sigma})|^2}{4\sqrt{\sigma}(\sigma-\sigma_n)}\,d\sigma\,.
\end{equation}
Letting $\delta$ tend to zero does not alter the value of the integral, and
then we may use the Plemelj-Sokhotski formula\footnote{This is merely the distributional identity 
\begin{displaymath}
(x\pm i0)^{-1}={\rm P.\,V.\,}x^{-1}\mp i\pi\delta(x)\,.
\end{displaymath}
} to evaluate the terms with $\sigma_n>0$ to find
\begin{equation}
\overline{\gamma_0+\gamma_2^+} = i\overline{M_2} -i\Lambda_2 +\Gamma_2\,,
\end{equation}
where
\begin{equation}
\Lambda_2 \doteq  \sum_{n=-\infty}^{n_0-1}\int_0^\infty \frac{|N_{1,n}(\sqrt{\sigma})|^2\,d\sigma}{4\sqrt{\sigma}(\sigma-\sigma_n)} +\sum_{n=n_0}^\infty
{\rm P.\,V.\,}\int_0^\infty\frac{|N_{1,n}(\sqrt{\sigma})|^2\,d\sigma}{4\sqrt{\sigma}(\sigma-\sigma_n)}\,,
\label{eq:Lambshift}
\end{equation}
and
\begin{equation}
\Gamma_2\doteq\frac{\pi}{4}\sum_{n=n_0}^\infty \frac{|N_{1,n}(\sqrt{\sigma_n})|^2}{\sqrt{\sigma_n}}\,.
\label{eq:gamma2}
\end{equation}
Finally, consider the term $\gamma_1^+(T_0)$.
Its time integral, calculated term by term, is
\begin{equation}
\begin{array}{rcl}
\displaystyle \int_0^{T_0}\gamma_1^+(s)\,ds&=&\displaystyle
-\sum_{n,k=-\infty}^\infty
\int_{C_+^\delta} \frac{N_{1,n}(\sqrt{\sigma})N_{1,k}(\sqrt{\sigma})^*}
{8\sqrt{\sigma}(\sigma-\sigma_k)(\sigma-\sigma_n)}\left[e^{-2i(\sigma-\sigma_n)T_0}-1\right]\,d\sigma\\\\
&=&\displaystyle
\sum_{n,k=-\infty}^\infty
\int_{C_+^\delta} \frac{N_{1,n}(\sqrt{\sigma})N_{1,k}(\sqrt{\sigma})^*}
{8\sqrt{\sigma}(\sigma-\sigma_k)(\sigma-\sigma_n)}\,d\sigma \\\\
&&\displaystyle\hspace{0.2 in}-\,\,
\sum_{n,k=-\infty}^\infty
\int_{C_+^\delta} \frac{N_{1,n}(\sqrt{\sigma})N_{1,k}(\sqrt{\sigma})^*
}
{8\sqrt{\sigma}(\sigma-\sigma_k)(\sigma-\sigma_n)}e^{-2i(\sigma-\sigma_n)T_0}\,d\sigma\,.
\end{array}
\end{equation}
The first term is independent of $T_0$ (and also of $\delta>0$, since the integral converges and the integrand is analytic).  In the second term,
the real part of the exponent is negative for
$T_0>0$, so for $T_0$ large and positive, the integrand is
exponentially small except in a small neighborhood of $\sigma=0$.
This small neighborhood gives a leading contribution to the integrand
that is $O(T_0^{-1/2})$, and in particular is bounded for large $T_0>0$.  

Putting these results together, we find that for large $T_0>0$, we have
\begin{equation}
\int_0^{T_0}\gamma(s)\,ds = (i\overline{M_2} -i\Lambda_2 + \Gamma_2)T_0 + O(1)\,.
\end{equation}
Going back to (\ref{eq:formulah}), it is clear that choosing
\begin{equation}
\partial_{T_2}D=-(i\overline{M_2} -i\Lambda_2 +\Gamma_2)D\,,
\end{equation}
will lead to a solution $h_2(T_0)$ that is uniformly bounded
for all $T_0>0$.  Also note that the first term in $\gamma(T_0)$
contributes a term
\begin{equation}
h_{2,M}\doteq-i\int_0^{T_0}\left(M_2(s)-\overline{M_2}\right)\,ds\,,
\end{equation}
to the expression for $h_2(T_0)$.  We write
$h_2(T_0)=h_{2,M}(T_0)+\tilde{h}_2(T_0)$.

To find the behavior of $\gamma(T_0)$ and its time integral as $T_0$
tends to $-\infty$, we repeat the above steps, this time deforming the
integration contour into the upper half plane to facilitate the
steepest descents argument.  The path of integration is now
$C^\delta_-$ (see Figure~\ref{fig:cdelta}).  The only difference is in
the sign of $\Gamma_2$; the correct choice for a bounded solution for
all $T_0$ is
\begin{equation}
\partial_{T_2}D=-(i\overline{M_2} -i\Lambda_2 +{\rm sgn}(T_0)\Gamma_2)D\,.
\end{equation}

Thus, the method of multiple scales gives the following approximation to
the bound state mode amplitude:
\begin{equation}
\begin{array}{l}
\displaystyle
A_{\rm b}(t)=A_{{\rm b}0}e^{2i\beta_{\rm b}t}
e^{-i(\epsilon\overline{M_1}+\epsilon^2\overline{M_2})t}
e^{i\epsilon^2\Lambda_2 t}
e^{-\epsilon^2\Gamma_2|t|}
\Bigg(1-i\epsilon\int_0^t\left(M_1(s)-\overline{M_1}\right)\,ds \\\\
\displaystyle
\hspace{0.2 in}
-\,\,
\frac{\epsilon^2}{2}\left(\int_0^t\left(M_1(s)-\overline{M_1}\right)\,ds
\right)^2 -i\epsilon^2\int_0^t\left(M_2(s)-\overline{M_2}\right)\,ds
+ \epsilon^2 \tilde{h}_2(t)+O(\epsilon^3)\Bigg)\,.
\end{array}
\end{equation}
It is not hard to see that an asymptotically equivalent expression is just
\begin{equation}
A_{\rm b}(t)=A_{{\rm b}0}e^{2i\beta_{\rm b}t}e^{-\epsilon^2\Gamma_2|t|}e^{i\epsilon^2\Lambda_2
t}e^{-i\int_0^tM(s)\,ds}\left(1+O(\epsilon^2)\right)\,.
\label{eq:asympformula}
\end{equation}
This asymptotic formula is expected to be uniformly valid as
$\epsilon$ tends to zero for all $|t|<K\epsilon^{-2}$ for any constant
$K$.

So the behavior of the bound state amplitude under the influence of a
periodic perturbation, as predicted by the multiple scale theory, is
dominated by two effects, a shift in frequency accompanied by
exponential decay.  The shift in frequency is an order $O(\epsilon)$
effect, coming from $\overline{M}$.  This shift can be traced back to
the influence of the perturbation directly on the bound state; there
is no coupling to any other modes in this term.  The order
$O(\epsilon^2)$ effects include both a further adjustment to the
frequency through the quantity $\epsilon^2\Lambda_2$, the {\em Lamb
shift}, and exponential decay through the quantity
$\epsilon^2\Gamma_2$.  Clearly these two numbers are the real and
imaginary parts of the same complex frequency.  Unlike the leading
order phase shift, both of these effects are clearly due to the
resonant coupling between the bound state and the continuum that is
introduced and mediated by the periodic perturbation.  Due to the
exponential decay, the lifetime of the bound state is seen to be
approximately $\epsilon^{-2}/\Gamma_2$, which is quite long for small
$\epsilon$.  For this reason, under small perturbations of the
potential energy the state is called {\em metastable}.

\noindent{\bf Remark:} The validity of this expansion procedure is
clearly called into question if any of the resonances $\sigma_n$ are
very close to zero, in which case the complex frequency
$\Lambda_2+i\Gamma_2$ is potentially large.  The breakdown of the
expansion in this case indicates the presence of a {\em parametric
zero energy resonance}.  Note, however, that in the odd case the
matrix element $N^{\rm (p)}(t,\lambda)$ vanishes as $\lambda$ tends to
zero, and therefore so do the corresponding Fourier coefficients (and
in particular they vanish to leading order in $\epsilon$, that is,
$N_{1,n}(\lambda)$ vanishes at $\lambda=0$ for all $n$).  This
suggests that the expansion (\ref{eq:asympformula}) continues to hold
in the odd case as the parameters $\rho_1$ and $\rho_2$ of the
two-soliton potential are varied so as to cause a resonance
$\sigma_n(\rho_1,\rho_2)$ to pass through zero.  In the even case,
however, behavior possibly very different from that predicted by the
formula (\ref{eq:asympformula}) is expected if a resonance is close to
zero.  We plan to investigate this phenomenon analytically; however in
this paper we will demonstrate the effects of parametric zero energy
resonance in both the even and odd cases with numerical simulations.
Sudden changes in the behavior of a simple model for atomic ionization
as a parameter is smoothly varied, causing the system to pass through
a zero energy resonance, have recently been observed and compared with
experiment by Costin, Lebowitz and Rokhlenko \cite{CLR99}. \pfend

\subsection{Rigorous analysis and infinite time results.}
The multiple-scale analysis of the preceding section leads to an
asymptotic formula for the decaying bound state amplitude that is  
valid on time intervals of order $\epsilon^{-2}$.  In this section,
we will establish the validity of the asymptotic formula
(\ref{eq:asympformula}) in certain circumstances using the results of 
Kirr and Weinstein \cite{KW}.  When applicable, these results also yield
a detailed description of the solution as $t\rightarrow\pm\infty$.

More precisely, we now study the perturbed periodic system in the form
obtained by use of Floquet factorization of the time-periodic
unperturbed Hamiltonian $\op{H}_0(t)$:
\begin{equation}
i\partial_ty - \op{B}y =\tilde{\op{W}}(t)y\,.
\label{eq:KWstartingpoint}
\end{equation}
The self-adjoint operator $\op{B}:L^2(\R)\rightarrow L^2(\R)$ defined
in \S \ref{sec:periodic} can be thought of as a time-independent
Hamiltonian, and the idea is to apply the theory of periodic (or
almost periodic) perturbations of autonomous linear Hamiltonian
systems as developed in \cite{KW} and \cite{ionization} directly to
the problem in this form.

As we did in the multiple scales analysis, we will restrict attention
to the special case of periodically perturbed even two-soliton
periodic potentials.  As we know, in this case the operator $\op{B}$
has exactly two $L^2$ eigenfunctions, one an even function of $x$ and
the other an odd function of $x$.  Since $L^2$ is the direct sum of
its two subspaces $L^2_{\rm (e,o)}$ of even and odd functions, and
since $\op{B}$ leaves each subspace invariant, we may study the
problem (\ref{eq:KWstartingpoint}) restricted to one subspace at a
time.  This reduction results in an unperturbed problem with a single
bound state, and to such problems the results described in \cite{KW}
and \cite{ionization} can be applied without modification.

On each subspace $L^2_{\rm (e,o)}(\R)$, the operator $\op{B}$ is
explicitly given by:
\begin{equation}
\op{B}f(x)=\langle\Psi_{\rm b}^{\rm (e,o)}(\cdot,0),-2\beta_{\rm b}f(\cdot)
\rangle\Psi_{\rm b}^{\rm (e,o)}(x,0) +\int_0^\infty\langle\Psi_{\rm d}^{\rm (e,o)}(\cdot,0,\lambda),2\lambda^2f(\cdot)\rangle\Psi_{\rm d}^{\rm (e,o)}(x,0,\lambda)\,d\lambda\,,
\end{equation}
where the functions $\Psi_{\rm b}^{\rm (e,o)}(x,t)$ and $\Psi_{\rm
d}^{\rm (e,o)}(x,t,\lambda)$ are defined given the parameters $\rho_1$
and $\rho_2$ in \S \ref{sec:separable}.  The hypotheses required in
\cite{KW} of the even and odd restrictions of the operator $\op{B}$ are
reproduced here adapted to our application:
\begin{itemize}
\item[\bf (H1)] The even and odd restrictions of $\op{B}$ are densely
defined on subspaces of $L^2_{\rm (e,o)}(\R)$ and have self-adjoint
extensions to all of $L^2_{\rm (e,o)}(\R)$.
\item[\bf (H2)] The spectrum of $\op{B}$ in each of $L^2_{\rm
(e,o)}(\R)$ consists of an absolutely continuous part $\sigma^{\rm
(e,o)}_{\rm cont}(\op{B})=[0,\infty]$ with associated spectral
projection $\op{P}_{\rm c}^{\rm (e,o)}$ and a single isolated
eigenvalue $\lambda_0=-2\beta_{\rm b}$ with corresponding normalized
eigenstate $\psi_0(x)=\Psi_{\rm b}^{\rm (e,o)}(x,0)$, so that
\begin{equation}
\op{B}\psi_0 = \lambda_0\psi_0\,,\hspace{0.3 in}\|\psi_0\|_2=1\,.
\end{equation}
\item[\bf (H3)] The {\em odd} restriction of $\op{B}$ satisfies two
{\em dispersive local decay estimates}.  There exist constants $C_{\rm
ns}$ and $C_{\rm s}$ such that
\begin{itemize}
\item[\bf (a)] The nonsingular local decay estimate
\begin{equation}
\|\langle \cdot\rangle^{-7/2}e^{-i\op{B}t}\op{P}_{\rm c}^{\rm (o)}f\|_2
\le C_{\rm ns}\langle t\rangle^{-3/2}\|\langle\cdot\rangle^{7/2}f\|_2\,,
\label{eq:nsld}
\end{equation}
holds for all $f\in L^2_{\rm (o)}(\R)$.
\item[\bf (b)] The singular local decay estimate
\begin{equation}
\|\langle\cdot\rangle^{-7/2}e^{-i\op{B}t}(\op{B}-2\mu-2i\kappa 0)^{-1}
\op{P}_{\rm c}^{\rm (o)}f\|_2\le C_{\rm s}\langle t\rangle^{-3/2}\|\langle\cdot
\rangle^{7/2}f\|_2\,,
\label{eq:sld}
\end{equation}
where $\kappa={\rm sgn}\,(t)$, holds uniformly for all $\mu$ satisfying
$|\mu|\ge\mu_{\rm min}>0$, that is, the constant $C_{\rm s}$ only
depends on $\mu_{\rm min}$.
\end{itemize}
\end{itemize}
The local decay hypotheses are established in \ref{app:localdecay}.
We remark here that due to a {\em zero energy resonance}, the decay
estimates that are established in \ref{app:localdecay} for the even
case are of the form (\ref{eq:nsld}) and (\ref{eq:sld}) but with decay
rate $\langle t\rangle^{-1/2}$ (this is a sharp estimate).
Unfortunately, this slower rate of decay precludes the direct
application of the results in \cite{KW} and \cite{ionization}.  On the
other hand, as long as the perturbation does not create a resonance
$\mu$ that is close to zero, we can expect similar results to hold in
the even case over time scales of length $|t|<K/\epsilon^2$, since
there is no obvious difficulty with the multiple scales analysis.

The application of the results of \cite{KW} and \cite{ionization} also
requires some hypotheses to be satisfied by the perturbation operator
$\tilde{\op{W}}(t)$ and its relation to the unperturbed Hamiltonian
$\op{B}$.  The perturbation operator acting on $L^2_{\rm (e,o)}(\R)$
takes the form
\begin{equation}
\begin{array}{rcl}
\displaystyle
\tilde{\op{W}}(t)f(x)&=&
\displaystyle
\Bigg(M(t)\langle \Psi_{\rm b}^{\rm (e,o)}(\cdot,0),f(\cdot)
\rangle \\\\
&&\displaystyle \hspace{0.4 in}
+\,\,
\int_0^\infty N^{\rm (p)}(t,\eta)\langle
\Psi_{\rm d}^{\rm (e,o)}(\cdot,0,\eta),f(\cdot)\rangle\,d\eta\Bigg)
\Psi_{\rm b}^{\rm (e,o)}(x,0)\\\\
&&\displaystyle\hspace{0.2 in}+\,\,
\int_0^\infty \Bigg(N^{\rm (p)}(t,\lambda)^*\langle \Psi_{\rm b}^{\rm (e,o)}
(\cdot,0),f(\cdot)\rangle \\\\
&&\displaystyle\hspace{0.4 in}+\,\,\int_0^\infty K^{\rm (p)}(t,\lambda,\eta)
\langle \Psi_{\rm d}^{\rm (e,o)}(\cdot,0,\eta),f(\cdot)\rangle\,d\eta
\Bigg)\Psi_{\rm d}^{\rm (e,o)}(x,0,\lambda)\,d\lambda\,.
\end{array}
\label{eq:oncedressedpert}
\end{equation}
We recall that the periodic ``matrix elements'' in the above
expression are defined in terms of either the odd or even modes by
(\ref{eq:MNKs}).  This operator, being periodic in $t$ with period
$L$, has a Fourier series expansion
\begin{equation}
\tilde{\op{W}}(t)=\sum_{k=-\infty}^\infty e^{2\pi ikt/L} \tilde{\op{W}}_k\,,
\label{eq:pertFourier}
\end{equation}
where each operator $\tilde{\op{W}}_k$ has the same form as
(\ref{eq:oncedressedpert}) with the functions $M(t)$, $N^{\rm
(p)}(t,\lambda)$, and $K^{\rm (p)}(t,\lambda,\eta)$ replaced by the
corresponding $k$th Fourier coefficients, $M_k$, $N_k(\lambda)$ and
$K_k(\lambda,\eta)$ respectively.  The hypotheses required in \cite{KW}
of the perturbation adapted to this context are:
\begin{itemize}
\item[\bf (H4)]  The operators $\tilde{\op{W}}_k$ satisfy $\tilde{\op{W}}_{-k}=
\tilde{\op{W}}_k^\dagger$ and 
\begin{equation}
\sum_{k=-\infty}^\infty \|\tilde{\op{W}}_k\|_{{\cal L}(L^2(\R))} <\infty\,,
\end{equation}
where $\|\cdot\|_{{\cal L}(L^2(\R))}$ denotes the uniform operator norm in
$L^2(\R)$.  Also,
\begin{equation}
|\|\tilde{\op{W}}(\cdot)\||\doteq
\sum_{k=-\infty}^\infty \left(\|\langle\cdot\rangle^{7/2} \tilde{\op{W}}_k
\|_{{\cal L}(L^2(\R))}+\|\langle\cdot\rangle^{7/2} \tilde{\op{W}}_k
\langle\cdot\rangle^{7/2}\|_{{\cal L}(L^2(\R))}\right)<\infty\,.
\end{equation}
\item[\bf (H5)]  The following {\em resonance condition} holds:
\begin{equation}
\Gamma\doteq\pi\sum_{n=n_0}^\infty
\langle \tilde{\op{W}}_n\psi_0,\delta(\op{B}-2\sigma_n)\tilde{\op{W}}_n
\psi_0\rangle >0\,,
\label{eq:abstractgamma}
\end{equation}
where the resonances are defined by $\sigma_n=(\lambda_0+2\pi
n/L)/2=-\beta_{\rm b}+n\pi/L$.  Here $n_0$ is the smallest positive
integer for which $\sigma_{n_0}>0$.  
%Note that in \cite{KW} a stronger
%version of this condition is hypothesized, namely that $\Gamma/
%|\|\tilde{\op{W}}(\cdot)\||^2$ is uniformly bounded away from zero for
%perturbations in a set of interest.  This stronger version is only
%important if constants that appear in subsequent estimates are desired
%to be uniform as well.  One has this uniformity automatically if one
%considers corrections to the potential energy of the form
%$W(x,t)=\epsilon W_1(x,t)$ for some $\epsilon\in\R$ and some given
%function $W_1(x,t)$.  
Note that since
\begin{equation}
\tilde{\op{W}}_n\psi_0=\tilde{\op{W}}_n\Psi_{\rm b}^{\rm (e,o)}(x,0)=
M_n\Psi_{\rm b}^{\rm (e,o)}(x,0)+\int_0^\infty N_n(\lambda)^*
\Psi_{\rm d}^{\rm (e,o)}(x,0,\lambda)\,d\lambda\,,
\end{equation}
and since for $\sigma>0$
\begin{equation}
\langle f(\cdot),\delta(\op{B}-2\sigma)f(\cdot)\rangle = \int_0^\infty
|(\op{V}\op{P}_{\rm c}^{\rm
(e,o)}f)(\lambda)|^2\delta(2\lambda^2-2\sigma)\,d\lambda =
\frac{|(\op{V}\op{P}_{\rm c}^{\rm
(e,o)}f)(\sqrt{\sigma})|^2}{4\sqrt{\sigma}}\,,
\end{equation}
The formula for $\Gamma$ can be written as
\begin{equation}
\Gamma=\frac{\pi}{4}\sum_{n=n_0}^\infty \frac{N_n(\sqrt{\sigma_n})}{\sqrt{\sigma_n}}\,.
\end{equation}
Note that if $N_n(\lambda)$ has an expansion in a small parameter
$\epsilon$, with leading term linear in $\epsilon$ as was assumed in
the multiple scale analysis, then the leading term of the
corresponding expansion for $\Gamma(\epsilon)$ is exactly
$\epsilon^2\Gamma_2$, where $\Gamma_2$ is correctly obtained by the
multiple scale analysis and is given by (\ref{eq:gamma2}).  The
constant $\Gamma$ is a decay rate associated with the bound state of
the unperturbed system.  The statement that the expression
(\ref{eq:abstractgamma}) should be positive for decay to occur as a
consequence of resonant coupling to the continuum is attributed to
Fermi and is known as the ``Fermi Golden Rule''.  Again, because the
decay constant $\Gamma$ is quadratic in the size of the perturbation,
the exponential decay process is very slow for small perturbations.
Thus, in the presence of a small perturbation $W(x,t)$, the bound
state is said to be metastable.
\item[\bf (H6)] There are no finite accumulation points of the
resonances $\sigma_n$, $n\ge n_0$.  This is satisfied automatically
because the Fourier expansion of $\tilde{\op{W}}(t)$ is that of a
periodic function.  The point here is that the results in \cite{KW}
are more general; for example this hypothesis is satisfied by finite
Fourier sums with incommensurate frequencies.  Yet further
generalizations can be found in \cite{KW}.
\end{itemize}
Verifying the hypothesis {\bf (H4)} would seem to require more
detailed information about the correction to the potential energy $W(x,t)$
than we have used thus far.  We merely point out at this time that
by elementary Cauchy-Schwarz arguments applied to the unitarily equivalent
operators $\op{V}\tilde{\op{W}}_n\op{V}^\dagger$, one finds the estimate
\begin{equation}
\|\tilde{\op{W}}_n\|_{{\cal L}(L^2(\R))}\le 2\sqrt{|M_k|^2 +
2\int_0^\infty |N_n(\lambda)|^2\,d\lambda +\int_0^\infty\int_0^\infty
|K_n(\lambda,\eta)|^2\,d\lambda\,d\eta}\,.
\end{equation}
Assuming these bounds all are finite, which is really a question of
the smoothness and decay of ``snapshots'' of the function $W(x,t)$ at
fixed $t$, we see that the first required bound in {\bf (H4)} will be
satisfied if the Fourier coefficients of the function $W(x,t)$ in $t$
decay faster than, say, $1/n$.  This is because the other periodic
contributions come from the analytic eigenfunctions, whose Fourier
coefficients decay faster than $1/n^p$ for any $p>0$.  Therefore, not much
beyond continuity in $t$ is required of $W(x,t)$, at least for this simpler
estimate.  More restrictions are certainly required to satisfy the second
estimate of {\bf (H4)}.

These hypotheses imply the following results:
\begin{prop}[\bf Theorem 2.1 of \cite{KW}]
Let $\op{B}$ and $\tilde{\op{W}}(t)$ satisfy the above hypotheses and
let an odd function $y_0(x)$ be given such that $\langle
x\rangle^{7/2}y_0(x)\in L^2_{\rm (o)}(\R)$.  Let $y(x,t)$ be the
solution of (\ref{eq:KWstartingpoint}) with initial condition
$y(x,0)=y_0(x)$.  Then if $|\|\tilde{\op{W}}(\cdot) \||$ is
sufficiently small there exists a constant $C$ such that
\begin{equation}
\|\langle \cdot\rangle^{-7/2} y(\cdot,t)\|_2\le C
\langle t\rangle^{-3/2}\|\langle\cdot\rangle^{7/2}y_0(\cdot)\|_2\,,
\end{equation}
holds for all $t\in\R$.
\label{prop:general}
\end{prop}
\begin{prop}[\bf Theorem 2.2 of \cite{KW}]
Assume the same hypotheses of $\op{B}$ and $\tilde{\op{W}}(t)$.  Then
if $|\|\tilde{\op{W}}(\cdot)\||$ is sufficiently small, the solution
$y(x,t)$ of (\ref{eq:KWstartingpoint}) corresponding to the odd
initial condition $y_0(x)$ with $\langle x\rangle^{7/2}y_0(x)\in
L^2_{\rm (o)}(\R)$ is of the form
\begin{equation}
\begin{array}{rcl}
y(x,t) &=&\displaystyle \left[\langle \Psi_{\rm b}^{\rm
(o)}(\cdot,0),y_0(\cdot)\rangle e^{2i\beta_{\rm b}t}e^{-\Gamma
|t|}e^{i\Lambda t}e^{-i\int_0^t
M(s)\,ds}e^{ir_1(t)}+r_2(t)\right]\Psi_{\rm b}^{\rm (o)}(x,0) \\\\
&&\displaystyle\hspace{0.2 in}+\,\, \left( e^{-it\op{B}}\op{P}_{\rm
c}^{\rm (o)}y_0(\cdot)\right)(x,t) + \tilde{y}(x,t)\,,
\end{array}
\end{equation}
where 
\begin{equation}
\Lambda\doteq \sum_{n=-\infty}^\infty \langle \tilde{\op{W}}_n\psi_0,
{\rm P.\,V.\,}(\op{B}-2\sigma_n)^{-1}\op{P}_{\rm c}^{\rm (o)}\tilde{\op{W}}_n
\psi_0\rangle\,,
\end{equation}
and where
\begin{itemize}
\item
The phase correction $r_1(t)$ is uniformly bounded and
$O(|\|\tilde{\op{W}}(\cdot)\||^2)$\,,
\item
The bound state amplitude error $r_2(t)$ is $O(|\|\tilde{\op{W}}(\cdot)
\||)$ uniformly for all $|t|<K/\Gamma$ for all fixed $K$ and decays for
large time as $O(\langle t\rangle^{-3/2})$,
\item
The correction $\tilde{y}(x,t)$ is orthogonal to the bound state: $\langle
\psi_0,\tilde{y}(\cdot,t)\rangle=0$ for all $t$, and satisfies the dispersive
decay estimate $\|\langle\cdot\rangle^{-7/2}\tilde{y}(\cdot,t)\|_2=O(\langle
t\rangle^{-3/2})$ for large $t$.
\end{itemize}
\label{prop:detail}
\end{prop}

\noindent{\bf Remark:} Propositions~\ref{prop:general} and
\ref{prop:detail} would appear to say that all initial conditions
decay exponentially and then algebraically.  However, a more careful
reading shows that it is possible for there to be a transient stage of
growth, before the decay ultimately sets in.  This is because the
error terms, although small when the perturbation is small, are not
{\em uniformly} small for all initial conditions $y_0(x)$ such that
$\langle x\rangle^{7/2} y_0(x)$ ranges over the unit sphere in
$L^2(\R)$.  So, for each fixed perturbation $W(x,t)$, no matter how
small, it is possible to find an initial condition $y_0(x,t)$ that
grows before it decays.  This is achieved by the following thought
experiment.  Suppose the periodic perturbation $W(x,t)$ is fixed and
even in time $t$.  Now pick any initial condition $y_0(x)$ so that
$\langle x\rangle^{7/2}y_0(x)\in L^2(\R)$.
Proposition~\ref{prop:general} guarantees that after a sufficiently
large number $N$ of periods, the size of the solution of
(\ref{eq:KWstartingpoint}) when measured in the weighted $L^2$ norm is
as small as we please.  Note that throughout this process, the
solution continues to satisfy $\|\langle
\cdot\rangle^{7/2}y(\cdot,t)\|_2<\infty$.  So now, start again at
$t=0$ with the new initial condition $y_0(x)=y(x,t=NL)^*$.  Since the
potential is real and even in time, integration of
(\ref{eq:KWstartingpoint}) with this new initial condition is, up to
complex conjugation, equivalent to integration {\em backwards in time}
from the time $t=NL$ with the initial condition $y(x,t=NL)$.  So we
know that for this very small initial condition, the weighted $L^2$
norm must first grow to an order one size at time $t=NL$ as the decay
process transiently reverses itself, before ultimately giving way to
decay over longer times.  The existence of such solutions does not
violate the statement of Proposition~\ref{prop:general} or
Proposition~\ref{prop:detail} because if one keeps the same initial
condition and then makes the perturbation smaller yet again, the
connection with the time-reversed problem is lost for this initial
condition, and decay occurs sooner. \pfend

%The condition $\|\langle
%\cdot\rangle^{7/2}y_0(\cdot)\|_2<\infty$ is very important.  Assume
%that $\tilde{\op{W}}(-t)=\tilde{\op{W}}(t)$ and consider the thought
%experiment of taking $y_0(x)$ in the Schwartz class, so the condition
%is satisfied, and solving for the corresponding $y(x,t)$.  From
%Proposition~\ref{prop:detail}, the bound state amplitude decays at
%first exponentially, and then algebraically, until it is as small as
%we please.  So pick some $n$ with $T=nL>0$ so that the bound state has
%decayed by some amount independent of $|\|\tilde{\op{W}}(\cdot)\||$,
%and reset $y_0(x)=y(x,T)^*$.  Solving the problem
%(\ref{eq:KWstartingpoint}) with this new initial condition is
%equivalent to simply integrating backwards; in particular at time $T$,
%we must recover the original order one bound state amplitude, which is
%in apparent contradiction to Proposition~\ref{prop:detail}.  The
%resolution of the ``paradox'' is simply that once the radiative decay
%process begins, the presence of packets of dispersive waves at
%$x=\pm\infty$ implies that the function $\langle x\rangle^{7/2}y(x,t)$
%is no longer in $L^2(\R)$ as a function of $x$ for any $t>0$.  So,
%this condition is designed to exclude the possibility of radiation
%entering the domain from infinity and pumping the bound state.  \pfend

By the same arguments applied in the above discussion of the decay constant
$\Gamma$, it follows that there is an alternative formula for $\Lambda$:
\begin{equation}
\Lambda=\sum_{n=-\infty}^{n_0-1}
\int_0^\infty\frac{|N_n(\sqrt{\sigma})|^2d\sigma}
{4\sqrt{\sigma}(\sigma-\sigma_n)} - \sum_{n=n_0}^\infty
{\rm P.\,V.\,}\int_0^\infty\frac{|N_n(\sqrt{\sigma})|^2d\sigma}
{4\sqrt{\sigma}(\sigma-\sigma_n)}\,.
\end{equation}
Again, if $N_n(\lambda)$ has an expansion in a small parameter
$\epsilon$ of the form $N_n(\lambda)=\epsilon N_{1,n}(\lambda) +
O(\epsilon^2)$ then the leading term of $\Lambda$ is of the form
$\epsilon^2\Lambda_2$, where $\Lambda_2$ as given by
(\ref{eq:Lambshift}) was resolved by the multiple scales analysis.  This
frequency shift associated with the decay of the bound state is the
Lamb shift.

From these results, one recovers the true dynamics by setting $f(x,t)=
(\op{P}(t)y(\cdot,t))(x,t)$, where $\op{P}(t)$ is the periodic
operator that appeared in the Floquet factorization of the propagator
$\op{U}(t)$ for the periodic unperturbed Hamiltonian $\op{H}_0(t)$.
Since
\begin{equation}
(\op{P}(t)e^{2i\beta_{\rm b}t}\Psi_{\rm b}^{\rm (o)}(\cdot,t))(x,t)=
(\op{P}(t)e^{-it\op{B}}\Psi_{\rm b}^{\rm (o)}(\cdot,t))(x,t)=
(\op{U}(t)\Psi_{\rm b}^{\rm (o)}(\cdot,t))(x,t)=\Psi_{\rm b}^{\rm (o)}(x,t)\,,
\end{equation}
it follows that the time-dependent projection of $f(x,t)$ onto the
bound state Bloch function $\Psi_{\rm b}^{\rm (o)}(x,t)$ is uniformly
approximated by
\begin{equation}
B_{\rm b}(t)\doteq\langle \Psi_{\rm b}^{\rm (o)}(\cdot,t),f(\cdot,t)\rangle \sim
\langle \Psi_{\rm b}^{\rm (o)}(\cdot,0),f(\cdot,0)\rangle
e^{-\Gamma |t|}e^{i\Lambda t}e^{-i\int_0^t M(s)\,ds}\,.
\end{equation}

For the system restricted to the odd part of $L^2(\R)$, these theorems
provide justification for the formal multiple scales analysis carried
out above, and more.  They globally describe the decay process for all
time, where the multiple scales calculation only attempts to capture
the dynamics over time scales of length $\Gamma^{-1}$.  On the other
hand, since the rate of free dispersive decay is not sufficient in the
even case to apply this detailed theory, we must settle for the
multiple scale expansions.

\section{Applications in Planar Waveguide Optics}
\label{sec:application}
\setcounter{equation}{0} In this section, we present a physical
application of the kinds of perturbed time-dependent Schr\"odinger
equations we have been studying in detail.  This will provide a
concrete family of perturbations $W(x,t)$ that we can use in
subsequent numerical experiments.
\subsection{Time-dependent Schr\"odinger equations in waveguide optics.}
For completeness, we present here a brief derivation of the
time-dependent Schr\"odinger equation as it occurs in the paraxial
theory of monochromatic waveguide optics.  Consider Maxwell's wave
equation for the electric field vector $\vc{E}(\vc{x},t)$ in a planar
($\vc{x}=(y,z)$) dielectric medium with isotropic, inhomogeneous
linear susceptibility $\chi^{(1)}(\vc{x},t)$
\begin{equation}
\Delta\vc{E}-\frac{1}{c^2}\vc{E}_{tt}-\nabla(\nabla\cdot\vc{E})=
\frac{1}{c^2}[\chi^{(1)}(\vc{x},t)\ast\vc{E}]_{tt}\,.
\end{equation}
Here, the asterisk indicates convolution in time.  A Fourier transform
(denoted with the operator $\op{F}$) in $t$ with dual variable
$\omega$ (the optical frequency) leads to
\begin{equation}
\Delta\op{F}\vc{E}-\nabla(\nabla\cdot\op{F}\vc{E})+
\frac{\omega^2n^2(y,z,\omega)}
{c^2}\op{F}\vc{E}=0\,,
\label{eq:Maxwellwave}
\end{equation}
where the refractive index $n$ is defined by $n^2(y,z,\omega)\doteq 1+
(\op{F}\chi^{(1)})(y,z,\omega)$.  We now assume that the inhomogeneity
is weak, so that gradients of $(\op{F}\chi^{(1)})(y,z,\omega)$ are
small.  This implies that in the absence of any free charges, the
approximate relation $\nabla\cdot\op{F}\vc{E}\approx 0$ follows from
the exact relation for the electric displacement
$\nabla\cdot\op{F}\vc{D}=0$.  Neglecting the divergence term in
(\ref{eq:Maxwellwave}), one may then choose any unit vector $\vc{e}$
and set $(\op{F}\vc{E})(y,z,\omega)= \phi(y,z,\omega)\vc{e}$, which
gives the Helmholtz or scalar wave equation for $\phi$:
\begin{equation}
\phi_{zz}+\phi_{yy}+\frac{\omega^2n^2(y,z,\omega)}{c^2}\phi=0\,.
\end{equation}

In the design of integrated optical devices, the inhomogeneity in the
refractive index is a localized modulation of a ``background index''
$n_0(\omega)$.  Choose a fixed length scale $L_0$ and nondimensionalize
by setting $z/L_0=\delta^{-2}Z$ and $y/L_0=\delta^{-1}Y$, where $\delta$
is a dimensionless parameter, and $Y$ and $Z$ are dimensionless coordinates.
Setting 
$\phi(y,z)=f(Y,Z)\exp(i\beta Z/\delta^2)$,
where $\beta=L_0\omega n_0(\omega)/c$ is also dimensionless,
one arrives at
\begin{equation}
2i\beta \delta^{-2}f_{Z}+f_{ZZ}+\delta^{-2}f_{YY}+\beta^2\delta^{-4}
\left[\frac{n^2(YL_0\delta^{-1},ZL_0\delta^{-2},\omega)}{n_0^2(\omega)}-1
\right]f=0\,.
\end{equation}
With the definition
\begin{equation}
Q(Y,Z;\omega)\doteq -\frac{1}{2\delta^2}\left[\frac{n^2(YL_0\delta^{-1},
ZL_0\delta^{-2},\omega)}{n_0^2(\omega)}-1\right]\,,
\label{eq:Qdef}
\end{equation}
we see that the formal limit of $\delta\downarrow 0$ with $\beta$ and
$Q(Y,Z;\omega)$ held fixed yields the paraxial wave equation
\begin{equation}
i\beta f_Z + \frac{1}{2}f_{YY} - \beta^2 Q(Y,Z;\omega)f=0\,.
\label{eq:paraxial}
\end{equation}
The potential function $Q(Y,Z;\omega)$ vanishes as the refractive
index approaches its background value $n_0(\omega)$, say as $Y$ and
$Z$ go off to infinity (at least in most directions).  Given a fixed
function $Q(Y,Z;\omega)$, we see that the paraxial approximation made
here ($\delta\downarrow 0$) is valid if the modulation in the
refractive index is weak, slowly varying, and more slowly varying in
the $z$ direction than in the $y$ direction.  That is, a fixed
function $Q(Y,Z;\omega)$ provides an {\em asymptotic description} of a
{\em family} of physical refractive index profiles parametrized by
$\delta\ll 1$:
\begin{equation}
n^2(y,z,\omega;\delta)= n_0^2(\omega)-2\delta^2n_0^2(\omega)Q(\delta y/L_0,
\delta^2 z/L_0;\omega)\,.
\label{eq:physicalindex}
\end{equation}
Note that these assumptions about the refractive index justify {\em a
posteriori} our neglect of the term $\nabla(\nabla\cdot\op{F}\vc{E})$
in the original wave equation, because in the limit $\delta\rightarrow 0$,
gradients of $n^2(y,z,\omega;\delta)$ necessarily vanish.

\subsection{Spectral properties of paraxial waveguides.}
In optical waveguide theory, integration (numerical or otherwise) of
the linear Schr\"odinger equation (\ref{eq:paraxial}), also known as
the beam propagation method, is one of the main tools for studying the
optical properties of ``long'' planar structures like gradual fiber
tapers or channel waveguide junctions, in which backward reflecting
waves can be neglected.  In this connection, a common problem that
arises is the description of the change in behavior of a waveguiding
structure as the optical frequency is varied in the neighborhood of
some frequency $\omega_0$.  If the structure $n^2(y,z,\omega)$ is one
that admits the paraxial approximation, we can use the theory
described above as a model.  In this case, it is convenient to choose
the length scale $L_0$ so that at the frequency $\omega_0$ we have
$\beta=1$.  With this choice, we think of $\beta=\beta(\omega)$ as a
function of frequency satisfying $\beta(\omega_0)=1$.  With the
function $Q(Y,Z;\omega)$ chosen consistently, the problem becomes one
of studying the dependence of solutions of (\ref{eq:paraxial}) on the
frequency parameter $\omega$ near $\omega_0$.  With the change of
variables $x=Y\sqrt{\beta(\omega)}$ and $t=Z$, the equation
(\ref{eq:paraxial}) takes the form
\begin{equation}
if_t = \left(\ -\frac{1}{2}\partial_x^2 + V_0(x,t)\ \right) f\ +\ W(x,t)f\,,
\label{eq:schrod}
\end{equation}
where
\begin{equation}
V_0(x,t)=Q(x,t;\omega_0)\,,
\end{equation}
and the correction to the potential is given by
\begin{equation}
W(x,t)=\beta(\omega)Q(x/\sqrt{\beta(\omega)},t;\omega)-Q(x,t;\omega_0)\,.
\end{equation}
Setting $\epsilon=\omega/\omega_0-1$, we see that
$W(x,t)=W(x,t;\epsilon)$ is uniformly small in $\epsilon$ if
$Q(Y,Z;\omega)$ is in the Schwartz space with respect to $Y$.  We have
the expansion
\begin{equation}
W(x,t;\epsilon)= \epsilon
\omega_0\left[\beta'(\omega_0)\left(1-\frac{x}{2}\partial_x\right)
Q(x,t;\omega_0)
+ \partial_\omega Q(x,t;\omega_0)\right] + O(\epsilon^2)\,,
\label{eq:expansionofW}
\end{equation}
uniformly in $x$ and $t$.  If the frequency range of interest is
sufficiently small, then it is often a good approximation to consider
the problem to be {\em dispersionless}, so that the refractive index
$n(y,z,\omega)$ is independent of $\omega$.  In this paper, we will
accordingly consider the function $Q$ to be independent of $\omega$ in
which case $Q(x,t;\omega)=V_0(x,t)$ for all $\omega$ in the range of
interest, and we can drop the corresponding term in
(\ref{eq:expansionofW}).

Suppose now that we choose to study a refractive index profile
$n^2(y,z)$ that is even in $y$ and periodic in $z$, such that after
choosing a frequency $\omega_0$ and nondimensionalizing, the function
$V_0(x,t)$ is one of the separable potentials described in detail at
the end of \S \ref{sec:separable}.  Over length scales where the
paraxial approximation is valid, this periodically modulated channel
waveguide will actually have two ``breather modes'', approximately
described by the bound states $\Psi_{\rm b}^{\rm (e)}(x,t)$ and
$\Psi_{\rm b}^{\rm (o)}(x,t)$.  The effect of not being fully in the
paraxial limit (that is, $\delta$ is small but finite) is that the
modes will very slowly attenuate as they propagate forward due to a
small coupling to backward propagating fields.  This small attenuation
occurs at all frequencies near $\omega=\omega_0$ in a way that can be
quantified \cite{nonparaxial}.  However, the profile $n^2(y,z)$ is
very special in that at the frequency $\omega=\omega_0$ there is no
coupling between the bound modes and any forward propagating radiation
modes.  This additional coupling would indeed be present for
``typical'' $z$-periodic waveguide profiles $n^2(y,z)$.

In fact, the theory developed in \S \ref{sec:analysis} can be applied
to the perturbed problem (\ref{eq:schrod}) because the unperturbed
potential $V_0(x,t)$ and the perturbation $W(x,t;\epsilon)$ are both
even functions of $x$ that are periodic in $t$ with the same period
$L$.  This theory shows that the additional attenuation due to
coupling to forward propagating radiation, while completely suppressed
at the frequency $\omega_0$, reemerges upon detuning the frequency
slightly from $\omega_0$.  Suppose the waveguide is cleaved at $z=0$ and
is illuminated at this face with a broadband source consisting of many
frequencies $\omega$.  After some distance all of the frequencies will
have attenuated somewhat due to backscattering (weak non-paraxiality).
However, all frequencies except $\omega_0$ will {\em additionally}
decay by forward propagating radiation damping.  The waveguide
will therefore preferentially ``pass'' light of the
frequency\footnote{Actually, the ``background'' attenuation due to
nonparaxiality ($\delta\neq 0$) decreases slightly with increasing
frequency.  When this effect is combined with the frequency-dependent
decay calculated from the paraxial approximation, the preferred
frequency for which the loss is minimal is detuned slightly upward
\cite{nonparaxial} by an amount that is $O(\delta^2)$.} $\omega_0$.
These effects were observed numerically in \cite{BesleyOL,SAM}.

Note that from the point of view of optical waveguide theory, the
periodicity of the index $n(x,z)$ in $z$ is an important feature,
since it gives rise to an attenuated frequency response that is a
symmetric function of frequency $\omega$ in the neighborhood of
$\omega_0$.  Thus, attenuation occurs whether $\omega$ is less than or
greater than the frequency $\omega_0$ of structural instability.  By
constrast, channel waveguides, where $n(y,z)$ is independent of $z$,
also exhibit frequency-dependent structural instability at {\em cutoff
frequencies} where the number of bound states changes.  However, in
such waveguides the number of bound states (guided modes) is always an
increasing function of frequency \cite{M74}, which implies that an
input beam that matches onto a mode at its cutoff frequency $\omega_0$
will attenuate for $\omega$ slightly less than $\omega_0$ but will
remain bound and thus give rise to a significant transmission for
$\omega$ slightly greater than $\omega_0$.  Thus, whereas channel
waveguides with $z$-independent refractive index profiles can behave
as ``high-pass'' components, $z$-periodic waveguides that at
frequency $\omega_0$ are modeled by separable potentials can behave as
``band-pass'' components.

\section{Numerical Simulations}
\label{sec:numerics}
\setcounter{equation}{0} Here, we describe some numerical simulations
we performed to verify the analytical predictions where we expect them
to apply.  We also would like to explore the behavior of the perturbed
system in parameter regimes where we expect zero-energy resonances
(see the remark at the end of \S \ref{sec:sub-ms}) to prevent the theory
from applying in its current form.  For concreteness, we considered
periodic perturbations of two problems, each associated with a
particular two-soliton separable periodic potential.  The particular
perturbation we selected was exactly the type considered in \S
\ref{sec:application}, namely, given a separable two-soliton periodic
potential $V_0(x,t)$, we numerically integrated the equation
\begin{equation}
i\partial_t f +\frac{1}{2(1+\epsilon)}\partial_x^2 f - (1+\epsilon)V_0(x,t)f=0\,,
\label{eq:actuallyintegrated}
\end{equation}
for several small values of $\epsilon$.  This problem differs from the
type to which the theory developed above applies only by a rescaling
of $x$; in particular, the time scale is unaffected.

Let us give some details about our numerical scheme.  We used a
Fourier split-step method with a local truncation error of $O(\Delta
t^3)$ \cite{NewellMoloney}.  The spatial domain $[x_L,x_R]$ of
$[-80,80]$ in the ``non zero-energy resonance'' case and $[-40,40]$ in
the ``zero-energy resonance'' case (see below for more details about
these two cases) was discretized into $1024$ points.  The scheme
splits the Hamiltonian into two parts: $\op{H}(t)=\op{H}_1 +
\op{H}_2(t)$, where
\begin{equation}
\op{H}_1 \doteq -\frac{1}{2(1+\epsilon)}\partial_x^2\,,\hspace{0.2 in}
\op{H}_2(t) \doteq (1+\epsilon)V_0(x,t)\,.
\end{equation}
Let $\op{U}^\epsilon(t,s)$ denote the propagator associated with 
(\ref{eq:actuallyintegrated}).  Let $\op{U}^\epsilon_1(t-s)$ and 
$\op{U}^\epsilon_2(t,s)$ be those associated with $\op{H}_1$ and 
$\op{H}_2(t)$.  Then, the numerical scheme approximates the true integration
over a time step of size $\Delta t$ as follows:
\begin{equation}
\op{U}^\epsilon(t+\Delta t,t)\sim\op{U}_1^\epsilon(\Delta t/4)
\op{U}_2^\epsilon (t+3\Delta t/4,t+\Delta t/4)
\op{U}_1^\epsilon(\Delta t/4)\,,
\end{equation}
which has an error of order $\Delta t^3$.  It is easy to see that,
after getting started with a quarter-step, and until finishing with a
quarter-step, iterating this approximation to the propagator
$\op{U}^\epsilon(t,s)$ over many steps amounts to simply alternating
between $\op{U}^\epsilon_1$ and $\op{U}^\epsilon_2$ each acting over
a half-step of length $\Delta t/2$.

So, in each half-step, only one of the two parts is integrated.  The
half-step involving $\op{H}_1$ is carried out in the Fourier transform
domain where one multiplies by the explicit exponential of the
operator.  This step is thus exact in time, so that the only error
appears in discretizing the Fourier transform and is smaller than any
power of $\Delta x$ if the functions to be differentiated are taken to
be arbitrarily smooth.  The half-step involving $\op{H}_2(t)$ is done
{\em exactly} because we have explicit formulas for $V_0(x,t)$ and it
is possible to find an {\em explicit} exponential of $\op{H}_2(t)$.
That is, we can write down a formula for the multiplication operator:
\begin{equation}
\op{U}_2^\epsilon (t,t_0)=\exp\left(i(1+\epsilon)\int_{t_0}^t
V_0(x,s)\,ds\right)\,,
\end{equation}
and use it in the code.  Since the temporal gradients of $V_0(x,t)$
can be large in some parts of each period and small in others, we
adjusted the time step throughout the period.  

We expect the perturbation to generate radiation from the central bound
region of the potential, and we need to remove this radiation from the problem
as it moves to large $|x|$.  To take care of this we used a ``sponge layer''
in which we effectively add a term of the form
\begin{equation}
-id\cdot\left[\exp\left(-\left(\frac{x-x_R}{w}\right)^2\right) + 
\exp\left(-\left(\frac{x-x_L}{w}\right)^2\right)\right]f\,,
\end{equation}
to the right-hand side of (\ref{eq:actuallyintegrated}) for a positive
damping factor $d$ and width $w$.  These parameters were adjusted
heuristically until it was observed, roughly speaking, that no energy
was being artificially drawn out of the center and that no energy that
was radiated outward was either reflected or transmitted through to
the other side of the periodic domain.

We integrated for $50$ periods.  In all the experiments it was
arranged that the fundamental period was $L=2\pi$.  We initialized the
field $f$ at $t=0$ to be a snapshot of either the odd or the even mode
of the unperturbed problem.  Then, after integrating, we calculated
the projection of the numerical solution onto the exact solution of
the unperturbed problem, defining:
\begin{equation}
B_{\rm b}(t)\doteq \langle \Psi_{\rm b}^{\rm
(e,o)}(\cdot,t),f(\cdot,t)\rangle\,.
\end{equation}
We verified the accuracy of the code by checking that for $\epsilon=0$
we had $B_{\rm b}(t)\equiv 1$ to several digits, even in the presence of the
damping in the sponge layer.  Note that the function $B_{\rm b}(t)$ is related
to $A_{\rm b}(t)$ by the simple relation:
\begin{equation}
B_{\rm b}(t)=A_{\rm b}(t)e^{-2i\beta_{\rm b}t}\,.
\end{equation}

\subsection{Away from parametric zero-energy resonance.}
For the first experiments, we selected $\rho_1=1/4$ and $\rho_2=3/4$ as
the parameters of the function $V_0(x,t)$.  It is easy to check that the
period is $L=2\pi$, and that the Floquet exponent of both odd and even
bound states may be taken to be $\beta_{\rm b}=\rho_1^2=1/16$.  Therefore
the resonances are explicitly given by
\begin{equation}
\sigma_n=\frac{n}{2}-\frac{1}{16}\,,
\end{equation}
none of which are equal to zero.  This means that there is no
parametric zero-energy resonance, although in the even case there
still is a zero-energy resonance corresponding to insufficient
dispersive decay.  In this case, the formula for the decay constant
$\Gamma$ makes sense for both odd and even parity.  Furthermore, for
odd parity, we have a proof that the asymptotic expansion obtained
previously is indeed valid.

\begin{figure}[h]
\begin{center}
\mbox{\psfig{file=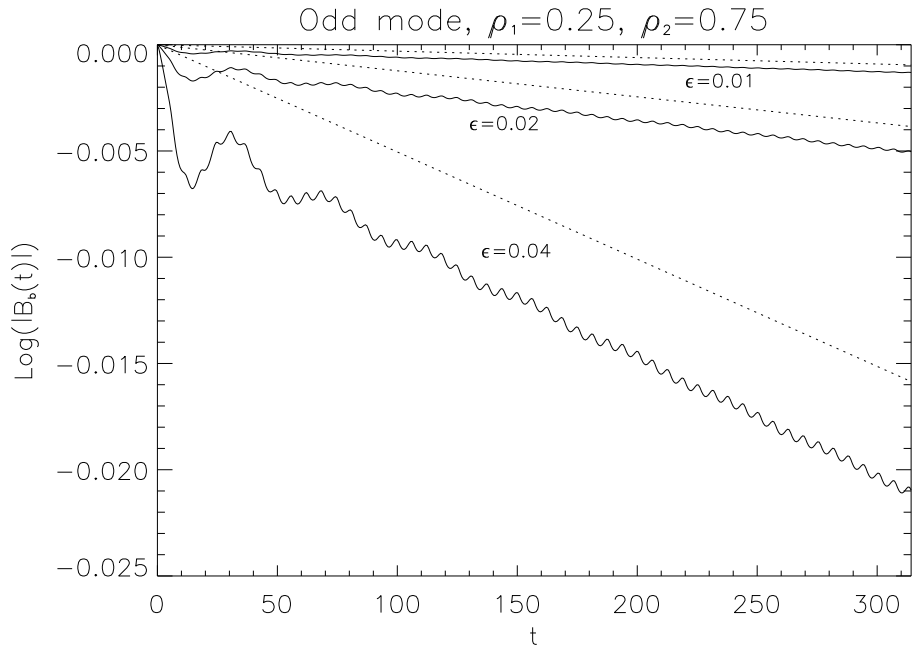,width=5 in}}
\end{center}
\caption{\em The magnitude of the projection of the solution onto the
bound state.  Odd mode.  $\rho_1=1/4$, $\rho_2=3/4$.  Solid lines are
the output from numerical simulations.  Dotted lines are the
analytical predictions.}
\label{fig:R9oddpertabs}
\end{figure}

\begin{figure}[h]
\begin{center}
\mbox{\psfig{file=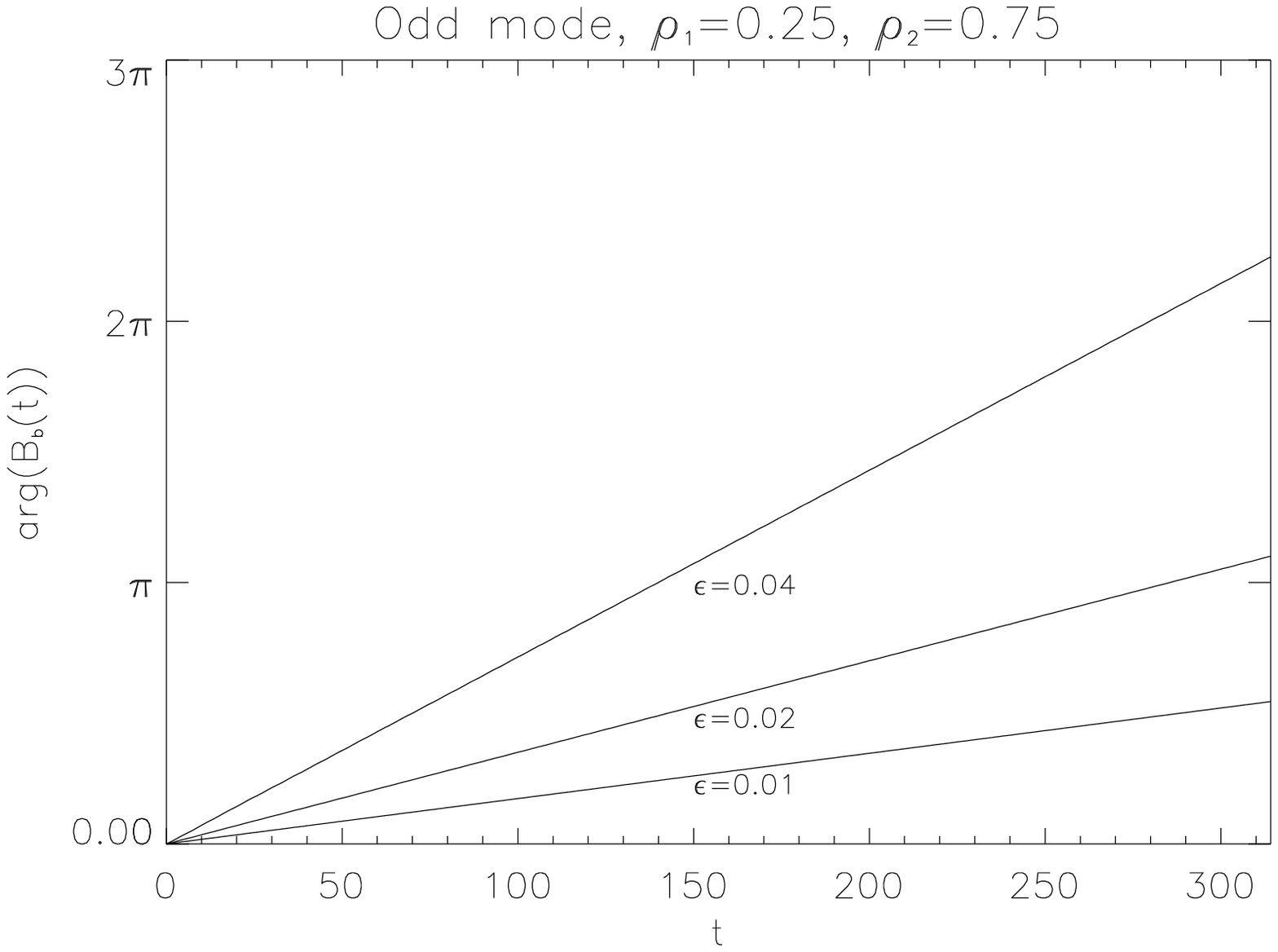,width=5 in}}
\end{center}
\caption{\em The phase of the projection of the solution onto the
bound state.  Odd mode.  $\rho_1=1/4$, $\rho_2=3/4$.}
\label{fig:R9oddpertphase}
\end{figure}

In Figure~\ref{fig:R9oddpertabs}, we show plots of $\log(|B_{\rm b}(t)|)$ for
$\epsilon=0.04$, $\epsilon=0.02$, and $\epsilon=0.01$ for an initial
condition of odd parity.  The numerical results are plotted with solid
curves, and superimposed are corresponding graphs of $-\Gamma |t|$
calculated from the analytical formula, the analogue of Fermi's golden
rule, and shown with dotted lines.
The main observation here is that the graphs follow the corresponding
straight lines, which have slopes that scale like $\epsilon^2$, as
expected.  The deviation from the straight lines appears to scale like
$\epsilon^2$ as well, and to decay in time.  In 
Figure~\ref{fig:R9oddpertphase}, we give corresponding plots of the argument
of $B_{\rm b}(t)$ for an initial condition of odd parity.  
In these plots, it is easy to see that the phase grows roughly
linearly in time, with slope that is $O(\epsilon)$.  This is the
contribution to the frequency shift of the term $\overline{M}$, which
is indeed $O(\epsilon)$.

\begin{figure}[h]
\begin{center}
\mbox{\psfig{file=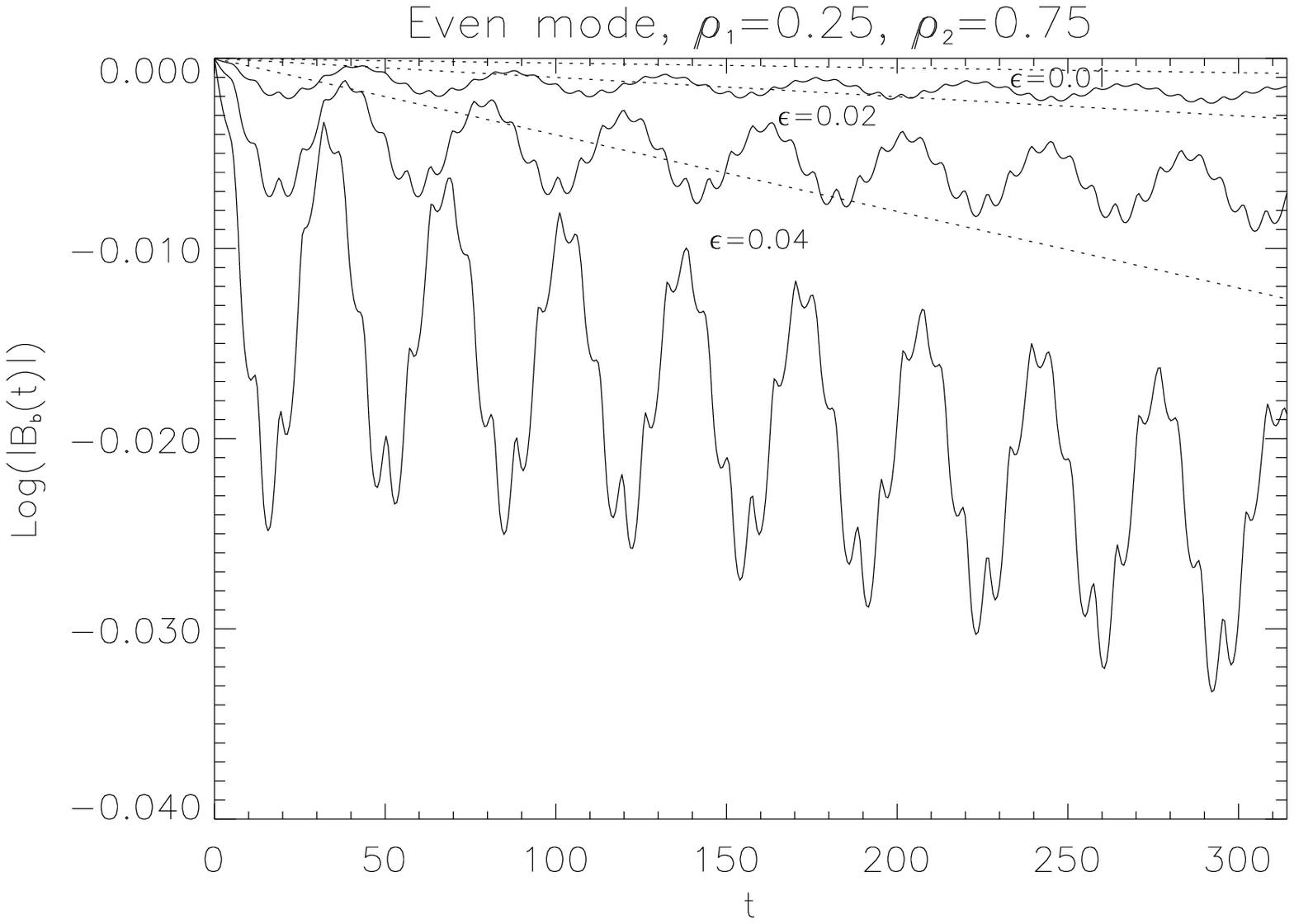,width=5 in}}
\end{center}
\caption{\em The magnitude of the projection of the solution onto the
bound state.  Even mode.  $\rho_1=1/4$, $\rho_2=3/4$.  Solid lines are
the output from numerical simulations.  Dotted lines are the
analytical predictions.}
\label{fig:R9evenpertabs}
\end{figure}

\begin{figure}[h]
\begin{center}
\mbox{\psfig{file=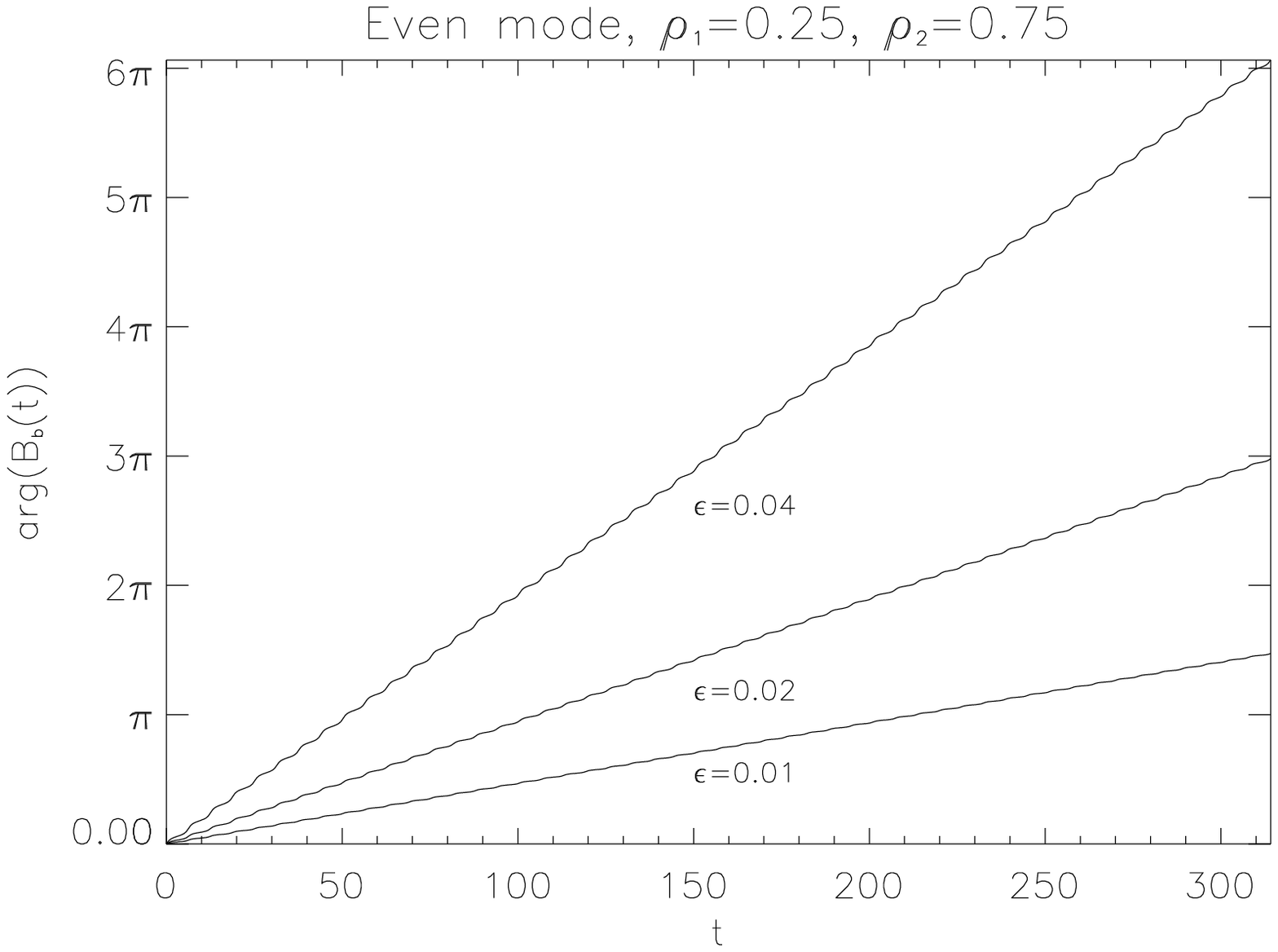,width=5 in}}
\end{center}
\caption{\em The phase of the projection of the solution onto the
bound state.  Even mode.  $\rho_1=1/4$, $\rho_2=3/4$.}
\label{fig:R9evenpertphase}
\end{figure}

Now, we consider an even initial condition, with corresponding
projection $B_{\rm b}(t)$ onto the even mode of the exact solution for
$\epsilon=0$.  Figure~\ref{fig:R9evenpertabs} contains plots of
$\log(|B_{\rm b}(t)|)$ as calculated from the numerical data for
$\epsilon=0.04$, $\epsilon=0.02$, and $\epsilon=0.01$ shown in solid
curves.  
Also plotted are the corresponding decay curves $-\Gamma |t|$
shown with dotted lines.  Although for even parity there is
insufficient dispersive decay for the results of \cite{KW} to 
apply, the decay constant $\Gamma$ (or more precisely as it is
obtained in the multiple scale analysis, $\epsilon^2\Gamma_2$) is
finite because none of the resonances $\sigma_n$ are zero, and we see
that the multiple scale theory accurately predicts the rate of decay
of the bound state even in this case.  
The plots of the phase of $B_{\rm b}(t)$
are shown in Figure~\ref{fig:R9evenpertphase}.
Again, one sees that the rate of drift of the phase is $O(\epsilon)$,
as predicted by the multiple scale theory.

The significant new feature apparently contributed by the lack of
sufficient dispersive decay for initial conditions of even parity
appears to be the quality of the deviations in $|B_{\rm b}(t)|$ from the
``backbone'' decay $e^{-\Gamma |t|}$.  Not only are they larger for
fixed $\epsilon$ than for initial conditions of odd parity, but they
have an undulatory character that suggests the possible contribution
of {\em subharmonic} frequencies to the dynamics.  The period of the
undulations superimposed on the decay appears to be long compared with
$L$, the fundamental period of the problem, but also appears to be
more or less independent of $\epsilon$.

\subsection{At parametric zero-energy resonance.}
As a second set of experiments, we considered a potential energy
function $V_0(x,t)$ obtained from the parameters $\rho_1=1/\sqrt{2}$
and $\rho_2=1$.  In this case, the period is again $L=2\pi$, and the
Floquet exponents of both modes are $\beta_{\rm b}=\rho_1^2=1/2$.  
The corresponding family of resonances is
\begin{equation}
\sigma_n=\frac{n}{2}-\frac{1}{2}\,.
\end{equation}
One of these values is equal to zero.  This condition for parametric
zero-energy resonance always goes hand-in-hand with another property
of this potential, namely that the Floquet multipliers of both modes
are equal to $1$.  Thus, both odd and even modes are actually periodic
functions of $t$ with period $L$.

\begin{figure}[h]
\begin{center}
\mbox{\psfig{file=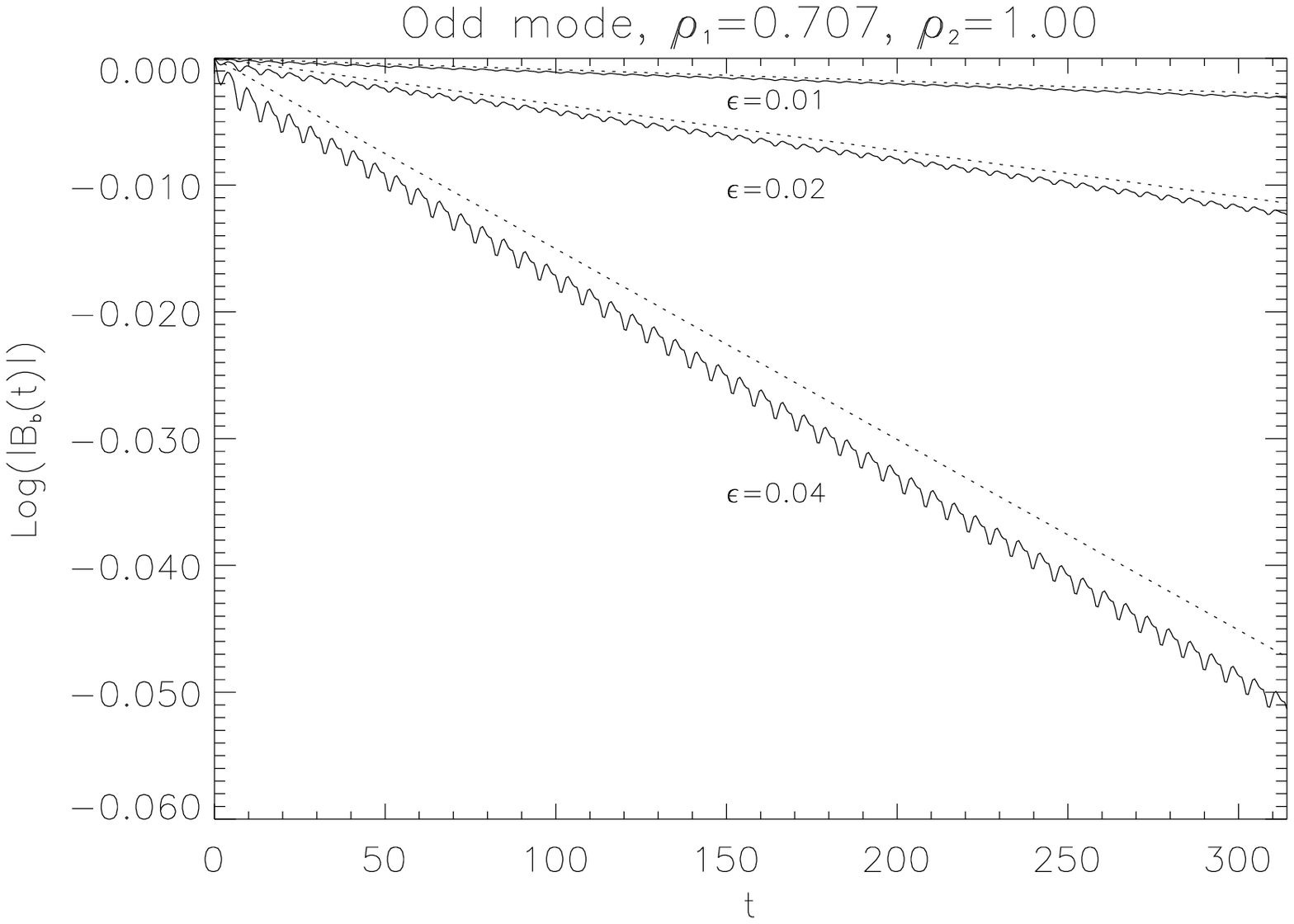,width=5 in}}
\end{center}
\caption{\em The magnitude of the projection of the solution onto the
bound state.  Odd mode.  $\rho_1=1/\sqrt{2}$, $\rho_2=1$.  Solid lines are
the output from numerical simulations.  Dotted lines are the
analytical predictions.}
\label{fig:R2oddpertabs}
\end{figure}

\begin{figure}[h]
\begin{center}
\mbox{\psfig{file=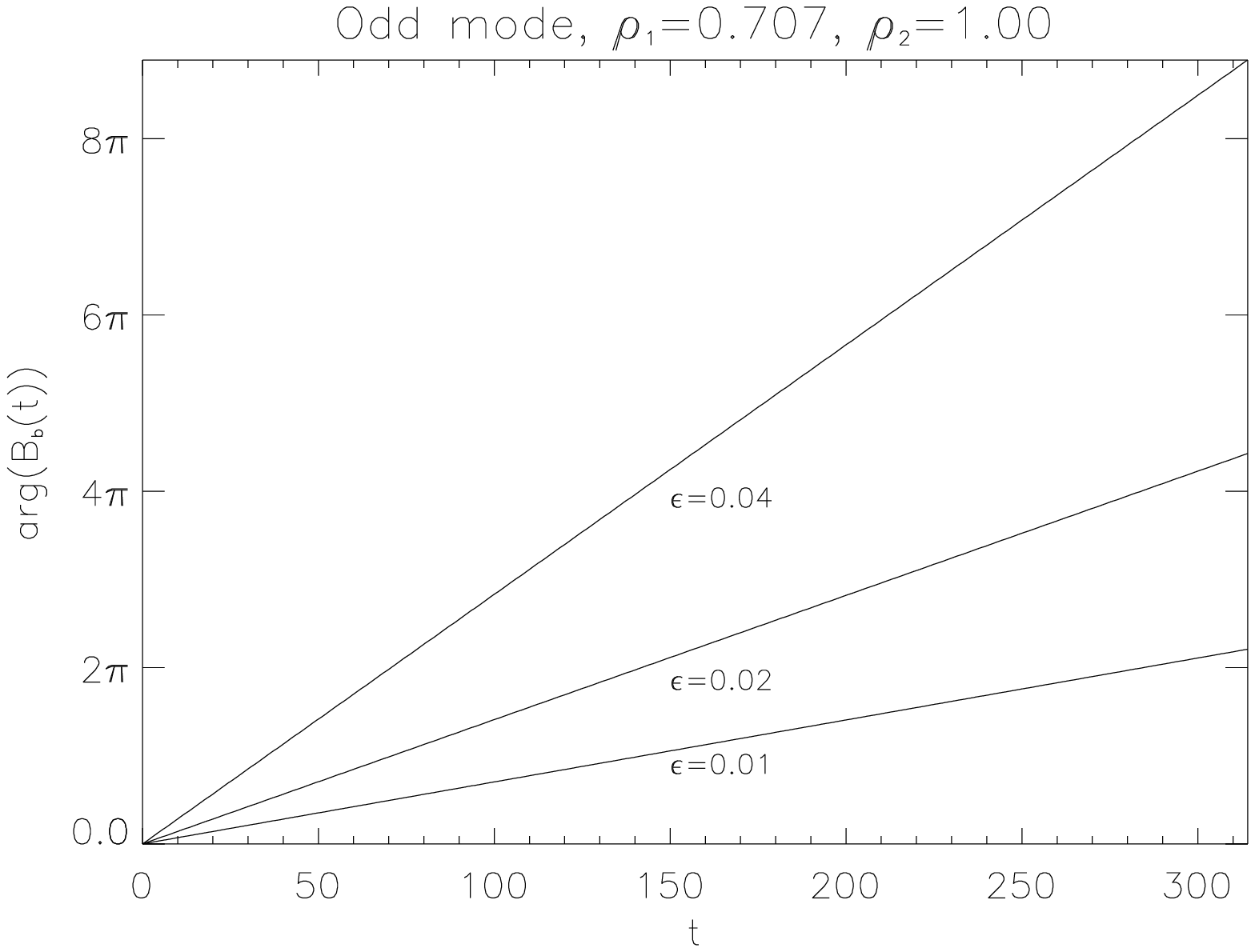,width=5 in}}
\end{center}
\caption{\em The phase of the projection of the solution onto the
bound state.  Odd mode.  $\rho_1=1/\sqrt{2}$, $\rho_2=1$.}
\label{fig:R2oddpertphase}
\end{figure}

At a parametric zero-energy resonance, the dispersive local decay
estimates fail to be sufficient to guarantee the applicability of the
theory in \cite{KW}, {\em also for initial conditions of odd parity}.
However, in the odd case, the formulas for the decay constant $\Gamma$
and the Lamb shift $\Lambda$ are finite because there is sufficient
vanishing in the numerator coming from the missing generalized
eigenfunction at $\sigma=0$ to cancel and overcome the weaker
vanishing of the denominator.  Plots of $\log(|B_{\rm b}(t)|)$ for odd
parity corresponding to $\epsilon=0.04$, $\epsilon=0.02$, and
$\epsilon=0.01$ are shown in Figure~\ref{fig:R2oddpertabs} along with
dotted lines indicating the analytical prediction of decay.
The prediction of the theory appears to be very accurate indeed.  The
plots of the phase of $B_{\rm b}(t)$ are shown in
Figure~\ref{fig:R2oddpertphase}.  
They show the frequency shift scaling like $\epsilon$, as we expect
from the contribution of the term $\overline{M}$.  So it appears that
for initial conditions of odd parity, there is little if any effect of
the parametric zero-energy resonance, although the rate of dispersive
decay is smaller here than at more generic parameter values.

\begin{figure}[h]
\begin{center}
\mbox{\psfig{file=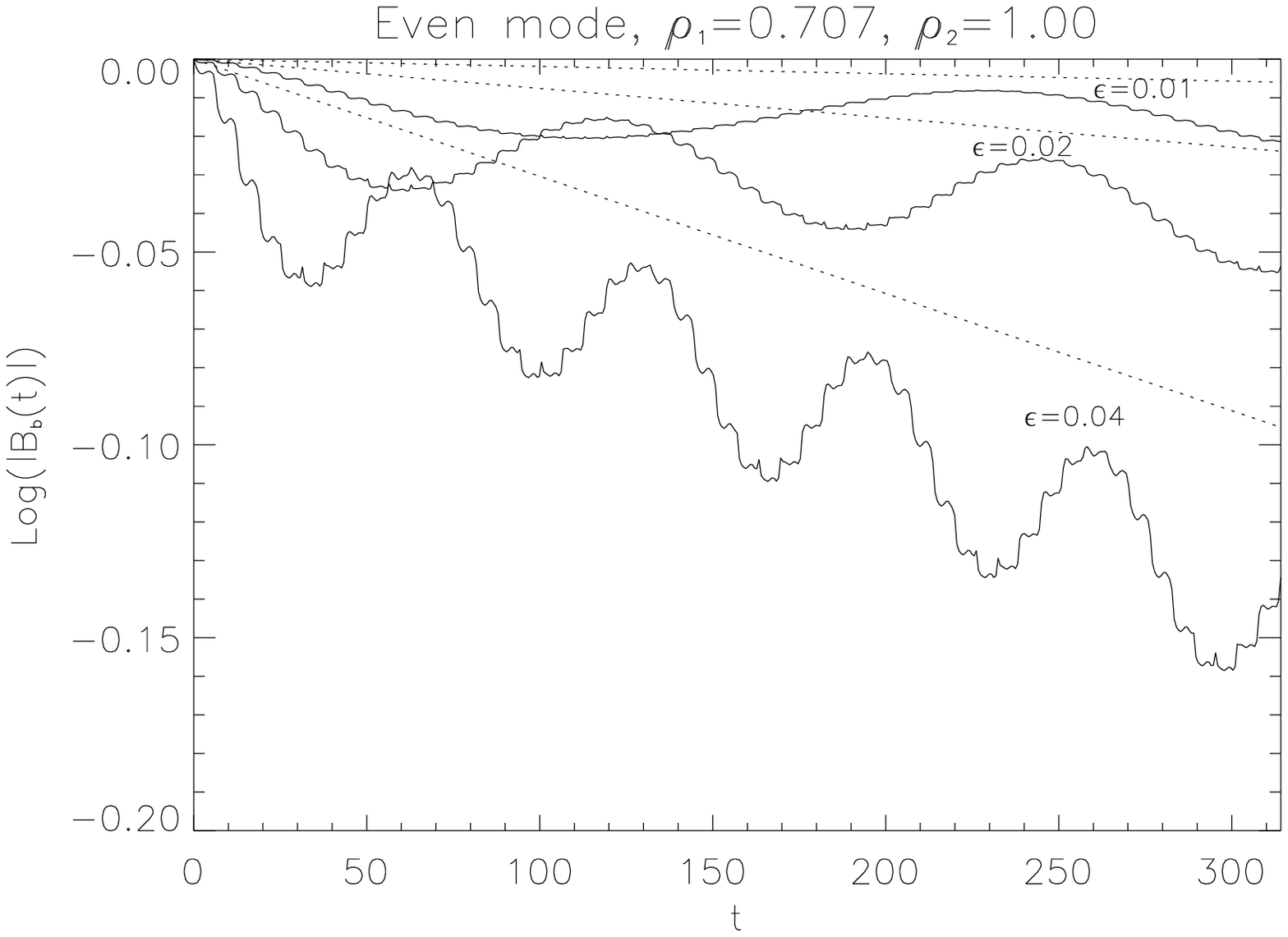,width=5 in}}
\end{center}
\caption{\em The magnitude of the projection of the solution onto the
bound state.  Even mode.  $\rho_1=1/\sqrt{2}$, $\rho_2=1$.  Solid
lines are the output from numerical simulations.  Dotted lines are the
``renormalized'' analytical predictions.}
\label{fig:R2evenpertabs}
\end{figure}

\begin{figure}[h]
\begin{center}
\mbox{\psfig{file=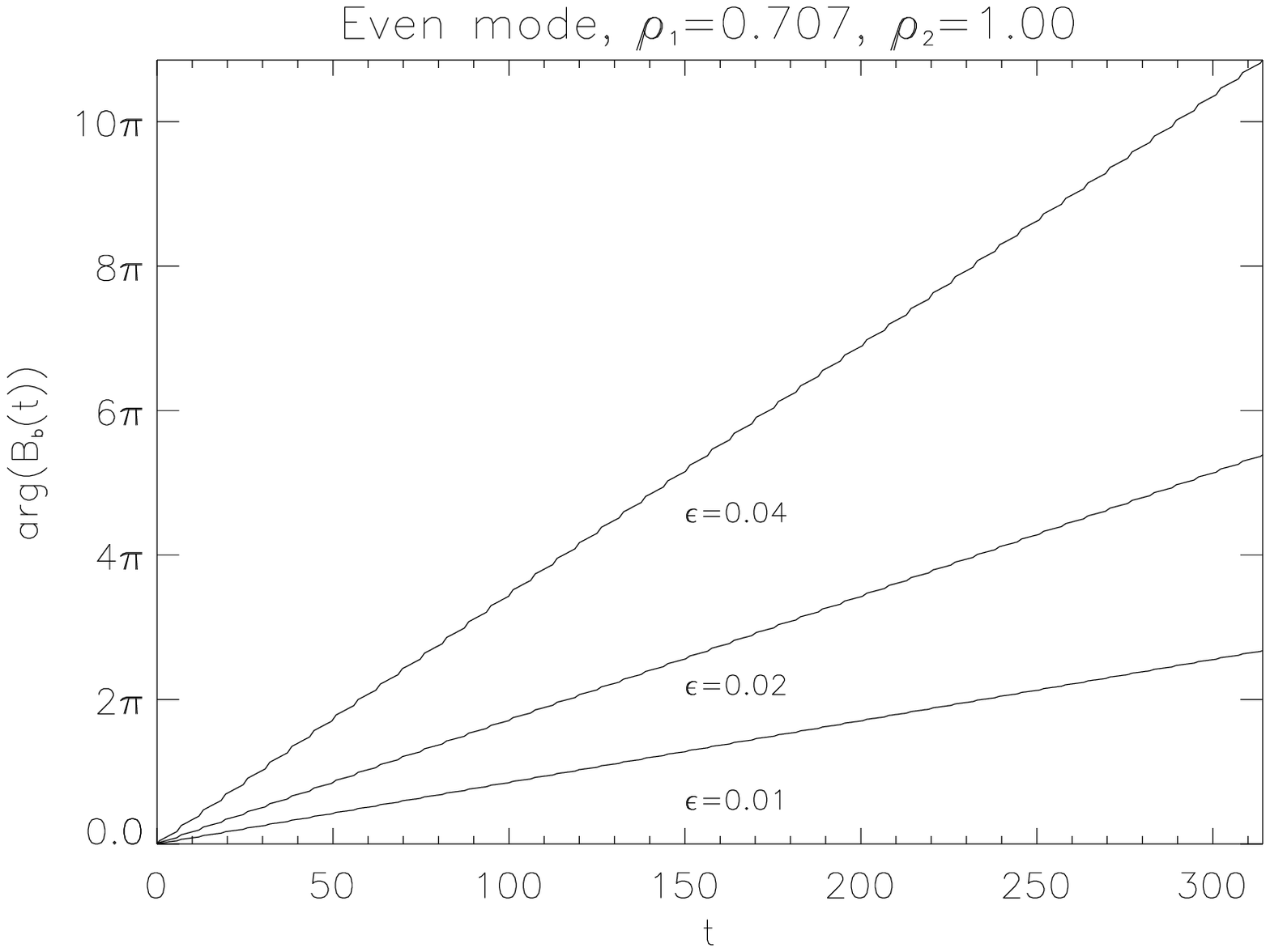,width=5 in}}
\end{center}
\caption{\em The phase of the projection of the solution onto the
bound state.  Even mode.  $\rho_1=1/\sqrt{2}$, $\rho_2=1$.}
\label{fig:R2evenpertphase}
\end{figure}

Finally, let us examine the behavior of initial conditions of even
parity.  For such initial conditions and for these parameter values,
we have {\em both} a simple zero-energy resonance (as one has in the
even case for all parameter values) and a parametric zero-energy
resonance (as occurs only for very special parameter values).  It is
easy to see that both $\Gamma$ and $\Lambda$ are infinite in this
case, and clearly one cannot expect the multiple scale analysis to be
valid.  So what can one expect?  In Figure~\ref{fig:R2evenpertabs} we
plot $\log(|B_{\rm b}(t)|)$ for $\epsilon=0.04$, $\epsilon=0.02$, and
$\epsilon=0.01$, as before.  
This time, rather than superimposing the straight lines $-\Gamma|t|$,
we might try to compare with a ``renormalized'' rate of decay given by
the formula for $\Gamma$ with the term coming from $\sigma=0$ simply
dropped.  The straight lines calculated from the renormalized version
of $\Gamma$ appear dotted on the plots.  We {\em still} see quite good
agreement at the level of a mean drift of $|B_{\rm b}(t)|$.  As in the
previous experiment with even parity, we see subharmonic undulations
about this mean drift.  However, a key point is that whereas
previously the period of these undulations appeared to be more or less
indpendent of $\epsilon$, in this case we note that the period appears
to scale like $\epsilon^{-1}$.  Thus, there is a ``slow'' dynamical
process involving variations of the amplitude that is completely
missed by the multiple scales analysis in its current form.  We must
expect that whatever rescalings are required to balance the blowing up
of $\Gamma$ in the vicinity of a parametric zero-energy resonance will
also introduce interesting subharmonic dynamics on the scale of
$T_1=\epsilon t$ that will reproduce the effects we are seeing
numerically.  As a final remark, the phase of $B_{\rm b}(t)$, as shown in
Figure~\ref{fig:R2evenpertphase}
exhibits no particularly different behavior than was seen in any of
the other experiments.  The frequency adjustment continues to be
dominated by the relatively large term $\overline{M}$ and is therefore
order $\epsilon$.

\section{Conclusions}
In studying the propagation of waves in time-periodic potentials,
considering the problem at hand to be a perturbation of a separable
periodic problem is evidently as easy as, and in many cases more
convenient than working with periodic perturbations of stationary
potential problems.  A particular application to the theory of
periodically modulated optical waveguides in planar dielectric media
allows one to study frequency dependent attenuation properties of
certain optical waveguides.

Many of the difficulties described in our paper concern the influence
of zero-energy resonances.  These are generically not present (that
is, for most separable periodic potentials, as for most stationary
potentials), but are {\em always} present when the potential has
sufficient symmetry, as in the evenness considered above.  Many
problems would therefore vanish upon dropping the symmetry.  From one
point of view, this introduces the additional complication of having
multiple bound states that are essentially coupled to one another by
the perturbation.  The study of perturbed multimode problems arises
naturally in the theory of light propagation in optical fibers having
large effective cross-sections.  Some of the necessary modifications
in the theory described in \cite{KW} are described by the same authors
in \cite{KWmultimode}.

Of course another point of view is to keep the symmetry, and hence the
possibility of zero-energy resonance, and study the effect of the
resonance in more detail.  Our numerical experiments suggest that the
effects of such a resonance are most dramatic when the expressions for
$\Gamma$ and $\Lambda$ blow up, but we also see significant effects,
presumably coming simply from the lack of sufficient long time decay
of freely dispersing waves, when these quantities are finite.  An
asymptotic perturbation theory for small $\epsilon$ should be
uniformly valid with respect to parameters like $\rho_1$ and $\rho_2$,
and we plan to investigate zero-energy resonances with such a goal in
mind in future work.

\section{Acknowledgements}
P.~D.~Miller is grateful for the support of the NSF under grant number
DMS 9304580 while at the Institute for Advanced Study.  A.~Soffer is
supported in part by a FAS-Rutgers grant and by NSF grant number DMS
9706780.  M.~I.~Weinstein is supported in part by NSF grant number DMS
9500997.  Collaboration began while M.~I.~Weinstein visited the IAS in
March 1998 as part of the Program in Geometric Partial Differential
Equations organized by Karen Uhlenbeck, whom the authors thank for her
support.

\appendix
\renewcommand{\theequation}{\Alph{section}.\arabic{equation}}
\renewcommand{\thesection}{Appendix \Alph{section}}
\section{Some theory of separable potentials.}
\label{app:separable}
\renewcommand{\thesection}{\Alph{section}} \setcounter{equation}{0}
For completeness, we here give a self-contained description of the
separable potentials for the linear Schr\"odinger equation that are
connected with the soliton theory of vector nonlinear Schr\"odinger
equations.  However, the material is auxiliary and all needed facts
are reproduced in the main text.  The results here are not new
\cite{physd} but some arguments are carried out here in more detail.

Let $q_1(x,t),\dots,q_N(x,t)$ be given smooth bounded complex
functions of real $x$ and $t$, and let $A$ be the vector space of
differentiable $\C^{N+1}$-valued functions of $x$ and $t$.  Let $\lambda$ be a
complex parameter, and consider the two linear operators acting in
$A$:
\begin{equation}
\begin{array}{rcl}
\op{X}(\lambda,\vc{q})&\doteq &\displaystyle
\partial_x-\left[\begin{array}{cc}
-2i\lambda & \vc{q}^T\\-\vc{q}^*&\mat{0}\end{array}\right]\\\\
&=&\partial_x - (-2i\lambda \mat{E}+\mat{U}(\vc{q}))\,,
\end{array}
\label{eq:X}
\end{equation}
where $\mat{E}$ is a matrix whose elements are given by
$E_{ij}=\delta_{i1}\delta_{j1}$ and
\begin{equation}
\mat{U}(\vc{q})=\left[\begin{array}{cc}0 & \vc{q}^T\\-\vc{q}^* & \mat{0}
\end{array}\right]\,,
\end{equation}
and
\begin{equation}
\begin{array}{rcl}
\op{T}(\lambda,\vc{q})&\doteq &\displaystyle
\partial_t-\left[\begin{array}{cc}
-2i\lambda^2+i\vc{q}^T\vc{q}^*/2 & \lambda \vc{q}^T + i\partial_x\vc{q}^T/2
\\-\lambda\vc{q}^*+i\partial_x\vc{q}^*/2 & -i\vc{q}^*\vc{q}^T/2
\end{array}\right]\\\\
&=&\displaystyle \partial_t-\left(-2i\lambda^2\mat{E}
+\lambda\mat{U}(\vc{q}) + \frac{i}{2}\mat{V}(\vc{q})\right)\,,
\end{array}
\label{eq:T}
\end{equation}
where 
\begin{equation}
\mat{V}(\vc{q})=\left[\begin{array}{cc}\vc{q}^T\vc{q}^* &
\partial_x\vc{q}^T\\\partial_x\vc{q}^* & -\vc{q}^*\vc{q}^T
\end{array}\right]\,.
\label{eq:matrixV}
\end{equation}
Here $\vc{q}$ denotes the column vector of the functions $q_k(x,t)$
and $\mat{0}$ denotes the $N\times N$ zero matrix.  Along with these two
operators, we consider their nullspaces, $K_\op{X}(\lambda,\vc{q})\subset
A$ and $K_\op{T}(\lambda,\vc{q})\subset A$.  For generic $\lambda$, these
subspaces are $N+1$-dimensional, and if restricted to generic fixed $x$, $t$,
and $\lambda$ span $\C^{N+1}$.

If the functions $q_k(x,t)$ are chosen just right, then the subspaces
$K_\op{X}(\lambda,\vc{q})$ and $K_\op{T}(\lambda,\vc{q})$ may coincide
for all complex $\lambda$: $K_\op{X}=K_\op{T}\doteq K$.  If this is
the case, then the common nullspace will certainly be contained in the
nullspace of the commutator: $K\subset K_{[\op{X},\op{T}]}$.  As is
easily checked, the commutator $[\op{X},\op{T}]$ is not a
differential operator, but is merely a matrix multiplication operator,
with entries depending on $x$ and $t$ through the $q_k(x,t)$.  Since
the kernel of the commutator contains a subspace $K$ of dimension
$N+1$ for most $\lambda$, $x$, and $t$, this implies that the
operators $\op{X}$ and $\op{T}$ commute.  It is easily checked that the
compatibility condition $[\op{X},\op{T}]=\mat{0}$ is equivalent to
the vector nonlinear Schr\"odinger equation:
\begin{equation}
i\partial_t\vc{q}+\frac{1}{2}\partial_x^2\vc{q}+(\vc{q}^T\vc{q}^*)\vc{q}=0
\,.
\label{eq:vecNLS}
\end{equation}

It is therefore necessary that (\ref{eq:vecNLS}) be satisfied by the
functions $q_k(x,t)$ if we are to have a basis of simultaneous
nullvectors in the common nullspace $K$.  When they exist, we can
collect all these linearly independent column vectors into a square
matrix $\mat{F}(x,t,\lambda)$.  These ideas admit a natural geometric
interpretation in the trivial frame bundle $E\rightarrow \R^2$ with
fiber $GL(N+1,\C)$.  Here, $\op{X}$ and $\op{T}$ are covariant
derivative operators for $E$ in the $x$ and $t$ directions, and the
condition $[\op{X},\op{T}]= \mat{0}$ means that the curvature of the
affine connection specified by $\op{X}$ and $\op{T}$ is zero.  This
implies the existence of parallel global sections
$\mat{F}(x,t,\lambda)$ of the bundle $E$.

Finding a global section $\mat{F}(x,t,\lambda)$ of $E$ given
$\vc{q}(x,t)$, (that is, a matrix of simultaneous solution vectors) is
not always easy and for this reason, we will adopt a different point
of view below.  However, it is clear from (\ref{eq:X}) and
(\ref{eq:T}) that, given bounded functions $q_k(x,t)$ satisfying
(\ref{eq:vecNLS}) it is possible to develop an asymptotic expansion
for $\mat{F}(x,t,\lambda)$ in the limit $\lambda\rightarrow \infty$.
The expansion may be sought in the form:
\begin{equation}
\mat{F}(x,t,\lambda)=\left(c\bI_{N+1}+\lambda^{-1}\mat{F}^{(-1)}(x,t) +
\lambda^{-2}\mat{F}^{(-2)}(x,t)+\dots\right)\left[\begin{array}{cc}
e^{-2i(\lambda x + \lambda^2 t)} & \vc{0}^T\\\\
\vc{0} & \bI_N\end{array}\right]\,.
\label{eq:asymptotics}
\end{equation}
Here, $\bI_D$ denotes the $D\times D$ identity matrix, and $c$ is a
complex constant.  The coefficient matrices $\mat{F}^{(p)}(x,t)$ are
determined recursively in terms of $q_1(x,t),\dots,q_N(x,t)$ and the
constant $c$ by collecting powers of $\lambda$ in the compatible
equations $\op{X}\mat{F}=\op{T}\mat{F}=0$.  There is some ambiguity in this
expansion procedure entering as integration constants at each order.
However, it is easy to see that
\begin{equation}
F^{(-1)}_{1,k+1}(x,t)=\frac{c}{2i} q_k(x,t)\,,\hspace{0.3
in}k=1,\dots,N\,.
\label{eq:funcdef}
\end{equation}
regardless of the values of the integration constants.

The implications of this compatible structure for linear Schr\"odinger
equations that are of interest to us in this paper are easily stated.
\begin{prop}
Suppose that (\ref{eq:vecNLS}) is satisfied, and let
$\vc{v}(x,t,\lambda)\in K$ be any simultaneous nullvector of
$\op{X}(\lambda,\vc{q})$ and $\op{T} (\lambda,\vc{q})$.  Let
$\op{P}:A\rightarrow \C(x,t)$ be the operator of projection onto the
first component.  Define the {\em self-consistent potential}
\begin{equation}
V_0(x,t)\doteq -\vc{q}(x,t)^T\vc{q}(x,t)^*\,,
\label{eq:self-consistent}
\end{equation}
and set $f(x,t,\lambda)=\op{P}\vc{v}(x,t,\lambda)$.  Then it follows that
\begin{equation}
i\partial_tf +\frac{1}{2}\partial_x^2f - V_0(x,t)f=0\,.
\label{eq:linear}
\end{equation}
\label{prop:solve}
\end{prop}
So, for each complex $\lambda$, the function $f(x,t,\lambda)$ is a
solution of the linear, time-dependent Schr\"odinger
equation with potential (\ref{eq:self-consistent}).
Solutions corresponding to different values of $\lambda$ are linearly
independent.  Given functions $q_k(x,t)$ satisfying the nonlinear
system (\ref{eq:vecNLS}), one can look to the common nullspace $K$ of
the linear operators $\op{X}(\lambda,\vc{q})$ and
$\op{T}(\lambda,\vc{q})$ as a source of many solutions of the linear
equation (\ref{eq:linear}).

\noindent{\bf Remark:} Let us try to put these facts in a larger
context, and incidentally give the proof of
Proposition~\ref{prop:solve}.  It is part of the lore of integrable
systems theory that linearized evolution equations connected with
integrable systems are solvable in terms of ``squared eigenfunctions''
coming from the auxiliary linear problems making up the Lax pair for
the integrable system.  The integrable system (\ref{eq:vecNLS}) is the
compatibility condition for the equations ${\cal X}{\bf F}={\bf 0}$
and ${\cal T}{\bf F}={\bf 0}$.  By a change of variables (gauge
transformation) ${\bf F}={\bf G}\exp(-i\lambda x -i\lambda^2 t)$ the
two equations take the more familiar form of the Lax pair for
(\ref{eq:vecNLS}) \cite{Manakov,FT}:
\begin{equation}
\partial_x{\bf G} = {\bf A}{\bf G}\,,\hspace{0.3 in}
\partial_t{\bf G} = {\bf B}{\bf G}\,,
\label{eq:LaxPair}
\end{equation}
where
\begin{equation}
{\bf A} = \left[\begin{array}{cc}-i\lambda & \vec{q}^T\\-\vec{q}^* &
i\lambda{\mathbb I}
\end{array}\right]\,,\hspace{0.2 in}
{\bf B} = \left[\begin{array}{cc}-i\lambda^2 + i\vec{q}^T\vec{q}^*/2 &
\lambda\vec{q}^T + i\partial_x\vec{q}^T/2 \\
-\lambda\vec{q}^* + i\partial_x\vec{q}^*/2 & i\lambda^2{\mathbb I}-i\vec{q}^*\vec{q}^T/2
\end{array}\right]\,.
\end{equation}
If ${\bf G}_\alpha$ and ${\bf G}_\beta$ are any two simultaneous matrix
solutions of the Lax pair (\ref{eq:LaxPair}), and if $\bf C$ is any
constant (that is, $x$ and $t$ independent) matrix, then by setting
${\bf Q}={\bf G}_\alpha{\bf C} {\bf G}_\beta^{-1}$, one easily obtains the
equations
\begin{equation}
\partial_x{\bf Q}=[{\bf A},{\bf Q}]\,,\hspace{0.3 in}
\partial_t{\bf Q}=[{\bf B},{\bf Q}]\,.
\label{eq:QuadLaxPair}
\end{equation}
Equations of this form are called {\em Lax equations}, and the
elements of $\bf Q$ are the ``squared eigenfunctions''.  The
terminology becomes accurate in the scalar case $N=1$ when $\bf A$ and
$\bf B$ are in the Lie algebra $sl(2)$.  In this case the solutions
$\bf G$ of the Lax pair can be normalized to be in the Lie group
$SL(2)$ and therefore have determinant one.  Then, because ${\bf
G}_\beta$ is $2\times 2$ with determinant one, the elements of $\bf Q$
are seen to be {\em bona fide} quadratic forms in the solutions of the
Lax pair (\ref{eq:LaxPair}).  The emphasis in the literature on the
$sl(2)$-specific terminology of ``squared eigenfunctions'' for the
forms that satisfy the Lax equations (\ref{eq:QuadLaxPair}) no doubt
bears witness to the fact that so many of the famous integrable
equations (e.g. Korteweg-de Vries, scalar nonlinear Schr\"odinger,
sine-Gordon) are associated with $sl(2)$ representations.

If one introduces the splitting of a matrix into blocks: ${\bf M}={\bf
M}^{\rm D} + {\bf M}^{\rm OD}$ where ${\bf M}^{\rm D}$ consists of the
$1\times 1$ and $N\times N$ diagonal blocks of $\bf M$, and ${\bf
M}^{\rm OD}$ consists of the $1\times N$ and $N\times 1$ off-diagonal
blocks of $\bf M$, and if one introduces ${\bf A}_0 = {\bf
A}|_{\lambda=0}$ and ${\bf B}_0 = {\bf B}|_{\lambda=0}$, then it is an
exercise to check that the equations (\ref{eq:QuadLaxPair}) imply
\begin{equation}
\left[\begin{array}{cc}i & \vec{0}^T\\\vec{0}&-i{\mathbb I}\end{array}
\right]\partial_t{\bf Q}^{\rm OD} + \frac{1}{2}\partial_x^2{\bf Q}^{\rm OD} -
\left[\begin{array}{cc}i & \vec{0}^T\\\vec{0}&-i{\mathbb I}\end{array}
\right][{\bf B}_0^{\rm D},{\bf Q}^{\rm OD}]-\frac{1}{2}[
{\bf A}_0^{\rm OD},[{\bf A}_0^{\rm OD},{\bf Q}^{\rm OD}]] = {\bf 0}\,.
\end{equation}
If one writes
\begin{equation}
{\bf Q}^{\rm OD}=\left[\begin{array}{cc}0 & \vec{g}^T \\\vec{h} & {\bf 0}
\end{array}\right]\,,
\end{equation}
then one finds 
\begin{equation}
\begin{array}{rcccccccccl}
i\partial_t \vec{g}^T &+&\displaystyle\frac{1}{2}\partial_x^2 \vec{g}^T &+&
\vec{q}^T\vec{q}^*\vec{g}^T &+& \vec{q}^T\vec{h}\vec{q}^T &+&
\vec{g}^T\vec{q}^*\vec{q}^T &=&\vec{0}\\\\
-i\partial_t\vec{h} &+&\displaystyle\frac{1}{2}\partial_x^2\vec{h} &+&
\vec{q}^*\vec{q}^T\vec{h}&+&\vec{q}^*\vec{g}^T\vec{q}^*&+&
\vec{h}\vec{q}^T\vec{q}^*&=&\vec{0}\,.
\end{array}
\label{eq:complexlinearized}
\end{equation}
These linear equations for $\vec{g}$ and $\vec{h}$ are consistent with
the constraint $\vec{h}=\vec{g}^*$ at which point it becomes clear
that $\vec{g}(x,t)$ satisfies the linearization of the vector nonlinear
Schr\"odinger equation (\ref{eq:vecNLS}) about a solution
$\vec{q}(x,t)$.

Consider now a particular solution $\vec{q}(x,t)$ of (\ref{eq:vecNLS})
and by adjoining a new trivial component $q_{N+1}(x,t)\equiv 0$, view
it as a solution $\vec{q}'(x,t)$ of (\ref{eq:vecNLS}) in the $N+1$
component case.  From (\ref{eq:complexlinearized}) it is easily seen
that the corresponding components $g_{N+1}(x,t)$ and $h_{N+1}(x,t)$
satisfy
\begin{equation}
\begin{array}{rcccccl}
i\partial_tg_{N+1} &+&\displaystyle\frac{1}{2}\partial_x^2g_{N+1}&+&(\vec{q}^T\vec{q}^*)g_{N+1}&=&0\,,\\\\
-i\partial_th_{N+1}&+&\displaystyle\frac{1}{2}\partial_x^2h_{N+1}&+&(\vec{q}^T\vec{q}^*)h_{N+1}&=&0\,,
\end{array}
\end{equation}
where we have used the fact that
$\vec{q}'^T\vec{q}'^*=\vec{q}^T\vec{q}^*$.  Now considering the Lax
pair (\ref{eq:LaxPair}) for the primed potentials, it is easy to see
that there exists a nontrivial column vector solution of both
equations of the form $(\vec{0}_{N+1}^T,\exp(i\lambda x+i\lambda^2
t))^T$, and that further column vector solutions can then be chosen to
have a vanishing last component.  Taking the last column of the matrix
solution ${\bf G}_\beta$ to be this particular solution, and the first
$N+1$ columns all to have zeros in the final component, we see that
${\bf G}_\beta$ may be inverted in two independent blocks, and
therefore a solution of the linearized equation is given by
\begin{equation}
g_{N+1}(x,t)=Q_{1,N+2}(x,t,\lambda)=\exp(-i\lambda x-i\lambda^2 t)\sum_{k=1}^{N+2} C_{k,N+2}(\lambda)
G_{\alpha,1,k}(x,t,\lambda)\,.
\end{equation}
Since the matrix $\bf C$ is arbitrary, we may view the sum above as
the first component of an arbitrary column vector solution of the Lax
pair (\ref{eq:LaxPair}) with the primed potentials $\vec{q}'(x,t)$.
Moreover, since $q_{N+1}'(x,t)\equiv 0$, the first component of a
solution of the primed Lax pair is also the first component of a
solution of the unprimed Lax pair for the fully nontrivial potential
$\vec{q}(x,t)$.  Reversing the gauge transformation between solutions
$\bf G$ of the unprimed Lax pair (\ref{eq:LaxPair}) and solutions $\bf F$ of
${\cal X}{\bf F}={\cal T} {\bf F}={\bf 0}$ then establishes the
connection with Proposition~\ref{prop:solve}.

So, the procedure we are using for solving the time-dependent linear
Schr\"odinger equation is exactly the ``squared eigenfunction''
linearization of a certain $N+1$ component nonlinear Schr\"odinger
equation about a particular solution having $q'_{N+1}(x,t)\equiv 0$.
The ``squared eigenfunctions'' solving the linearized problem appear
to be linear in this special case because for $q'_{N+1}\equiv 0$ the
primed Lax pair becomes partly trivial, and the contribution of this
trivial part to the matrix $\bf Q$ is completely explicit (the
exponential function that we remove with a gauge transformation).
\pfend

We now return to the construction of self-consistent potentials and
the corresponding solutions of (\ref{eq:linear}).  The nonlinear
equation (\ref{eq:vecNLS}) is an integrable system by virtue of its
representation as the compatibility condition of two linear problems.
So there are many well-known ways to find functions $q_k(x,t)$ for
which the corresponding linear Schr\"odinger equation can be solved.
But as we are interested as much in the common nullspace of $\op{X}$
and $\op{T}$ as in the functions $q_k(x,t)$, we will now describe an
effective approach to finding both at the same time.  In this
approach, the object of fundamental importance is the common nullspace
$K$ itself.  We construct it first, with the functions $q_k(x,t)$
being chosen after the fact precisely so that for any basis matrix
$\mat{F}$ of $K$, we will have
$\op{X}(\lambda,\vc{q})\mat{F}=\op{T}(\lambda,\vc{q})\mat{F}=\mat{0}$.

What we know about $K$ is that whenever it exists by virtue of the
compatibility condition, the assumption that the functions $q_k(x,t)$
are bounded (this will be justified below) leads to expansions for
large $\lambda$ of a basis for $K$ of the form (\ref{eq:asymptotics}).
These expansions are generally only asymptotic; there is no guarantee
that there exists a choice of the integration constants such that the
expansion (\ref{eq:asymptotics}) converges for any $\lambda$ at all.
However, we now suppose that there exist solutions $q_k(x,t)$ of the
nonlinear system (\ref{eq:vecNLS}) for which an expansion
(\ref{eq:asymptotics}) not only converges in some deleted neighborhood
of $\lambda=\infty$, {\em but actually truncates}.  For such solutions
$q_k(x,t)$, if they exist, a basis of the subspace $K$ is given {\em
exactly} by an expression of the form
\begin{equation}
\mat{F}(x,t,\lambda)=\left(c\lambda^M\bI_{N+1}+
\sum_{p=0}^{M-1}\lambda^p\mat{F}^{(p)}(x,t)\right)
\left[\begin{array}{cc}\displaystyle e^{-2i(\lambda x + \lambda^2 t)} & \vc{0}^T \\\\
\vc{0} & \bI_N\end{array}\right]\,,
\label{eq:ansatz}
\end{equation}
for some positive integer $M$, where $c$ is a complex constant.  We
have multiplied by an explicit factor of $\lambda^M$ to bring the sum
into polynomial form.

Since we are not considering the functions $q_k(x,t)$ to be known, we
do not have the option of solving for the coefficient matrices
$\mat{F}^{(p)}(x,t)$ by substitution into the equations $\op{X}\mat{F}=
\op{T}\mat{F}=0$.  We therefore must consider them to be arbitrary
functions of $x$ and $t$ until we know otherwise.  Without any
constraints on the coefficients, we see that the differentiable matrix
functions of $x$, $t$, and $\lambda$ of the form (\ref{eq:ansatz}),
for given integer values of $M$ and $N$, form a vector space
$\Lambda_{N,M}$ over the complex numbers.

The space $\Lambda_{N,M}$ is very large.  If our claim --- that
appropriate solutions $q_k(x,t)$ of the nonlinear system
(\ref{eq:vecNLS}) exist --- is not vacuous, then $\Lambda_{N,M}$
should contain many proper subspaces that may be identified with the
common nullspace $K$ of $\op{X}(\lambda,\vc{q})$ and
$\op{T}(\lambda,\vc{q})$ for some $\vc{q}$.  If $\mat{F}(x,t,\lambda)$
is of the form (\ref{eq:ansatz}) and is a basis matrix of one of these
subspaces, then it must be determined modulo the constant $c$.  This
means that each such subspace of $\Lambda_{N,M}$ should ultimately be
isomorphic to $\C$, with the isomorphism being established via the
constant $c$.

We prepare to isolate the appropriate subspaces of $\Lambda_{N,M}$ by
defining a set of {\em discrete data}.  Let ${\cal D}$ denote an
$M$-tuple of pairs $(\lambda_k,\vc{g}^{(k)})$ where the $\lambda_k$
are distinct numbers in the complex upper half-plane and where the
$\vc{g}^{(k)}$ are vectors in $\C^N$.  From each vector
$\vc{g}^{(k)}$, we build $N+1$ vectors in $\C^{N+1}$ by setting
\begin{equation}
\vc{a}^{(k)}=(-1,g^{(k)*}_1,g^{(k)*}_2,\dots,g^{(k)*}_N)^T\in \C^{N+1}\,,
\end{equation}
and for $j=1,\dots,N$,
\begin{equation}
\vc{b}^{(k,j)}=(g^{(k)}_j,\vc{e}_j^T)^T\in\C^{N+1}\,,
\end{equation}
where $\vc{e}_j$ are the usual unit vectors in $\C^N$.

\begin{definition}
$\Lambda^{\cal D}_{N,M}$ is the subspace of $\Lambda_{N,M}$ whose
elements $\mat{F}(x,t,\lambda)$ satisfy the homogeneous linear conditions
\begin{equation}
\mat{F}(x,t,\lambda_k)\vc{a}^{(k)}=\vc{0}\,,
\label{eq:homog1}
\end{equation}
for $k=1,\dots,M$ and 
\begin{equation}
\mat{F}(x,t,\lambda_k^*)\vc{b}^{(k,j)}=\vc{0}\,,
\label{eq:homog2}
\end{equation}
for $k=1,\dots,M$ and $j=1,\dots,N$.
\end{definition}

It is not hard use dimension counting arguments to prove the following:
\begin{prop}
Let the discrete data $\cal D$ be given.  The set of solutions of
(\ref{eq:homog1}) and (\ref{eq:homog2}) forms a one-dimensional linear
subspace of $\Lambda_{N,M}$.  The general solution of
(\ref{eq:homog1}) and (\ref{eq:homog2}) is given by the one-parameter
family of matrices (\ref{eq:ansatz}), indexed by the complex parameter
$c$.  Thus, $\Lambda^{\cal D}_{N,M} \simeq \C$, with the isomorphism
being established via the complex constant $c$.  In particular, if $c$
is given, then the coefficient functions $\mat{F}^{(p)}(x,t)$ are
uniquely determined as functions of $x$ and $t$, and if $c$ is chosen
to be zero, then $\mat{F}(x,t,\lambda)$ is the zero matrix.
\end{prop}

This proposition allows us to index the elements of $\Lambda^{\cal
D}_{N,M}$ by the constant $c$ which is now a genuine coordinate for
the one-dimensional subspace $\Lambda^{\cal D}_{N,M}$.  We indicate
the dependence by writing $\mat{F}_{{\cal D},c}(x,t,\lambda)$ for the
matrices in this subspace.

This proposition is true even if homogeneous constraints less
structured than (\ref{eq:homog1}) and (\ref{eq:homog2}) are imposed.
In order for the dimension count to come out right it is sufficient to
choose $M\cdot (N+1)$ arbitrary complex numbers $\lambda_k$ along with
corresponding constant vectors $\vc{c}^{(k)}\in\C^{N+1}$ (the numbers
$\lambda_k$ need not all be distinct, as long as the vectors $\vc{c}$
belonging to each $\lambda_k$ are linearly independent) and to impose
$\mat{F}_{{\cal D},c}(x,t,\lambda_k)\vc{c}^{(k)}=\vc{0}$ for all
$k=1,\dots,M\cdot(N+1)$.  The additional structure in the constraints
(\ref{eq:homog1}) and (\ref{eq:homog2}) is needed for the following.

\begin{prop}
Let discrete data $\cal D$ be given, and let $\mat{F}_{{\cal D},c}(x,t,\lambda)
\in\Lambda^{\cal D}_{N,M}$.  Then
\begin{equation}
\frac{1}{c}F^{(M-1)}_{k+1,1}(x,t)=-\left(\frac{1}{c}F^{(M-1)}_{1,k+1}(x,t)\right)^*\,,\hspace{0.3
in}k=1,\dots,N\,.
\end{equation}
\label{prop:reality}
\end{prop}
We will have use for this symmetry property below.  Its proof is simple.

\noindent{\bf Proof of Proposition~\ref{prop:reality}:} It is
sufficient to consider the case of $c=1$, since the coefficient
matrices simply scale with $c$.  It will be convenient to introduce
the block form of the coefficient matrices:
\begin{equation}
\mat{F}^{(p)}(x,t)=\left[\begin{array}{cc} a^{(p)} & \vc{b}^{(p)T}\\
\vc{c}^{(p)} & \mat{D}^{(p)}\end{array}\right]\,,
\end{equation}
where $a^{(p)}(x,t)$ is a scalar, $\vc{b}^{(p)}(x,t)$ and
$\vc{c}^{(p)}(x,t)$ are $N$-component vectors, and
$\mat{D}^{(p)}(x,t)$ is an $N\times N$ matrix.  We will prove the
stronger result that for all $p=0,\dots,M-1$,
\begin{equation}
\vc{c}^{(p)}(x,t)=-\vc{b}^{(p)}(x,t)^*\,.
\end{equation}
In this form, the equations (\ref{eq:homog1}) and (\ref{eq:homog2}) take
the form of the system:
\begin{equation}
\begin{array}{rcl}\displaystyle
-\lambda_k^M +
\sum_{p=0}^{M-1}\lambda_k^p\left(\vc{b}^{(p)T}\vc{g}^{(k)*}
e^{2i(\lambda_k x +\lambda_k^2t)}-a^{(p)}\right)&=&0\,,\\\\
\displaystyle
\lambda_k^M\vc{g}^{(k)*}+\sum_{p=0}^{M-1}\lambda_k^p\left(\mat{D}^{(p)}
\vc{g}^{(k)*}-e^{-2i(\lambda_kx+\lambda_k^2t)}\vc{c}^{(p)}\right)&=&\vc{0}\,,\\\\
\displaystyle
\lambda_k^{*M}\vc{g}^{(k)}+\sum_{p=0}^{M-1}\lambda_k^{*p}\left(a^{(p)}\vc{g}^{(k)}+e^{2i(\lambda_k^*x+\lambda_k^{*2}t)}\vc{b}^{(p)}\right)&=&\vc{0}\,,\\\\
\displaystyle
\lambda_k^{*M}\bI_{N}+\sum_{p=0}^{M-1}\lambda_k^{*p}\left(e^{-2i(\lambda_k^*x+\lambda_k^{*2}t)}\vc{c}^{(p)}\vc{g}^{(k)T}+\mat{D}^{(p)}\right)&=&\mat{0}\,,
\end{array}
\end{equation}
where $k=1,\dots,M$.  From the first and fourth equations, we can
eliminate $a^{(p)}(x,t)$ and $\mat{D}^{(p)}(x,t)$, $p=0,\dots,M-1$ in
favor of the $\vc{b}^{(p)}(x,t)$ and $\vc{c}^{(p)}(x,t)$.  This
involves introducing the elements of the inverse $\mat{W}$ of the
Vandermonde matrix $\mat{V}$ having elements
$V_{jk}\doteq\lambda_j^{k-1}$, but it leads to two decoupled linear
systems, one for the $\vc{b}^{(p)}(x,t)$ and the other for the
$\vc{c}^{(p)}(x,t)$.  These systems are:
\begin{equation}
\begin{array}{rcl}\displaystyle
\sum_{r=1}^MH_{kr}\vc{c}^{(r-1)}&=&h_k\,,\\\\\displaystyle
\sum_{r=1}^MH^*_{kr}\vc{b}^{(r-1)}&=&-h_k^*\,,
\end{array}
\end{equation}
where
\begin{equation}
H_{kr}\doteq
V_{kr}e^{-2i(\lambda_kx+\lambda_k^2t)}+\sum_{s=1}^MV_{ks}\sum_{j=1}^MW_{sj}^*
V_{jr}^*e^{-2i(\lambda_j^*x+\lambda_j^{*2}t)}\vc{g}^{(j)T}\vc{g}^{(k)*}\,,
\end{equation}
and
\begin{equation}
h_k\doteq
\lambda_k^M\vc{g}^{(k)*}-\sum_{s=1}^M\sum_{j=1}^MV_{ks}W_{sj}^*\lambda_j^{*M}
\vc{g}^{(k)*}\,.
\end{equation}
It is then clear that $\vc{c}^{(p)}=-\vc{b}^{(p)*}$ for all $p$.
\pfend

So, the emphasis has changed with respect to these matrices and their
relation to the functions $q_k(x,t)$.  Rather than determining the
coefficient matrices $\mat{F}^{(p)}(x,t)$ from a given set of
functions $q_k(x,t)$ solving (\ref{eq:vecNLS}) by an asymptotic
expansion procedure, we are determining them from the discrete data
$\cal D$ and a choice of the constant $c$.  If there is to be any
consistency, then we must still have relations between the coefficient
matrices $\mat{F}^{(p)}(x,t)$ of $\mat{F}_{{\cal D},c}$ and the
functions $q_k(x,t)$; in particular, we can rewrite
(\ref{eq:funcdef}):
\begin{equation}
q_k(x,t)\doteq \frac{2i}{c}F^{(M-1)}_{1,k+1}(x,t)
\label{eq:funcdef2}
\end{equation}
and use it as a {\em definition} of some functions $q_k(x,t)$ in terms
of the discrete data $\cal D$ and the constant $c$.  Note that as long
as $c\neq 0$, then this definition is actually independent of $c$
because $\mat{F}_{{\cal D},c}$ is directly proportional to $c$.  The
fact that (\ref{eq:funcdef2}) is sensible as a definition of the
$q_k(x,t)$ is shown by:
\begin{prop}
Let the discrete data $\cal D$ be given and let the constant $c$ be
nonzero, and let the functions $q_1(x,t),\dots, q_N(x,t)$ be defined
(in terms of $\cal D$ alone) via (\ref{eq:funcdef2}).  This determines
the operators $\op{X}(\lambda,\vc{q})$ and $\op{T}(\lambda, \vc{q})$.
Then for any $\mat{F}_{{\cal D},c}(x,t,\lambda)\in\Lambda^{\cal
D}_{N,M}$,
\begin{equation}
\op{X}(\lambda,\vc{q})\mat{F}_{{\cal D},c}(\lambda,x,t)=
\op{T}(\lambda,\vc{q})\mat{F}_{{\cal D},c}(\lambda,x,t)=\mat{0}\,.
\end{equation}
For these $q_k(x,t)$ the columns of $\mat{F}_{{\cal
D},c}(x,t,\lambda)$ are generically linearly independent and therefore
form a basis of the common nullspace $K$ for almost all $\lambda$.
\label{keyprop}
\end{prop}

Recall that the commutator $[\op{X},\op{T}]$ is, for fixed $x$ and
$t$, a matrix multiplication operator in $\C^{N+1}$.  Thus, the
existence of the common nullspace $K$ of $\op{X}(\lambda,\vc{q})$
and $\op{T}(\lambda,\vc{q})$ of generic dimension $N+1$ for these
functions $q_k(x,t)$ implies the vanishing of the commutator and the
compatibility of the two linear problems.  Therefore we have
\begin{cor}
The functions $q_k(x,t)$ constructed from any set of discrete data
$\cal D$ satisfy the vector nonlinear Schr\"odinger equation (\ref{eq:vecNLS}).
\end{cor}
A time-dependent self-consistent potential function $V_0(x,t)$
generated from the functions $q_k(x,t)$ connected with a set of
discrete data $\cal D$ according to (\ref{eq:self-consistent}) will be
called a {\em separable potential} \cite{physd}.

\noindent{\bf Proof of Proposition~\ref{keyprop}:} Let $\mat{F}_{{\cal
D},c}(x, t,\lambda)\in\Lambda^{\cal D}_{N,M}$.  The proof begins with
the simple observation that, as a consequence of the vectors
$\vc{a}^{(k)}$ and $\vc{b}^{(k,j)}$ in the homogeneous relations
(\ref{eq:homog1}) and (\ref{eq:homog2}) satisfied by $\mat{F}_{{\cal
D},c}(x,t,\lambda)$ being independent of $x$ and $t$, these relations
are satisfied by $(\op{X}\mat{F}_{{\cal D},c})(x,t,\lambda)$ and
$(\op{T}\mat{F}_{{\cal D},c})(x,t,\lambda)$ as well.  For example,
with the operator $\op{X}$,
\begin{equation}
(\op{X}(\lambda,\vc{q})\mat{F}_{{\cal
D},c})(x,t,\lambda_k)\vc{a}^{(k)}=\op{X}(\lambda_k,\vc{q})(\mat{F}_{{\cal
D},c}(x,t,\lambda_k)\vc{a}^{(k)})=
\op{X}(\lambda_k,\vc{q})\mat{0}=\mat{0}\,,
\end{equation}
for $k=1,\dots,M$,
and
\begin{equation}
(\op{X}(\lambda,\vc{q})\mat{F}_{{\cal
D},c})(x,t,\lambda_k^*)\vc{b}^{(k,j)}=\op{X}(\lambda_k^*,\vc{q})(\mat{F}_{{\cal
D},c}(x,t,\lambda_k^*)\vc{b}^{(k,j)})=\op{X}(\lambda_k^*,\vc{q})\mat{0}=\mat{0}\,,
\end{equation}
for $k=1,\dots,M$ and $j=1,\dots,N$.  The argument is unchanged if
$\op{X}$ is replaced with $\op{T}$.

Next, we examine the form of the matrix
$(\op{X}(\lambda,\vc{q})\mat{F}_{{\cal D},c})(x,t,\lambda)$.  It is
straightforward to see that
\begin{equation}
\begin{array}{rcl}
\displaystyle
(\op{X}(\lambda,\vc{q})\mat{F}_{{\cal D},c})(x,t,\lambda)&=&\displaystyle
\Bigg\{\lambda^M\left(2i[\mat{E},\mat{F}^{(M-1)}]-c\mat{U}\right
) +\\\\
&&\displaystyle
\sum_{p=1}^{M-1}\lambda^p\left(
\partial_x\mat{F}^{(p)}+2i[\mat{E},\mat{F}^{(p-1)}]-\mat{U}
\mat{F}^{(p)}\right) + \\\\
&&\displaystyle
\left(\partial_x\mat{F}^{(0)}-\mat{U}\mat{F}^{(0)}\right)
\Bigg\}\exp(-2i(\lambda x+\lambda^2 t)\mat{E})\,.
\end{array}
\label{eq:XF}
\end{equation}
Now, as a consequence of the definition of the functions $q_k(x,t)$ in
terms of the discrete data $\cal D$ and the symmetry property
guaranteed by Proposition~\ref{prop:reality}, the leading term
vanishes identically, that is,
\begin{equation}
2i[\mat{E},\mat{F}^{(M-1)}(x,t)]=c\mat{U}(\vc{q})\,.
\end{equation}
This, along with the fact that $\op{X}(\lambda,\vc{q})\mat{F}_{{\cal
D},c}(x,t,\lambda)$ satisfies the homogeneous conditions
(\ref{eq:homog1}) and (\ref{eq:homog2}) means that
\begin{equation}
\op{X}(\lambda,\vc{q})\mat{F}_{{\cal
D},c}(x,t,\lambda)\in\Lambda^{\cal D}_{N,M}\,.
\end{equation}
Not only that, but for matrices in $\Lambda^{\cal D}_{N,M}$ the only
way that the coefficient of $\lambda^M$ can vanish is for the leading
constant to vanish.  Therefore, by the isomorphism between
$\Lambda^{\cal D}_{N,M}$ and $\C$ via the leading constant, it follows
that
\begin{equation}
\op{X}(\lambda,\vc{q})\mat{F}_{{\cal D},c}(x,t,\lambda)=\mat{0}\,.
\end{equation}

We now consider the form of $\op{T}(\lambda,\vc{q})\mat{F}_{{\cal
D},c}(x,t,\lambda)$:
\begin{equation}
\begin{array}{rcl}
\displaystyle \op{T}(\lambda,\vc{q})\mat{F}_{{\cal
D},c}(x,t,\lambda)&=&\displaystyle
\Bigg\{\lambda^{M+1}\left(2i[\mat{E},\mat{F}^{(M-1)}]-c\mat{U} \right)
+\\\\ &&\displaystyle \lambda^M\left(2i[\mat{E},\mat{F}^{(M-2)}]-
\mat{U} \mat{F}^{(M-1)}-\frac{ic}{2}\mat{V} \right)+\\\\
&&\displaystyle \sum_{p=2}^{M-1}\lambda^p\left(\partial_t\mat{F}^{(p)}
+2i[\mat{E},\mat{F}^{(p-2)}]-\mat{U}\mat{F}^{(p-1)}-
\frac{i}{2}\mat{V} \mat{F}^{(p)}\right)+\\\\ &&\displaystyle
\lambda\left(\partial_t\mat{F}^{(1)}-\mat{U} \mat{F}^{(0)}-
\frac{i}{2}\mat{V}\mat{F}^{(1)}\right)+\\\\ &&\displaystyle
\left(\partial_t\mat{F}^{(0)}-\frac{i}{2}\mat{V}\mat{F}^{(0)}\right)\Bigg\}\exp(-2i(\lambda
x+\lambda^2 t)\mat{E})\,.
\end{array}
\label{eq:TF}
\end{equation}
Once again, the definition of the functions $q_k(x,t)$ and the 
symmetry property of Proposition~\ref{prop:reality} guarantee that
the coefficient of $\lambda^{M+1}$ vanishes.  We shall now show
that the coefficient of $\lambda^M$ vanishes as well.  Begin by writing
$\mat{F}^{(M-1)}(x,t)$ in the block form:
\begin{equation}
\mat{F}^{(M-1)}(x,t)=\left[\begin{array}{cc}a^{(M-1)} & \vc{b}^{(M-1)T}
\\ \vc{c}^{(M-1)} &
\mat{D}^{(M-1)}\end{array}\right]\,.
\end{equation}
We already know by definition of the functions $q_k(x,t)$ and
Proposition~\ref{prop:reality}, that $\vc{b}^{(M-1)T}=\vc{q}^T/2i$
and $\vc{c}^{(M-1)}=-\vc{b}^{(M-1)*}$.  Making use of the fact that
all of the terms in (\ref{eq:XF}) vanish identically, we also have
\begin{equation}
\partial_x\mat{F}^{(M-1)}(x,t) + 2i[\mat{E},\mat{F}^{(M-2)}(x,t)]-
\mat{U}(\vc{q})\mat{F}^{(M-1)}(x,t)=\mat{0}\,.
\label{eq:alsohave}
\end{equation}
This implies that for the coefficient of $\lambda^M$ in (\ref{eq:TF})
to vanish, it will be enough to show that
\begin{equation}
c\mat{V}(\vc{q})=2i\partial_x\mat{F}^{(M-1)}(x,t)\,.
\label{eq:toshow}
\end{equation}
From (\ref{eq:matrixV}), it is clear (\ref{eq:toshow}) is satisfied in
the off-diagonal blocks.  To show that the diagonal blocks also vanish,
we write out the diagonal blocks of (\ref{eq:alsohave}):
\begin{equation}
\left[\begin{array}{cc}
\partial_x a^{(M-1)} & \vc{0}^T\\
\vc{0}&\partial_x\mat{D}^{(M-1)}\end{array}\right]-\left[\begin{array}{cc}
\vc{q}^T\vc{c}^{(M-1)} & \vc{0}^T\\
\vc{0} & -\vc{q}^*\vc{b}^{(M-1)T}\end{array}\right]=\mat{0}\,.
\end{equation}
Eliminating the derivatives of $a^{(M-1)}$ and $\mat{D}^{(M-1)}$
between this equation and the diagonal blocks of (\ref{eq:toshow}) and
comparing with the definition (\ref{eq:matrixV}) of $\mat{V}(\vc{q})$,
we finally see that (\ref{eq:toshow}) is satisfied identically.  By
similar arguments as we used above, it follows that
$\op{T}(\lambda,\vc{q})\mat{F}_{{\cal D},c}(x,t,\lambda)$ is also the
zero element of $\Lambda^{\cal D}_{N,M}$. This ends the proof of the
proposition.  \pfend

We now return to the problem of interest, namely the algebraic
construction of separable time-dependent potentials for the linear
Schr\"odinger equation and of a large number of exact solutions to
this linear equation.  From the construction of the subspace
$\Lambda_{N,M}^{\cal D}$ we can extract a simpler construction of the
quantities of immediate interest, and cast the whole procedure in the
form of an algorithm.  The key observation is that it is sufficient to
build from given discrete data $\cal D$ {\em only the first row} of a
matrix $\mat{F}_{{\cal D},c}(x,t,\lambda)$ in the space
$\Lambda_{N,M}^{\cal D}$.  This gives us both the functions $q_k(x,t)$
via the first row of the coefficient matrix $\mat{F}^{(M-1)}(x,t)$ from
which we find the potential $V_0(x,t)$ and also the image of the projection
operator $\op{P}$ that consists of solutions of the linear Schr\"odinger
equation with this potential.

So we consider the first row of $\mat{F}_{{\cal D},c}(x,t,\lambda)$
and impose the homogeneous linear constraints (\ref{eq:homog1}) and
(\ref{eq:homog2}).  Introducing
\begin{equation}
a(x,t,\lambda)=F_{11}(x,t,\lambda)=\left(\lambda^M +\sum_{p=0}^{M-1}\lambda^pa^{(p)}(x,t)\right)e^{-2i(\lambda x +\lambda^2 t)}\,,
\end{equation}
and
\begin{equation}
\vc{b}(x,t,\lambda)=(F_{12}(x,t,\lambda),\dots,F_{1,N+1}(x,t,\lambda))^T=
\sum_{p=0}^{M-1}\lambda^p\vc{b}^{(p)}(x,t)\,,
\end{equation}
the relations (\ref{eq:homog1}) and (\ref{eq:homog2}) take the simple
form:
\begin{equation}
\begin{array}{rcl}
a(x,t,\lambda_k)&=&\vc{g}^{(k)\dagger}\vc{b}(x,t,\lambda_k)\,,\\
\vc{b}(x,t,\lambda_k^*)&=&-a(x,t,\lambda_k^*)\vc{g}^{(k)}\,,
\end{array}
\end{equation}
where $k=1,\dots,M$.  Note that without loss of generality, we are
taking $c=1$.  Written out in its entirety, this is a square linear
system for the $M\cdot(N+1)$ unknowns, $a^{(p)}(x,t)$, and the $N$
elements of $\vc{b}^{(p)}(x,t)$ for $p=0,\dots,M-1$.  The matrix of
this system, and the right-hand side, are explicit functions of $x$
and $t$ through the exponential functions contributed by
$a(x,t,\lambda_k)$ and $a(x,t,\lambda_k^*)$.  

From the solution of this linear system, one computes the potential 
function as:
\begin{equation}
V_0(x,t)=-4\sum_{n=1}^N\Big|b_n^{(M-1)}(x,t)\Big|^2\,.
\end{equation}
Then, we see that $a(x,t,\lambda)$ and all the elements of 
$\vc{b}(x,t,\lambda)$ are solutions of the linear equation
\begin{equation}
i\partial_tf +\frac{1}{2}\partial_x^2f -V_0(x,t)f=0\,,
\label{eq:schrodagain}
\end{equation}
for fixed but arbitrary $\lambda$.  Being polynomial in $\lambda$,
each element of $\vc{b}(x,t,\lambda)$ sweeps out an $M$-dimensional
space of solutions as $\lambda$ varies.  The solutions contained in
$a(x,t,\lambda)$ are more interesting because the presence of the
exponential means that all of these solutions for real
$\lambda$ are linearly independent.  This immediately gives an
infinite-dimensional space of solutions to the linear Schr\"odinger
equation.

In fact, the function $a(x,t,\lambda)$ contains an $L^2(\R)$ basis of
solutions of the Schr\"odinger equation as the parameter $\lambda$ is
varied \cite{physd}.  In particular, the set \be
\{a(x,t,\lambda_1^*),\dots,a(x,t,\lambda_M^*),a(x,t,\lambda),
\lambda {\rm \,\,real}\} \label{eq:explicitbasis} \ee considered as
functions of $x$ for fixed $t$, is complete.  
For real $\lambda$, set
\begin{equation}
\Psi_{\rm
d}(x,t,\lambda)\doteq\left(\pi\prod_{k=1}^M|\lambda-\lambda_k|^2\right)^{
-1/2}a(x,t,\lambda)\,.
\end{equation}
The subscript ``d'' indicates solutions that superpose to form
dispersive waves.  For $\lambda$ and $\eta$ real we then have the
inner products
\begin{equation}
\langle\Psi_{\rm d}(\cdot,t,\lambda),\Psi_{\rm
d}(\cdot,t,\eta)\rangle=\delta(\lambda-\eta)\,,
\end{equation}
and for $k=1,\dots,M$,
\begin{equation}
\langle a(\cdot,t,\lambda_k^*),\Psi_{\rm d}(\cdot,t,\mu)\rangle =
0\,.
\end{equation}
Also, ${\rm dim\,\,\,span\,}\{a(x,t,\lambda_k^*),\,\,
k=1,\dots,M\}=M$ as functions of $x$ for fixed $t$.  So, let
$\{\Psi_{{\rm b},k}(x,t)\}$ be any basis of ${\rm
span\,}\{a(x,t,\lambda_k^*),\,\, k=1,\dots,M\}$ that is orthonormal
with respect to the inner product (say obtained by the Gram-Schmidt
procedure), so that
\begin{equation}
\langle\Psi_{{\rm b},j}(\cdot,t),\Psi_{{\rm
b},k}(\cdot,t)\rangle=\delta_{jk}\,.
\end{equation}
The subscript ``b'' indicates solutions that are bound and have finite
energy.  Note that this basis remains orthonormal because the time
evolution of these functions under (\ref{eq:schrodagain}) is unitary.

The completeness relation is generalized to $L^2(\R)$ from that proved
in \cite{physd} as:
\begin{prop}
Let discrete data $\cal D$ be given and let $t\in\R$ be fixed.  For 
all $f(x)\in L^2(\R)$, we have the expansion
\begin{equation}
f(x)=\int_{-\infty}^\infty f_{\rm d}(\lambda,t)\Psi_{\rm d}(x,t,\lambda)\,
d\lambda +\sum_{k=1}^M f_{{\rm b},k}(t)\Psi_{{\rm b},k}(x,t)\,,
\end{equation}
where the expansion coefficients are given by:
\begin{equation}
f_{\rm d}(\lambda,t)=\langle \Psi_{\rm d}(\cdot,t,\lambda),f(\cdot)\rangle\,,
\hspace{0.3 in}
f_{{\rm b},k}(t)=\langle \Psi_{{\rm b},k}(\cdot,t),f(\cdot)\rangle\,.
\end{equation}
\end{prop}
The orthogonality relations for the functions $\Psi_{{\rm b},k}(x,t)$
and $\Psi_{\rm d}(x,t,\lambda)$ are implied by this result.  Note that
if $f=f(x,t)$ satisfies (\ref{eq:schrodagain}) then the expansion
coefficients are independent of $t$ and can be constructed from the
initial data $f(x,0)$.  Thus one solves the initial value problem for
(\ref{eq:schrodagain}) in $L^2(\R)$.

\renewcommand{\thesection}{Appendix \Alph{section}}
\section{Dispersive Local Decay Estimates}
\label{app:localdecay}
\renewcommand{\thesection}{\Alph{section}} \setcounter{equation}{0}
Here, we establish several important properties of the unitary group
$e^{-it{\op{B}}}$.  We will consider even perturbations of the even
two-soliton periodic potentials, so we will work in either the even or
odd subspace of $L^2(\R)$.  For a given function $f(\cdot)$ in
$L^2(\R)$, the operator $\op{P}^{\rm (e,o)}_{\rm c}$ is defined as the
spectral projection onto the continuous part of the spectrum of
$\op{B}$:
\begin{equation}
(\op{P}^{\rm (e,o)}_{\rm c}f)(x)\doteq \int_0^\infty\langle\Psi^{\rm
(e,o)}_{\rm d}(\cdot,0,\lambda),f(\cdot)\rangle\Psi^{\rm (e,o)}_{\rm
d}(x,0,\lambda)\,d\lambda\,.
\end{equation}
As we will now see, the main difference between the even and odd cases
is in the rate of dispersive decay, and the difference can be directly
traced to the behavior of the dispersive eigenfunction $\Psi^{\rm
(e,o)}_{\rm d}(x,0,\lambda)$ in the vicinity of $\lambda=0$.  It is
easy to see from the explicit formulas that the eigenfunctions are
continuous in $\lambda$ at $\lambda=0$, and that the odd mode vanishes
there:
\begin{equation}
\Psi_{\rm d}^{\rm (o)}(x,0,\lambda=0)=0\,,
\end{equation}
while the even mode does not vanish, but is simply finite at
$\lambda=0$.  We say that the existence of a nontrivial eigenfunction
at $\lambda=0$, as in the even case, indicates a {\em zero-energy
resonance} of the system.  The ubiquitous effect of a zero-energy
resonance is to alter the rate of dispersive decay in the system.
However, more dramatic effects can appear if under the influence of a
perturbation, the zero-energy resonance is directly excited.  This
latter situation we refer to as a {\em parametric zero-energy
resonance}.  A system with a zero-energy resonance is ``primed'' to
feel the effects of a parametric zero-energy resonance in the presence
of an appropriate perturbation.

\subsection{Nonsingular local decay.}
First, we will prove the {\em nonsingular local decay estimate} for the
unitary group $e^{-it\op{B}}$.
\begin{prop}
\label{prop:nonsingular}
Fix $\sigma>5/2$.  There exist constants $L^{\rm (e,o)}>0$ such that
\begin{equation}
\|\langle\cdot\rangle^{-\sigma}\left(e^{-it\op{B}}\op{P}_{\rm c}^{\rm (e)}f\right)(\cdot)\|_2\le L^{\rm (e)}\langle t\rangle^{-1/2}\|\langle\cdot\rangle^{\sigma}f(\cdot)\|_2\,,
\end{equation}
and 
\begin{equation}
\|\langle\cdot\rangle^{-\sigma}\left(e^{-it\op{B}}\op{P}_{\rm c}^{\rm (o)}f\right)(\cdot)\|_2\le L^{\rm (o)}\langle t\rangle^{-3/2}\|\langle\cdot\rangle^{\sigma}f(\cdot)\|_2\,,
\end{equation}
for all $f\in L^2(\R)$ for which the right hand
side makes sense.
\end{prop}

The proof is based on a sequence of intermediate results.
First, from the simple chain of estimates:
\begin{equation}
\|\langle \cdot\rangle^{-\sigma}e^{-it\op{B}}\op{P}^{\rm (e,o)}_{\rm c} f(\cdot)\|_2  \le 
\|e^{-it\op{B}}\op{P}_{\rm c}^{\rm (e,o)}f(\cdot)\|_2 
 =
\|\op{P}_{\rm c}^{\rm (e,o)}f(\cdot)\|_2
\le  \|f(\cdot)\|_2 
\le 
\|\langle\cdot\rangle^\sigma f(\cdot)\|_2\,,
\end{equation}
we have
\begin{lemma}
\label{lemma:uniform}
For all $\sigma> 0$, we have the simple estimate
\begin{equation}
\|\langle \cdot\rangle^{-\sigma}e^{-it\op{B}}\op{P}^{\rm (e,o)}_{\rm c} f(\cdot)\|_2 \le\|\langle\cdot\rangle^\sigma f(\cdot)\|_2\,,
\end{equation}
for all $f\in L^2(\R)$ for which the right-hand side make sense.
\end{lemma}

We now want to refine the above uniform estimate to include a
multiplicative factor of $\|\langle\cdot\rangle^\sigma f(\cdot)\|_2$
that decays in $|t|$.  To this end, we fix $t\neq 0$ and observe that
by the definition of the operator $\op{B}$,
\begin{equation}
\label{eq:representation}
\langle x\rangle^{-\sigma}\left(e^{-it\op{B}}\op{P}_{\rm c}^{\rm
(e,o)}f\right)(x) =
\int_{-\infty}^\infty \langle y\rangle^\sigma f(y)
\,h(x,y;t)\,dy\,,
\end{equation}
where
\begin{equation}
h(x,y;t)\doteq
\langle x\rangle^{-\sigma}\langle y\rangle^{-\sigma}\int_0^\infty 
\Psi_{\rm d}^{\rm (e,o)}(y,0,\lambda)^*\Psi_{\rm d}^{\rm (e,o)}(x,0,\lambda)
e^{-2i\lambda^2 t}\,d\lambda\,.
\end{equation}
We note here that the integral in the definition of $h(x,y;t)$ is
improper; the integrand is not absolutely integrable, and the integral
from zero to infinity should be interpreted as the limit of the
integral from zero to $R$ as $R\uparrow\infty$.  This limit exists as
long as $t\neq 0$, and consequently the function $h(x,y;t)$ is
well-defined for $t\neq 0$.  The trouble with the function $h(x,y;t)$
at $t=0$ is not our concern here because we already have a uniform
estimate that holds for all $t$, and in particular for $t$ near zero.
Thus we will be thinking of $t$ as being large in what follows.

In any case, by Cauchy-Schwarz, we have
\begin{equation}
\left|\langle x\rangle^{-\sigma}\left(e^{-it\op{B}}\op{P}_{\rm c}^{\rm (e,o)}
f\right)(x)\right|\le 
\|h(x,\cdot;t)\|_2\,\,\|\langle\cdot\rangle^\sigma 
f(\cdot)\|_2\,.
\end{equation}
It follows that
\begin{equation}
\left\|\langle \cdot\rangle^{-\sigma}\left(e^{-it\op{B}}
\op{P}_{\rm c}^{\rm (e,o)}f\right)(\cdot)\right\|_2\le
\|h(\cdot,\cdot;t)\|_2\,\,\|\langle\cdot\rangle^\sigma f(\cdot)\|_2\,,
\end{equation}
an estimate that involves the Hilbert-Schmidt norm of the kernel
$h(x,y;t)$ for each fixed $t$.  The rest of our work will be to show 
$h(x,y;t)$ is in $L^2(\R^2)$ for each fixed $t$, with norm decaying in $|t|$.

First note that from the explicit formulas:
\begin{equation}
\Psi_{\rm d}^{\rm (o)}(x,0,\lambda)=
\frac{2\lambda a^{(1)}(x,0)\cos(2\lambda x)-2i(\lambda^2 + a^{(0)}(x,0))
\sin(2\lambda x)}{\sqrt{2\pi (\lambda^2 +\rho_1^2)(\lambda^2+\rho_2^2)}}\,,
\end{equation}
\begin{equation}
\Psi_{\rm d}^{\rm (e)}(x,0,\lambda)=
\frac{2(\lambda^2+a^{(0)}(x,0))\cos(2\lambda x)-2i\lambda a^{(1)}(x,0)
\sin(2\lambda x)}{\sqrt{2\pi (\lambda^2 +\rho_1^2)(\lambda^2 +\rho_2^2)}}\,,
\end{equation}
where $a^{(0)}(x,t)$ and $a^{(1)}(x,t)$ are bounded
analytic functions of $x$, we obtain
\begin{lemma}
Let the parameters $\rho_1$ and $\rho_2$ be fixed.  The function
defined by
\begin{equation}
q(\lambda)\doteq\Psi^{\rm (e,o)}_{\rm d}(x,0,\lambda)\Psi^{\rm
(e,o)}_{\rm d}(y,0,\lambda)^*
\end{equation}
is in $C^k(\R_+)$ for all $k\ge 0$.  In
particular all derivatives with respect to $\lambda$ are uniformly
bounded functions of $\lambda$.  The norms $\|q^{(k)}(\cdot)\|_\infty$
are homogeneous polynomials in $|x|$ and $|y|$ of degree $k$, with
nonnegative coefficients that depend only on $\rho_1$ and $\rho_2$.
Also, in the odd case, we have $q(\lambda)=O(\lambda^2)$ for $\lambda$ near
zero, while in the even case $q(\lambda)=O(1)$.
\label{lemma:qproperties}
\end{lemma}

In showing that $h(x,y;t)$ is $L^2(\R^2)$ with norm decaying in $t$,
we will find that the main contribution for large $t$ comes from the
part of the integral near $\lambda=0$.  To see this, we first separate
the contributions near and away from zero.  Let $g_\Delta(\lambda)$ be
a nonnegative ``bump function'', infinitely differentiable for
real $\lambda>0$, identically equal to $1$ for $0\le\lambda\le
\Delta/2$ and identically equal to zero for $\lambda\ge 3\Delta/2$.
Let $\tilde{g}_\Delta(\lambda)\doteq 1-g_\Delta(\lambda)$.  Then
\begin{equation}
h(x,y;t)=h_\Delta(x,y;t) + \tilde{h}_\Delta(x,y;t)\,,
\end{equation}
where
\begin{equation}
h_\Delta(x,y;t)\doteq \langle x\rangle^{-\sigma}\langle y\rangle^{-\sigma}
\int_0^{3\Delta/2}\Psi_{\rm d}^{\rm (e,o)}(y,0,\lambda)^*\Psi_{\rm d}^{\rm (e,o)}(x,0,\lambda)g_\Delta(\lambda)e^{-2i\lambda^2 t}\,d\lambda\,,
\end{equation}
and
\begin{equation}
\tilde{h}_\Delta(x,y;t)
\doteq \langle x\rangle^{-\sigma}\langle y\rangle^{-\sigma}
\int_{\Delta/2}^\infty\Psi_{\rm d}^{\rm (e,o)}(y,0,\lambda)^*\Psi_{\rm d}^{\rm (e,o)}(x,0,\lambda)\tilde{g}_\Delta(\lambda)e^{-2i\lambda^2 t}\,d\lambda\,.
\end{equation}

First, we will show that away from $\lambda=0$, we can get arbitrary
decay in time.  
\begin{lemma}
Fix $L>0$.  For some $k\ge 2$, suppose that $f(\lambda)$ is in
$C^n([L,\infty])$ for all $n=0,1,\dots,k$.  Suppose that $f(L)=f'(L)=\dots 
=f^{(k-1)}(L)=0$ and that the limit
\begin{equation}
\lim_{R\uparrow\infty}\int_L^R f(\lambda)e^{-2i\lambda^2 t}\,d\lambda
\end{equation}
exists for $t\neq 0$.  Then
\begin{equation}
\left|\lim_{R\uparrow\infty}\int_L^R f(\lambda)e^{-2i\lambda^2t}\,d\lambda
\right|\le\frac{1}{L\cdot 4^k\cdot |t|^k}
\sup_{\lambda > L}\left|\lambda^2 (\op{A}^kf)(\lambda)\right|\,,
\end{equation}
where the operator $\op{A}$ is defined by
\begin{equation}
(\op{A}f)(\lambda)\doteq \frac{\partial}{\partial\lambda}\left(\frac{f(\lambda)}{\lambda}\right)\,.
\end{equation}
\end{lemma}
\noindent{\bf Proof:} Integrating by parts $k$ times, 
\begin{equation}
\begin{array}{rcl}
\displaystyle
\lim_{R\uparrow\infty}\int_L^Rf(\lambda)e^{-2i\lambda^2 t}\,d\lambda &=&
\displaystyle
\lim_{R\uparrow\infty}\Bigg[\sum_{n=0}^{k-1}
\left(\frac{i}{4t}\right)^{n+1}(-1)^n\lambda^{-1}e^{-2i\lambda^2t}
(\op{A}^nf)(\lambda)\Bigg|_{\lambda=L}^{\lambda=R} \\\\
&&\displaystyle\hspace{0.2 in}+\,\, \left(
\frac{-i}{4t}\right)^k
\int_L^R (\op{A}^kf)(\lambda)e^{-2i\lambda^2t}\,d\lambda\Bigg]\,.
\end{array}
\end{equation}
The boundary terms at $\lambda=L$ vanish identically, and those at
$\lambda=R$ tend to zero as $R\uparrow\infty$.  These facts prove the
existence of the limit of the integral in the second line, and we find
\begin{equation}
\begin{array}{rcl}
\displaystyle
\left|\lim_{R\uparrow\infty}\int_L^R f(\lambda)e^{-2i\lambda^2 t}\,d\lambda
\right|&=&\displaystyle\frac{1}{4^k\cdot |t|^k}
\left|\lim_{R\uparrow\infty}\int_L^R (\op{A}^kf)(\lambda)e^{-2i\lambda^2t}\,d\lambda\right|\\\\
&\le&\displaystyle\frac{1}{4^k\cdot |t|^k}\lim_{R\uparrow\infty}\int_L^R
\left|\lambda^2 (\op{A}^kf)(\lambda)\right|\cdot\frac{d\lambda}{\lambda^2}\\\\
&\le &\displaystyle
\frac{1}{L\cdot 4^k\cdot |t|^k}\sup_{\lambda>L}\left|\lambda^2 (\op{A}^kf)(\lambda)\right|\,.
\end{array}
\end{equation}
The bound is finite for $k\ge 2$. \pfend

We can now apply this result to estimate $\tilde{h}_\Delta(x,y;t)$.
\begin{lemma}
Fix an integer $k\ge 2$, and let $\sigma > k+1/2$.  Then, the function
$\tilde{h}_\Delta(x,y;t)$ 
is in $L^2(\R^2)$ as a function of $x$ and $y$, with norm
decaying as $|t|^{-k}$.
\label{lemma:tildehdelta}
\end{lemma}
{\bf Proof:}  We apply the above lemma with $L=\Delta/2$ and $f(\lambda)=
\Psi^{\rm (e,o)}_{\rm d}(x,0,\lambda)\Psi^{\rm (e,o)}_{\rm d}(y,0,\lambda)^*
\tilde{g}_{\Delta}(\lambda)$.  This gives the pointwise estimate
\begin{equation}
|\tilde{h}_\Delta(x,y;t)|
\le \frac{2\langle x\rangle^{-\sigma}\langle y\rangle^{-\sigma}}
{\Delta\cdot 4^k\cdot |t|^k}
\sup_{\lambda>\Delta/2}\left|\lambda^2 
(\op{A}^k \Psi^{\rm (e,o)}_{\rm d}(x,0,\cdot)
\Psi^{\rm (e,o)}_{\rm d}(y,0,\cdot)^*\tilde{g}_{\Delta}(\cdot))(\lambda)
\right|\,.
\end{equation}
The operator $\op{A}^k$ acting on the right-hand side makes the
supremum bound a polynomial in $|x|$ and $|y|$ of degree $k$.
Therefore for $\tilde{h}_\Delta(x,y;t)$ to lie in $L^2(\R^2)$ as a function
of $x$ and $y$, it is sufficient to take $\sigma>k+1/2$.  The claimed
time decay of the $L^2$ norm is then obvious.  Note that each
derivative of $\tilde{g}_\Delta(\lambda)$ contributes a factor of order
$O(\Delta^{-1})$, so the overall bound on the $L^2$ norm of $\tilde{h}_\Delta
(x,y;t)$
scales like $\Delta^{-(k+1)}$.  \pfend

Now, we move on to consider the part of $h(x,y;t)$ contributed by
the neighborhood of $\lambda=0$.  We again need some technical lemmas.
\begin{lemma}
For all $\mu\in\R$,
\begin{equation}
\left|\int_0^\mu e^{-2i\zeta^2}\,d\zeta\right|\le\sqrt{3}\,.
\end{equation}
\end{lemma}
\noindent{\bf Proof:}  First, note that
\begin{equation}
\left|\int_0^\mu e^{-2i\zeta^2}\,d\zeta\right|\le\int_0^\mu |d\zeta|=|\mu|\,.
\end{equation}
This estimate is useful for bounded $\mu$.  Suppose $\mu > M > 0$.  Then,
\begin{equation}
\left|\int_0^\mu e^{-2i\zeta^2}\,d\zeta\right|\le M + \left|\int_M^\mu 
e^{-2i\zeta^2}\,d\zeta\right|\,.
\end{equation}
Changing variables to $\tau=\zeta^2$ and integrating by parts, one finds
\begin{equation}
\left|\int_M^\mu e^{-2i\zeta^2}\,d\zeta\right|=\left|\frac{ie^{-2i\mu^2}}{4\mu}
-\frac{ie^{-2iM^2}}{4M}+\frac{i}{8}\int_{M^2}^{\mu^2}e^{-2i\tau}\tau^{-3/2}\,
d\tau\right|\le\frac{3}{4M}\,.
\end{equation}
Therefore, for $\mu>M>0$, we have the estimate
\begin{equation}
\left|\int_0^\mu e^{-2i\zeta^2}\,d\zeta\right|\le M+\frac{3}{4M}\,.
\end{equation}
The right hand side takes its smallest value, $\sqrt{3}$, for
$M_{\rm min}=\sqrt{3}/2$.  Since for $0<\mu\le M_{\rm min}$, we have
\begin{equation}
\left|\int_0^\mu e^{-2i\zeta^2}\,d\zeta\right|\le |\mu|\le M_{\rm min}\le 2M_{\rm min} = \sqrt{3}\,,
\end{equation}
the lemma is established uniformly for all positive $\mu$.  By
symmetry, the same estimate holds for $\mu<0$.  \pfend

\begin{lemma}
\label{lemma:nearzero3/2}
Fix $L>0$ and suppose $f(\lambda)$ is twice continuously differentiable,
with $f(0)=f'(0)=0$, and $f(L)=f'(L)=0$.  Then
\begin{equation}
\left|\int_0^L f(\lambda)e^{-2i\lambda^2t}\,d\lambda\right|\le 
\frac{L\sqrt{3}}{4|t|^{3/2}}\sup_{0<\lambda < L}\left|
\frac{\partial^2}{\partial\lambda^2}\left(\frac{f(\lambda)}{\lambda}\right)
\right|\,.
\end{equation}
\end{lemma}
\noindent{\bf Proof:} Integrating by parts using the boundary
conditions (evaluations at the lower boundary of $\lambda=0$ are
interpreted in the sense of the limit $\lambda\downarrow 0$, that is,
from above), we have
\begin{equation}
\int_0^L f(\lambda)e^{-2i\lambda^2 t}\,d\lambda=
\frac{i}{4t}\int_0^L\frac{f(\lambda)}{\lambda}
\frac{\partial}{\partial\lambda}\left(e^{-2i\lambda^2 t}\right)\,d\lambda=
-\frac{i}{4t}\int_0^L\frac{\partial}{\partial\lambda}\left(\frac{f(\lambda)}
{\lambda}\right)e^{-2i\lambda^2 t}\,d\lambda\,.
\end{equation}
Write 
\begin{equation}
e^{-2i\lambda^2 t}=\frac{\partial}{\partial\lambda}
\int_0^\lambda e^{-2i\sigma^2 t}\,d\sigma\,,
\end{equation}
and integrate by parts again making use of the boundary conditions
(with the same caveat as above), to find
\begin{equation}
\int_0^L f(\lambda)e^{-2i\lambda^2 t}\,d\lambda=
\frac{i}{4t}\int_0^L\frac{\partial^2}{\partial\lambda^2}\left(\frac{f(\lambda)}
{\lambda}\right)\int_0^\lambda e^{-2i\sigma^2 t}\,d\sigma\,d\lambda\,.
\end{equation}
With a change of variables to $\zeta=|t|^{1/2}\sigma$, this becomes
\begin{equation}
\int_0^L f(\lambda)e^{-2i\lambda^2 t}\,d\lambda=
\frac{i}{4t\cdot|t|^{1/2}}\int_0^L\frac{\partial^2}{\partial\lambda^2}
\left(\frac{f(\lambda)}{\lambda}\right)\int_0^{|t|^{1/2}\lambda}e^{-2i\zeta^2}\,d\zeta\,d\lambda\,.
\end{equation}
Estimating the $\lambda$ integral in the obvious way using the uniform
bound of the $\zeta$ integral by $\sqrt{3}$ establishes the claimed
estimate.  \pfend

Without the vanishing boundary conditions at $\lambda=0$, one gets less decay
in time.
\begin{lemma}
\label{lemma:nearzero1/2}
Let $f(\lambda)$ be absolutely continuous $0\le\lambda\le L$, so that
$f'(\lambda)\in L^1([0,L])$.  Then
\begin{equation}
\left|\int_0^L f(\lambda)e^{-2i\lambda^2 t}\,d\lambda\right|\le
\left(|f(0)| + 2\int_0^L |f'(\lambda)|\,d\lambda\right)\frac{\sqrt{3}}{|t|^{1/2}}\,,
\end{equation}
an order $O(|t|^{-1/2})$ bound.
\end{lemma}
\noindent{\bf Proof:}  Separate off the slow decay by writing
\begin{equation}
\int_0^Lf(\lambda)e^{-2i\lambda^2 t}\,d\lambda=f(0)\int_0^Le^{-2i\lambda^2 t}
\,d\lambda + \int_0^L(f(\lambda)-f(0))e^{-2i\lambda^2t}\,d\lambda=I_A+I_B\,.
\end{equation}
The first integral is easily transformed:
\begin{equation}
I_A=f(0)\int_0^Le^{-2i\lambda^2t}\,d\lambda=\frac{f(0)}{t^{1/2}}\int_0^{Lt^{1/2}}e^{-2i\zeta^2}\,d\zeta\,,
\end{equation}
and therefore easily uniformly estimated 
\begin{equation}
|I_A|\le\frac{\sqrt{3}|f(0)|}{|t|^{1/2}}\,.
\end{equation}
In the second integral, one integrates by parts to find
\begin{equation}
I_B=\int_0^L f'(\lambda)\int_\lambda^L e^{-2i\mu^2 t}\,d\mu\,d\lambda\,.
\end{equation}
Therefore,
\begin{equation}
\begin{array}{rcl}
\displaystyle
|I_B|&\le &\displaystyle
\sup_{0<\lambda <L}\left|\int_\lambda^L e^{-2i\mu^2 t}\,d\mu\right|
\cdot \int_0^L |f'(\lambda)|\,d\lambda\\\\
&\le&\displaystyle
\left(\left|\int_0^L e^{-2i\mu^2 t}\,d\mu\right|+
\sup_{0<\lambda <L}\left|\int_0^\lambda e^{-2i\mu^2 t}\,d\mu\right| 
\right)
\int_0^L |f'(\lambda)|\,d\lambda\\\\
&\le& \displaystyle
\frac{2\sqrt{3}}{|t|^{1/2}}\int_0^L |f'(\lambda)|\,d\lambda\,.
\end{array}
\end{equation}
Combining the estimates for $I_A$ and $I_B$ establishes the claimed
result.  \pfend

We now want to use these results to estimate $h_\Delta(x,y;t)$.  To do
this, we want to apply Lemma~\ref{lemma:nearzero3/2} or
Lemma~\ref{lemma:nearzero1/2} with $f(\lambda)=\Psi^{\rm (e,o)}_{\rm
d}(x,0,\lambda)\Psi^{\rm (e,o)}_{\rm
d}(y,0,\lambda)^*g_\Delta(\lambda)$.  Now, from
Lemma~\ref{lemma:qproperties}, it is clear that the hypotheses of
Lemma~\ref{lemma:nearzero3/2} concerning the behavior of $f$ at
$\lambda=0$ will only be satisfied in the odd case.  Here, we obtain
the following.
\begin{lemma}
\label{lemma:hdeltaodd}
Consider the odd case, and let $\sigma>5/2$.  
Then $h_\Delta(x,y;t)$ is in $L^2(\R^2)$ as a function of $x$ and $y$ with
norm decaying like $|t|^{-3/2}$.
\end{lemma}
\noindent{\bf Proof:}
We have the
pointwise estimate
\begin{equation}
|h_\Delta(x,y;t)|\le \frac{3\sqrt{3}\Delta
\langle x\rangle^{-\sigma}\langle y\rangle^{-\sigma}}
{8|t|^{3/2}}\sup_{0<\lambda<3\Delta/2}\left|\frac{\partial^2}{\partial\lambda^2}\left(\frac{\Psi^{\rm (o)}_{\rm d}
(x,0,\lambda)\Psi^{\rm (o)}_{\rm d}(y,0,\lambda)^*g_\Delta(\lambda)}{\lambda}
\right)\right|\,.
\end{equation}
From Lemma~\ref{lemma:qproperties} we have that the right-hand side is
a quadratic polynomial in $|x|$ and $|y|$.  Therefore for
$h_\Delta(x,y;t)$ to be in $L^2(\R^2)$ as a function of $x$ and $y$ it
is sufficient to take $\sigma>5/2$.  The time decay of the $L^2$ norm
is then obvious.  Note that each derivative of $g_\Delta(\lambda)$
contributes a factor that is $O(\Delta^{-1})$ so the bound on the
$L^2$ norm scales like $\Delta^{-1}$.  \pfend

In the even case, we have a zero-energy resonance, and this means that
the integrand near $\lambda=0$ is not small enough to allow decay as
rapid as in the odd case.  In this case, we can only apply
Lemma~\ref{lemma:nearzero1/2} to find the following.
\begin{lemma}
\label{lemma:hdeltaeven}
Consider the even case, and let $\sigma>3/2$.  Then $h_\Delta(x,y;t)$
is in $L^2(\R^2)$ as a function of $x$ and $y$, with norm decaying
like $|t|^{-1/2}$.
\end{lemma}
\noindent{\bf Proof:} Using $f(\lambda)=\Psi^{\rm (e)}_{\rm
d}(x,0,\lambda)\Psi^{\rm (e)}_{\rm d}(y,0,\lambda)^*g_\Delta(\lambda)$
and $L=3\Delta/2$ in Lemma~\ref{lemma:nearzero1/2},
we have the pointwise estimate
\begin{equation}
|h_\Delta(x,y;t)|\le \frac{\sqrt{3}\langle x\rangle^{-\sigma}\langle
y\rangle^{-\sigma}}{|t|^{1/2}}
\left(|f(0)|+2\int_0^{3\Delta/2}|f'(\lambda)|\,d\lambda\right)\,.
\end{equation}
Since the derivative with respect to $\lambda$ results in at most
linear growth in $x$ and $y$, taking $\sigma>3/2$ is sufficient to
ensure that $h_\Delta(x,y;t)$ is in $L^2(\R^2)$ as a function of $x$
and $y$.  Clearly, for large $t$, the $L^2$ norm is $O(|t|^{-1/2})$.
Note that the estimate is also $O(\Delta^{-1})$ due to differentiation
of the bump function $g_\Delta(\lambda)$.  \pfend

In both odd and even cases, the contribution of $h_\Delta(x,y;t)$ to
the $L^2$ norm of $h(x,y;t)$ dominates for large time that of
$\tilde{h}_\Delta(x,y;t)$, for which we had arbitrary decay.
According to Lemma~\ref{lemma:tildehdelta}, for $\sigma>5/2$ this latter
decay is at least as fast as $|t|^{-2}$.  These results imply the
following.
\begin{lemma}
Fix $\sigma>5/2$.  Then, for $t$ sufficiently large, we have the
estimates:
\begin{equation}
\left\|\langle\cdot\rangle^{-\sigma}\left(e^{-it\op{B}}
\op{P}_{\rm c}^{\rm (e)}f\right)(\cdot)\right\|_2\le\frac{K^{\rm (e)}}{|t|^{1/2}}\|\langle\cdot\rangle^\sigma f(\cdot)\|_2\,,
\end{equation}
and
\begin{equation}
\left\|\langle\cdot\rangle^{-\sigma}\left(e^{-it\op{B}}
\op{P}_{\rm c}^{\rm (o)}f\right)(\cdot)\right\|_2\le\frac{K^{\rm (o)}}{|t|^{3/2}}\|\langle\cdot\rangle^\sigma f(\cdot)\|_2\,,
\end{equation}
where $K^{\rm (e)}$ and $K^{\rm (o)}$ are some positive constants.
\end{lemma}

This result, taken together with the elementary time-independent bound
established in Lemma~\ref{lemma:uniform} completes the proof of
Proposition~\ref{prop:nonsingular}.  

\subsection{Singular local decay.}
Now we prove the {\em singular local decay estimate} for the unitary group
$e^{-it\op{B}}$.
\begin{prop}
\label{prop:singular}
Let $|\mu|\ge\mu_{\rm min}>0$.  Fix $\sigma>7/2$.  Let $t=\kappa r$ with
$r\ge 0$ and $\kappa=\pm 1$.  Then, there exist constants $M^{\rm
(e,o)}>0$, such that
\begin{equation}
\left\|\langle \cdot\rangle^{-\sigma}\lim_{\delta\downarrow 0}(
(\op{B}-2\mu-2i\kappa\delta)^{-1}e^{-it\op{B}}\op{P}_{\rm c}^{\rm
(e)}f)(\cdot)\right\|_2\le M^{\rm (e)}\langle
r\rangle^{-1/2}\|\langle\cdot\rangle^\sigma f(\cdot)\|_2\,,
\end{equation}
and
\begin{equation}
\left\|\langle \cdot\rangle^{-\sigma}\lim_{\delta\downarrow 0}(
(\op{B}-2\mu-2i\kappa\delta)^{-1}e^{-it\op{B}}\op{P}_{\rm c}^{\rm
(o)}f)(\cdot)\right\|_2\le M^{\rm (o)}\langle
r\rangle^{-3/2}\|\langle\cdot\rangle^\sigma f(\cdot)\|_2\,.
\end{equation}
The constants $M^{\rm (e,o)}$ depend only on $\mu_{\rm min}$, so the bounds
are uniform for large $|\mu|$.
\end{prop}

The proof of this proposition begins with a representation similar to
(\ref{eq:representation}),
\begin{equation}
\langle x\rangle^{-\sigma}\lim_{\delta\downarrow 0}
\left((\op{B}-2\mu-2i\kappa
\delta)^{-1}e^{-it\op{B}}\op{P}_{\rm c}^{\rm (e,o)}f
\right)(x)=\int_{-\infty}^\infty \langle y\rangle^\sigma f(y) k(x,y;t)\,dy\,,
\end{equation}
where
\begin{equation}
k(x,y;t)=\langle x\rangle^{-\sigma}\langle
y\rangle^{-\sigma}\lim_{\delta\downarrow 0}\int_0^\infty
\frac{\Psi_{\rm d}^{\rm (e,o)}(x,0,\lambda)\Psi_{\rm d}^{\rm (e,o)}(y,0,\lambda)^*}{2\lambda^2 - 2\mu-2i\kappa\delta}e^{-2i\lambda^2t}\,d\lambda\,.
\end{equation}
Perhaps despite appearances, the kernel $k(x,y;t)$ is somewhat more
amenable to analysis than the kernel $h(x,y;t)$ that appeared in the
nonsingular case.  This is because for each finite $\delta$ the
integrand is absolutely integrable as a consequence of the uniform
boundedness in $\lambda$ of $\Psi_{\rm d}^{\rm (e,o)}(x,0,\lambda)$ as
guaranteed by Lemma~\ref{lemma:qproperties} and the large $\lambda$
behavior of the denominator.

Again, the goal is to show that the kernel $k(x,y;t)$ is in
$L^2(\R^2)$ as a function of $x$ and $y$, with norm that is decaying
in time, although in this case we will only obtain the decay for $t$
of a particular sign.  First, we show that the $L^2$ norm exists and
is finite near $t=0$.
\begin{lemma}
\label{lemma:singularsmallt}
Fix $\sigma>3/2$, $\mu$ with $|\mu|\ge\mu_{\rm min}>0$, and $t$ with
$|t|<T$.  Then there exist constants $C^{\rm (e,o)}>0$ depending on
$\mu_{\rm min}$ and $T$ such that
\begin{equation}
\|k(\cdot,\cdot;t)\|_2\le C^{\rm (e,o)}\,.
\end{equation}
Since the bounds only depend on $\mu$ via $\mu_{\rm min}$, they are uniform
for large $|\mu|$.
\end{lemma}
\noindent{\bf Proof:} Begin by setting $f(\lambda)=\Psi_{\rm d}^{\rm
(e,o)}(x,0,\lambda)\Psi_{\rm d}^{\rm
(e,o)}(y,0,\lambda)^*e^{-2i\lambda^2t}$.
First we consider $\mu\le-\mu_{\rm min}<0$, in which case we have
\begin{equation}
k(x,y;t)=\frac{\langle x\rangle^{-\sigma}\langle y\rangle^{-\sigma}}{2}\int_0^\infty\frac{f(\lambda)\,d\lambda}{\lambda^2-\mu}\,,
\end{equation}
since there is no singularity for $\mu<0$.  We immediately get the
pointwise estimate
\begin{equation}
\begin{array}{rcl}
|k(x,y;t)|&\le&\displaystyle
\frac{\langle x\rangle^{-\sigma}\langle y\rangle^{-\sigma}}{2}\sup_{\lambda>0}|f(\lambda)|\int_0^\infty \frac{d\lambda}{\lambda^2-\mu}\\\\
&=&\displaystyle
\frac{\pi\langle x\rangle^{-\sigma}\langle y\rangle^{-\sigma}}{4\sqrt{-\mu}}\sup_{\lambda>0}|f(\lambda)|\\\\
&\le & \displaystyle
\frac{\pi\langle x\rangle^{-\sigma}\langle y\rangle^{-\sigma}}{4\sqrt{-\mu_{\rm min}}}\sup_{\lambda>0}|f(\lambda)|
\,.
\end{array}
\end{equation}

Now we consider $\mu\ge\mu_{\rm min}>0$.  Pick some positive $G$ less
than $\sqrt{\mu_{\rm min}}$.  Then
\begin{equation}
k(x,y;t)=\frac{\langle x\rangle^{-\sigma}\langle y\rangle^{-\sigma}}{2}
(I_{\rm ns}+I_{\rm s})\,,
\end{equation}
where
\begin{equation}
\begin{array}{rcl}
I_{\rm ns}&\doteq &\displaystyle
\int_0^{\sqrt{\mu}-G}\frac{f(\lambda)\,d\lambda}{\lambda^2-\mu}+
\int_{\sqrt{\mu}+G}^\infty\frac{f(\lambda)\,d\lambda}{\lambda^2-\mu}\,,\\\\
I_{\rm s}&\doteq &\displaystyle
\lim_{\delta\downarrow 0}\int_{\sqrt{\mu}-G}^{\sqrt{\mu}+G}
\frac{f(\lambda)\,d\lambda}
{(\lambda+\sqrt{\mu+i\kappa\delta})(\lambda-\sqrt{\mu+i\kappa\delta})}
\,,
\end{array}
\end{equation}
where the principal branch of the square root is understood, so that
the square root is nearly a positive number for $\delta$ small.  It is easy
to find
\begin{equation}
\begin{array}{rcl}
|I_{\rm ns}|&\le &\displaystyle
\sup_{\lambda>0}|f(\lambda)|\left(\int_0^{\sqrt{\mu}-G}\frac{d\lambda}{\mu-\lambda^2}+\int_{\sqrt{\mu}+G}^\infty\frac{d\lambda}{\lambda^2-\mu}\right)\\\\
&\le &\displaystyle
\frac{\sup_{\lambda>0}|f(\lambda)|}{\sqrt{\mu}}
\left({\rm arctanh}\left(\frac{\sqrt{\mu}}{\sqrt{\mu}+G}\right)+{\rm arctanh}\left(\frac{\sqrt{\mu}-G}{\sqrt{\mu}}\right)\right)\\\\
&\le &\displaystyle
\sup_{\lambda>0}|f(\lambda)|\sup_{\mu>\mu_{\rm min}}\left
(\frac{1}{\sqrt{\mu}}\left({\rm arctanh}\left(\frac{\sqrt{\mu}}{\sqrt{\mu}+G}\right)+{\rm arctanh}\left(\frac{\sqrt{\mu}-G}{\sqrt{\mu}}\right)\right)
\right)\,.
\end{array}
\end{equation}

For the singular part, we find
\begin{equation}
\begin{array}{rcl}
I_{\rm s}&=&\displaystyle
\lim_{\delta\downarrow 0}\frac{f(\sqrt{\mu})}{\sqrt{\mu}+\sqrt{\mu+i\kappa\delta}}\int_{\sqrt{\mu}-G}^{\sqrt{\mu}+G}\frac{d\lambda}{\lambda-\sqrt{\mu+i\kappa\delta}} +
\int_{\sqrt{\mu}-G}^{\sqrt{\mu}+G}\left(\frac{f(\lambda)}{\lambda+\sqrt{\mu}}-\frac{f(\sqrt{\mu})}{2\sqrt{\mu}}
\right)\frac{d\lambda}{\lambda-\sqrt{\mu}}\\\\
&=&\displaystyle
\frac{i\pi\kappa f(\sqrt{\mu})}{2\sqrt{\mu}} +\int_{\sqrt{\mu}-G}^{\sqrt{\mu}+G}\left(\frac{f(\lambda)}{\lambda+\sqrt{\mu}}-\frac{f(\sqrt{\mu})}{2\sqrt{\mu}}\right)\frac{d\lambda}{\lambda-\sqrt{\mu}}\,.
\end{array}
\end{equation}
Therefore, 
\begin{equation}
\begin{array}{rcl}
|I_{\rm s}|&\le &\displaystyle
\frac{\pi|f(\sqrt{\mu})|}{2\sqrt{\mu}} + 2G\sup_{|\lambda-\sqrt{\mu}|<G}\left|
\left(\frac{f(\lambda)}{\lambda+\sqrt{\mu}}-\frac{f(\sqrt{\mu})}{2\sqrt{\mu}}\right)\frac{1}{\lambda-\sqrt{\mu}}\right|\\\\
&\le &\displaystyle
\frac{\pi|f(\sqrt{\mu})|}{2\sqrt{\mu}} + 2G\sup_{|\lambda-\sqrt{\mu}|<G}\left|
\frac{\partial}{\partial\lambda}\left(\frac{f(\lambda)}{\lambda+\sqrt{\mu}}\right)
\right|\,.
\end{array}
\end{equation}
Again, the bounds are uniform in $\mu$ for large $\mu$.

Now, we simply note that the pointwise bounds for $I_{\rm ns}$ and
$I_{\rm s}$ are themselves bounded by functions of $x$, $y$, and $t$
that grow linearly at worst, as a consequence of differentiation of
$f(\lambda)$ with respect to $\lambda$
(c.f. Lemma~\ref{lemma:qproperties} for the growth in $x$ and $y$,
while the growth in $t$ comes from the factor $e^{-2i\lambda^2 t}$).
Thus, to have $k(x,y;t)\in L^2(\R^2)$ as a function of $x$ and $y$, it
is sufficient to take $\sigma > 3/2$, and then the norm will be
bounded by a linear function of $|t|$, and therefore uniformly for
$|t|<T$.  The bound is also uniform in $\mu$ for $|\mu|\ge\mu_{\rm
min}>0$.  \pfend

We now note that proving the decay for large $r=|t|$ for
$\mu\le-\mu_{\rm min}<0$ amounts to recalling the nonsingular local
decay estimate.  The integral is not really singular: 
\begin{equation}
\lim_{\delta\downarrow 0}\int_0^\infty\frac{\Psi_{\rm d}^{\rm (e,o)}(x,0,\lambda)\Psi_{\rm d}^{\rm (e,o)}(y,0,\lambda)}{\lambda^2-\mu-i\kappa\delta}
e^{-2i\lambda^2t}\,d\lambda=\int_0^\infty
\frac{\Psi_{\rm d}^{\rm (e,o)}(x,0,\lambda)\Psi_{\rm d}^{\rm (e,o)}(y,0,\lambda)}{\lambda^2-\mu}
e^{-2i\lambda^2t}\,d\lambda\,.
\end{equation}
Using the same arguments as used to prove the nonsingular local decay
estimate one gets a pointwise bound for this integral that is at most
quadratically growing in $x$ and $y$ and decaying like $|t|^{-1/2}$ in
the even case and $|t|^{-3/2}$ in the odd case.  Since the estimates
involve up to two derivatives of the quotient in the integrand, the
bounds will be uniform in $\mu$ for $\mu\le -\mu_{\rm min}<0$.

To prove the decay for large $r=|t|$ in the truly singular case when
$\mu\ge\mu_{\rm min}>0$, we split $k(x,y;t)$ into three parts.  Let
$g_\Delta(\lambda)$ and $\tilde{g}_\Delta(\lambda)$ be as before, and
introduce the new ``bump'' functions $g_G(\lambda)$ and
$\tilde{g}_G(\lambda)=1-g_G(\lambda)$, both infinitely differentiable
and nonnegative, with $g_G(\lambda)$ identically equal to zero outside
of the interval $(\sqrt{\mu}-3G/2,\sqrt{\mu}+3G/2)$ and identically
equal to one inside of the interval $(\sqrt{\mu}-G/2,\sqrt{\mu}+G/2)$.
Set
\begin{equation}
k(x,y;t)=\frac{\langle x\rangle^{-\sigma}\langle y\rangle^{-\sigma}}{2}
\left(I_0+I_\mu+\tilde{I}\right)\,,
\end{equation}
where
\begin{equation}
I_0\doteq \int_0^{3\Delta/2}\frac{\Psi_{\rm d}^{\rm (e,o)}(x,0,\lambda)
\Psi_{\rm d}^{\rm (e,o)}(y,0,\lambda)^*}{\lambda^2-\mu}g_{\Delta}(\lambda)
e^{-2i\lambda^2t}\,d
\lambda\,,
\end{equation}
\begin{equation}
I_\mu\doteq \lim_{\delta\downarrow
0}\int_{\sqrt{\mu}-3G/2}^{\sqrt{\mu}+3G/2}\frac{ \Psi_{\rm d}^{\rm
(e,o)}(x,0,\lambda)\Psi_{\rm d}^{\rm (e,o)}(y,0,\lambda)^*}{
\lambda^2-\mu-i\kappa\delta}g_G(\lambda)e^{-2i\lambda^2t}\,d\lambda\,,
\end{equation}
and
\begin{equation}
\begin{array}{rcl}
\tilde{I}&\doteq &\displaystyle
\int_{\Delta/2}^{\sqrt{\mu}-G/2}\frac{\Psi_{\rm d}^{\rm (e,o)}(x,0,\lambda)
\Psi_{\rm d}^{\rm (e,o)}(y,0,\lambda)^*}{\lambda^2-\mu}
\tilde{g}_\Delta(\lambda)\tilde{g}_G(\lambda)e^{-2i\lambda^2t}\,d\lambda \\\\
&&\displaystyle\hspace{0.2 in}+\,\,
\int_{\sqrt{\mu}+G/2}^\infty
\frac{\Psi_{\rm d}^{\rm (e,o)}(x,0,\lambda)
\Psi_{\rm d}^{\rm (e,o)}(y,0,\lambda)^*}{\lambda^2-\mu}
\tilde{g}_G(\lambda)e^{-2i\lambda^2t}\,d\lambda\,.
\end{array}
\end{equation}
Note that in keeping the contributions near zero and near $\mu$
distinct, we are assuming that $3\Delta/2 < \sqrt{\mu_{\rm min}}-G$.
The analysis of $I_0$ and $\tilde{I}$ proceeds exactly as in the proof
of the nonsingular local decay estimate.  The results are almost
identical.  For $\tilde{I}$ one can integrate by parts as many times
as one likes, and therefore one gets a pointwise estimate with
arbitrary decay in time of order $O(|t|^{-k})$ for $k\ge 2$, but at
the cost of polynomial growth in $x$ and $y$ of degree $k$.  For
$I_0$, one gets a pointwise estimate that decays like $O(|t|^{-1/2})$
in the even case and $O(|t|^{-3/2})$ in the odd case, at the cost of
quadratic growth in $x$ and $y$.

The estimates of $I_0$ and $\tilde{I}$ are uniform for large $\mu$.
The pointwise bounds for $\tilde{I}$ involve supremum bounds over the range
of integration of the quantity
\begin{equation}
\lambda^2 \left(\op{A}^k \frac{f(\cdot)}{(\cdot)^2-\mu}\right)(\lambda)\,,
\end{equation}
with $f(\lambda)$ given by $\Psi_{\rm d}^{\rm
(e,o)}(x,0,\lambda)\Psi_{\rm d}^{\rm
(e,o)}(y,0,\lambda)^*\tilde{g}_G(\lambda)$ for
$\lambda>\sqrt{\mu}+G/2$ and with $f(\lambda)$ given by $\Psi_{\rm
d}^{\rm (e,o)}(x,0,\lambda)\Psi_{\rm d}^{\rm
(e,o)}(y,0,\lambda)^*\tilde{g}_\Delta(\lambda)\tilde{g}_G(\lambda)$
for $\Delta/2<\lambda<\sqrt{\mu}-G/2$.  In particular, we will need
these bounds for $k=2$, in which case
\begin{equation}
\lambda^2\left(\op{A}^k \frac{f(\cdot)}{(\cdot)^2-\mu}\right)(\lambda)=
\left(\frac{15\lambda^4-10\mu\lambda^2+3\mu}{\lambda^2(\lambda^2-\mu)^3}
\right)f(\lambda)+\left(\frac{-7\lambda^2+3\mu}{\lambda(\lambda^2-\mu)^2}
\right)f'(\lambda)+\frac{f''(\lambda)}{\lambda^2-\mu}\,.
\end{equation}
For the part of $\tilde{I}$ involving $\lambda>\sqrt{\mu}+G/2$, it is easy to
check that the three coefficients above are monotonic functions of
$\lambda$ for $\lambda>\sqrt{\mu}$ that decay for large $\lambda$ with $\mu$
fixed.  Therefore each coefficient is bounded by its magnitude at the
lower endpoint $\lambda=\sqrt{\mu}+G/2$.  With $G$ held fixed, these bounds
are then seen to be decaying functions of $\mu$.  For the part of
$\tilde{I}$ involving $\lambda\in(\Delta/2,\sqrt{\mu}-G/2)$, it is easy to
see that the coefficients blow up at both endpoints. Therefore, for
$\Delta$ and $G$ sufficiently small but independent of $\mu$, the
coefficients will be bounded by the maximum of their values at
$\lambda=\Delta/2$ and $\lambda=\sqrt{\mu}-G/2$.  Again, holding $\Delta$ and
$G$ fixed, one sees that the bounds are uniform for large $\mu$.  This
direct argument shows that, at least for $k=2$, the pointwise bound
for $\tilde{I}$ is uniform in $\mu$.  Establishing the uniformity of
the pointwise estimate for $I_0$ is easier; the denominator
$\lambda^2-\mu$ plays no essential role for $\lambda<3\Delta/2$ for
$\mu$ sufficiently large.

The new term that must be handled
differently is $I_\mu$.
\begin{lemma}
For all $k\ge 2$, the integral $I_\mu$ satisfies the pointwise
estimate 
\begin{equation}
|I_\mu|\le \frac{2}{k-1}\frac{P_k(x,y)}{|t|^{k-1}}\,,
\end{equation}
where $P_k(x,y)$ is a polynomial in $|x|$ and $|y|$ of degree $k$ with
positive coefficients that are uniform in $\mu$.
\end{lemma}
\noindent{\bf Proof:}
Consider first $t>0$.  Then, the quantity to estimate is
\begin{equation}
\begin{array}{rcl}
I_\mu&=&\displaystyle
\lim_{\delta\downarrow 0}
\int_{\sqrt{\mu}-3G/2}^{\sqrt{\mu}+3G/2}
\frac{\Psi_{\rm d}^{\rm (e,o)}(x,0,\lambda)
\Psi_{\rm d}^{\rm (e,o)}(y,0,\lambda)^*}
{\lambda^2-\mu-i\delta}g_G(\lambda)e^{-2i\lambda^2t}\,d\lambda\\\\
&=&\displaystyle
2i
e^{-2i\mu t}\lim_{\delta\downarrow 0} e^{2\delta t}
\int_{\sqrt{\mu}-3G/2}^{\sqrt{\mu}+3G/2}
\Psi_{\rm d}^{\rm (e,o)}(x,0,\lambda)
\Psi_{\rm d}^{\rm (e,o)}(y,0,\lambda)^*
g_G(\lambda)\int_t^\infty 
e^{-2i(\lambda^2-\mu-i\delta)s}\,ds\,d\lambda\\\\
&=&\displaystyle
2ie^{-2i\mu t}\int_t^\infty e^{2i\mu s}
\left[\int_{\sqrt{\mu}-3G/2}^{\sqrt{\mu}+3G/2}
\Psi_{\rm d}^{\rm (e,o)}(x,0,\lambda)\Psi_{\rm d}^{\rm (e,o)}(y,0,\lambda)^*
g_G(\lambda)e^{-2i\lambda^2s}\,d\lambda\right]\,ds\,.
\end{array}
\end{equation}
Now, with $g_G(\lambda)$ vanishing to all orders at the integration
endpoints, it is possible to bound the integral in square brackets by
iterated integration by parts.  The bound is $O(|s|^{-k})$ and grows
in $x$ and $y$ like a polynomial $P_k(x,y)$ of degree $k$.  This bound
is uniform in $\mu$, since the only place $\mu$ appears is in the
range of integration over which bounds are required, and from
Lemma~\ref{lemma:qproperties} we know that these bounds are uniform for all
$\lambda$.  Therefore we have
\begin{equation}
|I_\mu|\le 2\int_t^\infty P_k(x,y)|s|^{-k}\,ds=\frac{2}{k-1}\frac{P_k(x,y)}{|t|^{k-1}}\,,
\end{equation}
which establishes the lemma for $t>0$.  For $t<0$, one gets an
integral from $-\infty$ to $t$ in the second step above, and
ultimately obtains the same bound.  \pfend

Finally, we put the pieces together to complete the proof of 
Proposition~\ref{prop:singular}.
For the odd case, we want decay
of order $O(|t|^{-3/2})$.  For $\langle x\rangle^{-\sigma}\langle
y\rangle^{-\sigma}I_0$ to be in $L^2(\R^2)$ with this decay rate, we
need $\sigma>5/2$.  With this bound on $\sigma$, we can get $\langle
x\rangle^{-\sigma}\langle y\rangle^{-\sigma}\tilde{I}$ being in
$L^2(\R^2)$ with decay bounded by $O(|t|^{-2})=o(|t|^{-3/2})$, but no
better.  Finally, for $\langle x\rangle^{-\sigma}\langle
y\rangle^{-\sigma} I_\mu$ to be in $L^2(\R^2)$ with decay
$O(|t|^{-2})$ we now see that we need to localize a bit more in space
by taking $\sigma>7/2$.  Combining these large time estimates with the
finite time bound of Lemma~\ref{lemma:singularsmallt} establishes the
proposition in the odd case.  Similar arguments for the even case give
an $L^2(\R^2)$ norm that decays like $O(|t|^{-1/2})$ for $\sigma>7/2$.
This finishes the proof of Proposition~\ref{prop:singular}.

\noindent{\bf Remark:}  Evidently, the singular decay estimates blow up
when $\mu_{\rm min}$ approaches zero.  This is an essential phenemenon in
both the odd and even cases.  This is best seen by considering the singular
integral for the case $\mu=0$:
\begin{equation}
\int_0^\infty \frac{\Psi_{\rm d}^{\rm (e,o)}(x,0,\lambda)
\Psi_{\rm d}^{\rm (e,o)}(y,0,\lambda)^*}{\lambda^2-i\kappa\delta}
e^{-2i\lambda^2t}\,d\lambda\,.
\end{equation}
This integral blows up for all $x$, $y$, and $t$, as $\delta$ tends to
zero in the even case.  In the odd case there is sufficient vanishing
at $\lambda=0$ for the limit of $\delta\downarrow 0$ to exist for all
$x$, $y$, and $t$, but the limit only decays in $t$ like $|t|^{-1/2}$.
Thus it is not possible for estimates of the form derived for
$|\mu|\ge\mu_{\rm min}>0$ to hold uniformly in any neighborhood of
$\mu=0$.\pfend

\end{document}